\documentclass[%
 aip,
 amsmath,amssymb,
 nofootinbib,
 reprint,%
]{revtex4-1}

\usepackage{graphicx}
\usepackage{dcolumn}
\usepackage{bm}

\usepackage[utf8]{inputenc}
\usepackage[T1]{fontenc}
\usepackage{mathptmx}
\usepackage{etoolbox}

\usepackage{hyperref}
\usepackage{amssymb}
\usepackage{amsmath}
\usepackage{multirow}
\usepackage{cleveref}
\usepackage{booktabs}

\usepackage{xcolor}
\usepackage{xargs}
\usepackage[colorinlistoftodos,textsize=footnotesize]{todonotes}

\DeclareSymbolFont{newfont}{OML}{cmm}{m}{it}
\DeclareMathSymbol{\Epsilon}{3}{newfont}{15}
\DeclareMathSymbol{\Varrho}{3}{newfont}{37}

\newcommandx{\unsure}[2][1=]{\todo[linecolor=red,backgroundcolor=red!25,bordercolor=red,#1]{#2}}
\newcommandx{\change}[2][1=]{\todo[linecolor=blue,backgroundcolor=blue!25,bordercolor=blue,#1]{#2}}
\newcommandx{\info}[2][1=]{\todo[linecolor=OliveGreen,backgroundcolor=OliveGreen!25,bordercolor=OliveGreen,#1]{#2}}

\newcommand{\ten}[1]{\ensuremath{\mathbf{#1}}}

\newcommand{\nam}[1]{\texttt{{#1}}}

\usetikzlibrary{arrows.meta, patterns}
\usetikzlibrary{shapes, calc,quotes,angles}

\makeatletter
\def\@email#1#2{%
 \endgroup
 \patchcmd{\titleblock@produce}
  {\frontmatter@RRAPformat}
  {\frontmatter@RRAPformat{\produce@RRAP{*#1\href{mailto:#2}{#2}}}\frontmatter@RRAPformat}
  {}{}
}%
\makeatother
\begin{document}

\preprint{AIP/123-QED}

\title[Techniques for second order convergent weakly-compressible
smoothed particle hydrodynamics schemes without boundaries]{Techniques for second order convergent weakly-compressible
smoothed particle hydrodynamics schemes without boundaries}
\author{Pawan Negi}
 \email{pawan.n@aero.iitb.ac.in}
\author{Prabhu Ramachandran}%
 \email{prabhu@aero.iitb.ac.in}
\affiliation{Department of Aerospace Engineering, Indian Institute of
Technology Bombay, Powai, Mumbai 400076
}%

\date{\today}

\begin{abstract}
  \begin{sloppypar} \noindent Despite the many advances in the use of weakly-compressible smoothed
    particle hydrodynamics (SPH) for the simulation of incompressible fluid
    flow, it is still challenging to obtain second-order convergence even for
    simple periodic domains. In this paper we perform a systematic numerical
    study of convergence and accuracy of kernel-based approximation,
    discretization operators, and weakly-compressible SPH (WCSPH) schemes. We
    explore the origins of the errors and issues preventing second-order
    convergence despite having a periodic domain. Based on the study, we propose
    several new variations of the basic WCSPH scheme that are all second-order
    accurate. Additionally, we investigate the linear and angular momentum
    conservation property of the WCSPH schemes. Our results show that one may
    construct accurate WCSPH schemes that demonstrate second-order convergence
    through a judicious choice of kernel, smoothing length, and discretization
    operators in the discretization of the governing equations.  \end{sloppypar}
\end{abstract}

\maketitle

\section{Introduction}
\label{sec:intro}

Smoothed Particle Hydrodynamics has been used to simulate weakly-compressible
fluids since the pioneering work of \citet{sph:fsf:monaghan-jcp94}. Many
variations of the basic method have been proposed to create an entire class of
weakly-compressible SPH schemes (WCSPH). One particularly difficult challenge
has been the poor convergence displayed by the WCSPH methods making it one of
the SPH grand-challenge problems~\cite{vacondio_grand_2020}.

The SPH method works by using a smoothing kernel to approximate a function
wherein the choice of the kernel influences the accuracy of the method. The
length scale of the smoothing kernel is often termed the support radius or
smoothing length, $h$. A variety of kernels are used in the literature and the
smoothing length may either be fixed in space/time or varying. One can show
that for a symmetric kernel, the SPH kernel approximation is spatially second-order accurate in $h$. However, the particle discretization of this
approximation seldom achieves this and even first-order convergence often
requires care and tuning of the smoothing length. \citet{hernquist1989treesph}
proposed that the support radius, $h$ be increased such that $h \propto \Delta
s^{-1/3}$ in three dimensions where $\Delta s$ is the local inter-particle
separation. Subsequently, \citet{quinlan_truncation_2006} derived error
estimates for the standard SPH discretization and found that the ratio
$h/\Delta s$ must increase as the $h$ value is reduced to attain convergence;
this is because of error terms of the form $\left( \frac{\Delta s}{h}
\right)^{\beta + 2}$, where $\beta$ is a measure of the smoothness of the
kernel at the edge of its support. This is an issue because as $h$ increases, the number of neighbors for
each particle increases resulting in a prohibitive increase of computational
effort. Furthermore, increasing the smoothing radius also reduces the accuracy
of the method. This is the approach used in the work of
\citet{zhu2015numerical} who proposed that the number of neighbors $N_{nb}
\propto N^{0.5}$, where $N$ is the number of particles, in order to get
convergence using SPH kernels.

\citet{kiara_sph_2013,kiara_sph_2013-1} shows that when the particles are
distributed uniformly it is possible to obtain second-order convergence. The
results of \onlinecite{quinlan_truncation_2006} show that when using sufficiently
smooth kernels (where $\beta$ is large or infinite), one can obtain second-order convergence. Indeed, \citet{lind_high-order_2016} demonstrate that for
particle distributions on a Cartesian mesh one can obtain higher order
convergence using higher order kernels.

However, for kernels that are normally used in SPH, the SPH approximations of
derivatives become inaccurate even on a uniform grid unless a very large
smoothing radius is used. Many methods have been proposed to correct the
gradient approximation
\cite{bonet_lok:cmame:1999,liu_restoring_2006,rosswog2015boosting,huang2019kernel}.
These typically ensure that the derivative approximation of a linear function is
exact. This linear consistency is achieved by inverting a small matrix for each
particle and using this to correct the computed gradients. This makes the
derivative approximation second-order accurate but increases the computational
cost of the gradient computation two-fold.

In the context of incompressible fluid flows, the governing equations involve
the divergence, gradient, and Laplacian operators. These operators must be
discretized and used in the context of Lagrangian particles.  The divergence
operator is encountered in the continuity equation and the discretization
proposed by \onlinecite{sph:fsf:monaghan-jcp94} is widely used.  Rather than using a
continuity equation, some authors~\cite{Adami2013} prefer to use the summation
density formulation proposed in \onlinecite{monaghan1983shock} to directly evaluate
density. The gradient operator is encountered in the momentum equation. Many
authors prefer using a discretized form that manifestly preserves linear
momentum and as a result employ the symmetric form of the gradient operator
\cite{sph:fsf:monaghan-jcp94,bonet_lok:cmame:1999}.  The symmetrization can be
done in two different ways and \citet{violeau2012fluid} shows that the selection
of one form dictates the form to be used for divergence discretization in order
to conserve volume (energy) in phase space.

Given the inaccuracy of the SPH approximation in computing derivatives
accurately, the kernel corrections of
\onlinecite{bonet_lok:cmame:1999,liu_restoring_2006} maybe applied to obtain linear
consistency of the gradient and divergence operators. Unfortunately, the use
of the corrections implies that linear momentum is no longer manifestly
conserved. \citet{crksph:jcp:2017} propose a symmetrization of the corrected
kernel as originally suggested by \citet{dilts1999moving,dilts2000moving} to
conserve linear momentum but the symmetrization implies that the operator is
not first order consistent. Thus, the unfortunate consequence of demanding
linear consistency is lack of conservation and vice-versa. A conservative and
linear consistent gradient operator is currently not available.

The Laplacian is a challenging operator in the context of SPH. The simplest
method is the one where the double derivative of the kernel is employed.
However, the double derivatives of the kernel are very sensitive to any particle
disorder. \citet{chenGeneralizedSmoothedParticle2000a} propose an approach by
considering the inner product with each of the double derivatives and taking
into account the leading order error terms.
\citet{zhangModifiedSmoothedParticle2004} propose using the inner product with
all the derivatives of the kernel lower and equal to the required derivative.
This generates a system of 10 equations in two-dimensions.
\citet{korziliusImprovedCSPMApproach2017} propose an improvement over the method
of \onlinecite{chenGeneralizedSmoothedParticle2000a} to evaluate the correction term.
All of these methods require the computation of higher order kernel derivatives.
Many
authors~\cite{schwaigerImplicitCorrectedSPH2008,maciaBoundaryIntegralSPH2012}
proposed methods to correct the Laplacian near the boundary. In all of these
formulations linear momentum is not manifestly conserved.

The Laplacian may also be discretized using the first derivative of the kernel
using an integral approximation of the Laplacian. This was first suggested by
\citet{brookshaw_method_1985} and has been improved by
\citet{morris-lowRe-97,cleary1999conduction}. They employ a finite difference
approximation to evaluate the first order derivative and then convolve this
with the kernel derivative. This formulation was structured such that it
conserves linear momentum. However, these approximations do not converge as
the resolution increases especially in the context of irregular particle
distributions. \citet{fatehi_error_2011} propose an improved formulation by
accounting for the leading error term; this makes the method accurate and
convergent but makes the approximations non-conservative.

Another method to discretize the Laplacian is the repeated use of a first
derivative and this has been used by \citet{bonet_lok:cmame:1999}, and
\citet{nugent_liquid_2000}. The formulation is generally not popular since it
shows high frequency numerical oscillations when the initial condition is
discontinuous. Recently, \citet{biriukovStableAnisotropicHeat2019a} show that
these oscillations can be removed by employing smoothing near the discontinuity.

Various SPH schemes have been proposed that use the above methods for
discretization of the different operators. The simplest of the schemes is the
original weakly compressible SPH (WCSPH)
method~\cite{sph:fsf:monaghan-jcp94,wcsph-state-of-the-art-2010}. This method
is devised such that it conserves linear momentum as well as the Hamiltonian
of the system. However, as the particles move they become highly disorganized
and this significantly reduces the accuracy of the method. Many particle
regularization methods popularly known as particle shifting techniques (PST)
have been proposed which can be incorporated into WCSPH schemes
\cite{acc_stab_xu:jcp:2009,lind2012incompressible,sun_consistent_2019,huang2019kernel}.
These methods ensure that the particles are distributed more uniformly.
Instead of displacing the particles directly, \citet{Adami2013} propose to use
a transport velocity instead of the particle velocity to ensure a uniform
particle distribution. This approach is also framed in the context of
Arbitrary Lagrangian Eulerian (ALE) SPH schemes by \citet{oger_ale_sph_2016}. A
similar approach is used by \citet{sun_consistent_2019} to incorporate the
shifting velocity in the momentum equation with the $\delta$-SPH
scheme~\cite{marrone-deltasph:cmame:2011,antuono-deltasph:cpc:2010}.
\citet{edac-sph:cf:2019} also employ a transport velocity formulation and
additionally propose using the EDAC scheme \cite{Clausen2013} in the context
of SPH, which removes the need for an equation of state (EOS). The resulting
method is accurate but does not converge with an increase of resolution. An
alternative approach to ensure particle homogeneity is the approach of
remeshing proposed by \citet{remeshed_sph:jcp:2002} where the particles are
periodically interpolated into a regular Cartesian mesh. The method can be
accurate but the remeshing can be diffusive and makes the method reliant on a
Cartesian mesh. In a subsequent development,
\citet{hieberImmersedBoundaryMethod2008} employ remeshing but couple it with
an immersed boundary method to deal with complex solid bodies. Recently,
\citet{nasar2019} modify the method introduced by \citet{lind2012incompressible}
to devise an Eulerian WCSPH scheme that also uses ideas from immersed boundary
methods to handle complex geometry.

To summarize the discussion in the context of convergence, some authors
\cite{quinlan_truncation_2006,schwaigerImplicitCorrectedSPH2008,maciaBoundaryIntegralSPH2012}
demonstrate numerical convergence for the derivative and function
approximation. Many
authors~\cite{kiara_sph_2013-1,huang2019kernel,Adami2013,chenGeneralizedSmoothedParticle2000a,zhangModifiedSmoothedParticle2004,sun2019consistent,marrone-deltasph:cmame:2011}
only show convergence in the form of plots that approach an exact solution
with increasing resolution without formally computing the order of
convergence. Some authors demonstrate second order convergence for simpler
problems with a fixed particle configuration like the heat conduction
equation~\cite{cleary1999conduction,schwaigerImplicitCorrectedSPH2008,fatehi_error_2011,korziliusImprovedCSPMApproach2017},
the Poisson equation~\cite{maciaBoundaryIntegralSPH2012}, and the evolution of
an acoustic wave~\cite{crksph:jcp:2017}. Second order convergence has also
been demonstrated for Eulerian SPH methods where the particles are held
fixed~\cite{lind_high-order_2016,nasar2019} or where the particles are
re-meshed~\cite{hieberImmersedBoundaryMethod2008}. Some
authors~\cite{rosswog2015boosting,crksph:jcp:2017,acc_stab_xu:jcp:2009,diff_smoothing_sph:lind:jcp:2009,edac-sph:cf:2019,dehnen-aly-paring-instability-mnras-2012}
show first order convergence for Lagrangian SPH schemes but this does not
persist as the resolution is increased. Therefore, to the best of our
knowledge, none of the contemporary Lagrangian SPH schemes appear to
demonstrate a formal second order convergence for simple fluid mechanics
problems like the Taylor-Green vortex problem for which an exact solution is
known.

In this paper, we carefully construct a family of Lagrangian SPH schemes that
demonstrate second order numerical convergence for the classic Taylor-Green
vortex problem. We first study several commonly used SPH kernels in the
context of function and derivative approximation using particles that are
either in a Cartesian arrangement or in an irregular but packed configuration
of particles encountered when employing some form of a particle shifting
technique. We choose a suitable correction scheme that produces second order
approximations. We then select a suitable kernel and smoothing radius based on
this study. We then systematically study the various discretization operators
along with suitable corrections. Our investigations are in two-dimensions
although the results are applicable in three dimensions as well. Our numerical
investigation covers a wide range of resolutions with our highest resolution
using a quarter million particles with $\frac{L}{\Delta s}=500$, where $L=1m$
is the length of the domain. Once we have identified suitable second order
convergent operators we carefully construct SPH schemes that display a second
order convergence (SOC). We use the Taylor-Green vortex problem to demonstrate
this. We also compare our results with those of several established SPH
methods that are currently used. We study the accuracy, convergence, and also
investigate the computational effort required. We construct both Lagrangian
and Eulerian schemes that are fully second order convergent. We provide
schemes that use either an artificial compressibility in the form of an
equation state or using a pressure evolution equation.

Once we have demonstrated second order convergence for the Taylor-Green vortex
problem we proceed to investigate the Gresho-Chan vortex~\cite{gresho1990theory}
problem as well as an incompressible shear layer
problem~\cite{diMovingMeshFinite2005} and look at how the lack of manifest
conservation impacts the conservation of linear and angular momentum.  In the
interest of reproducibility, all the results shown in the paper are
automatically generated through the use of an automation
framework~\cite{pr:automan:2018}, and the source code for the paper is available
at \url{https://gitlab.com/pypr/convergence_sph}. In the next section we discuss
the SPH method briefly and then proceed to look at the SPH kernel interpolation.

\section{Second order convergent WCSPH schemes}
\label{sec:sph}

We define the SPH approximation of any scalar (vector) field $f$ ($\ten{f}$)
in a domain $\Omega$ by
\begin{equation}
  \left<f(\ten{x})\right> = \int_{\Omega} f(\tilde{\ten{x}}) W(\ten{x} - \tilde{\ten{x}}, h) d \tilde{\ten{x}},
  \label{eq:sph}
\end{equation}
where $\ten{x}, \tilde{\ten{x}} \in \Omega$, $W$ is the kernel function, and
$h$ is the support radius of the kernel. It is well
known~\cite{monaghan-review:2005, violeau2012fluid} that for a symmetric
kernel which satisfies $\int W(\ten{x}) d\ten{x}= 1$ that,
\begin{equation}
  \label{eq:sph_approx}
  f(\ten{x}) = \left< f(\ten{x})\right> + O(h^2).
\end{equation}
Some of the widely used kernels are Gaussian~\cite{monaghan-review:2005}, cubic
spline~\cite{sph:fsf:monaghan-jcp94}, quintic spline, and Wendland
quintic~\cite{wendland_piecewise_1995}. We note that in this work we take $h$ to
be a constant.

When the kernel support is completely inside the domain boundary then we
can evaluate the gradient of a function by taking the gradient of the
kernel inside the integral. This approximation is also second-order in
$h$~\cite{violeau2012fluid}. In order to compute gradient numerically, we
discretize the domain $\Omega$ using particles having mass $m$, and density
$\rho$. The discretization of the domain into particles introduces
additional error in the approximation and is discussed in
\citet{quinlan_truncation_2006}. We can approximate the gradient of $f$ as,
\begin{equation}
  \begin{split}
  \nabla f(\ten{x}_i) =& \sum_j f(\ten{x}_j) \nabla W_{ij} \omega_j + |\nabla^3 f(\ten{x}_i)|  O(h^2) +\\
  ~& |\nabla f(\ten{x}_i)| O\left(\left(
  \frac{\Delta s}{h} \right)^{\beta + 4} \right),\\
  \end{split}
  \label{eq:sph_df}
\end{equation}
where $W_{ij} = W(\ten{x}_i - \ten{x}_j, h)$, $\omega_j=\frac{m_j}{\rho_j}$
is a measure of the volume of the particle, $\beta$ is the smoothness of
the kernel at the edge of its support, and the sum is taken over all the
particles under the support of the kernel. The value of $\beta$ is defined
as the smallest order of derivative of the kernel at the edge of its
support that is non-zero. For example, $\beta=3$ for cubic spline kernel,
and $\beta=5$ for quintic spline kernel. We note that the kernel gradient
is second order accurate only for a uniform distribution of
particles~\cite{fatehi_error_2011,
bonet_lok:cmame:1999,liu_restoring_2006}. However, many
authors~\cite{bonet_lok:cmame:1999,liu_restoring_2006,dilts1999moving} have
proposed methods to obtain second order convergent approximation of the
gradient of a function irrespective of the particle distribution.
\citet{fatehi_error_2011} obtained the error in approximation given by

\begin{equation}
  \begin{split}
  \nabla f(\ten{x}_i) =& \sum_j f(\ten{x}_j) \tilde{\nabla} W_{ij} \omega_j
  + |\nabla^3 f(\ten{x}_i)| O(h^2) +\\
  ~&  | \nabla^2 f(\ten{x}_i)|
  O\left(\tilde{\ten{d}_i} \left( \frac{\Delta s}{h} \right)^{\beta + 4}
  \right),\\
  \end{split}
  \label{eq:sph_df_corr}
\end{equation}
where $\tilde{\nabla} W_{ij} = B_i \nabla W_{ij}$, where $B_i$ is the
correction matrix, and $\tilde{\ten{d}_i}$ is the deviation of particle $i$
from its unperturbed location. We note that the error due to the
quadrature rule is retained.

The numerical volume $\omega$ in \cref{eq:sph_df_corr} is an approximation
and solely depends upon the spatial distribution of the particles. The
density $\rho$ may be computed for a particle using the summation density
as,
\begin{equation}
  \rho_i = \sum_j m_j W_{ij}.
  \label{eq:num_den}
\end{equation}
Therefore, for constant mass we may write the volume as,
\begin{equation}
 \omega_i = \frac{1}{\sum_j W_{ij}}.
 \label{eq:num_vol}
\end{equation}

Since the function and its derivative approximation depend on the kernel $W$,
support radius $h$, and the scaling factor $h_{\Delta s}=h/\Delta s$, we
perform a numerical study of the effect of the kernel on convergence next.

\subsection{Selection of approximating kernel}
\label{sec:comp_kernel}

\begin{figure}[ht!]
  \centering
  \includegraphics[width=\linewidth]{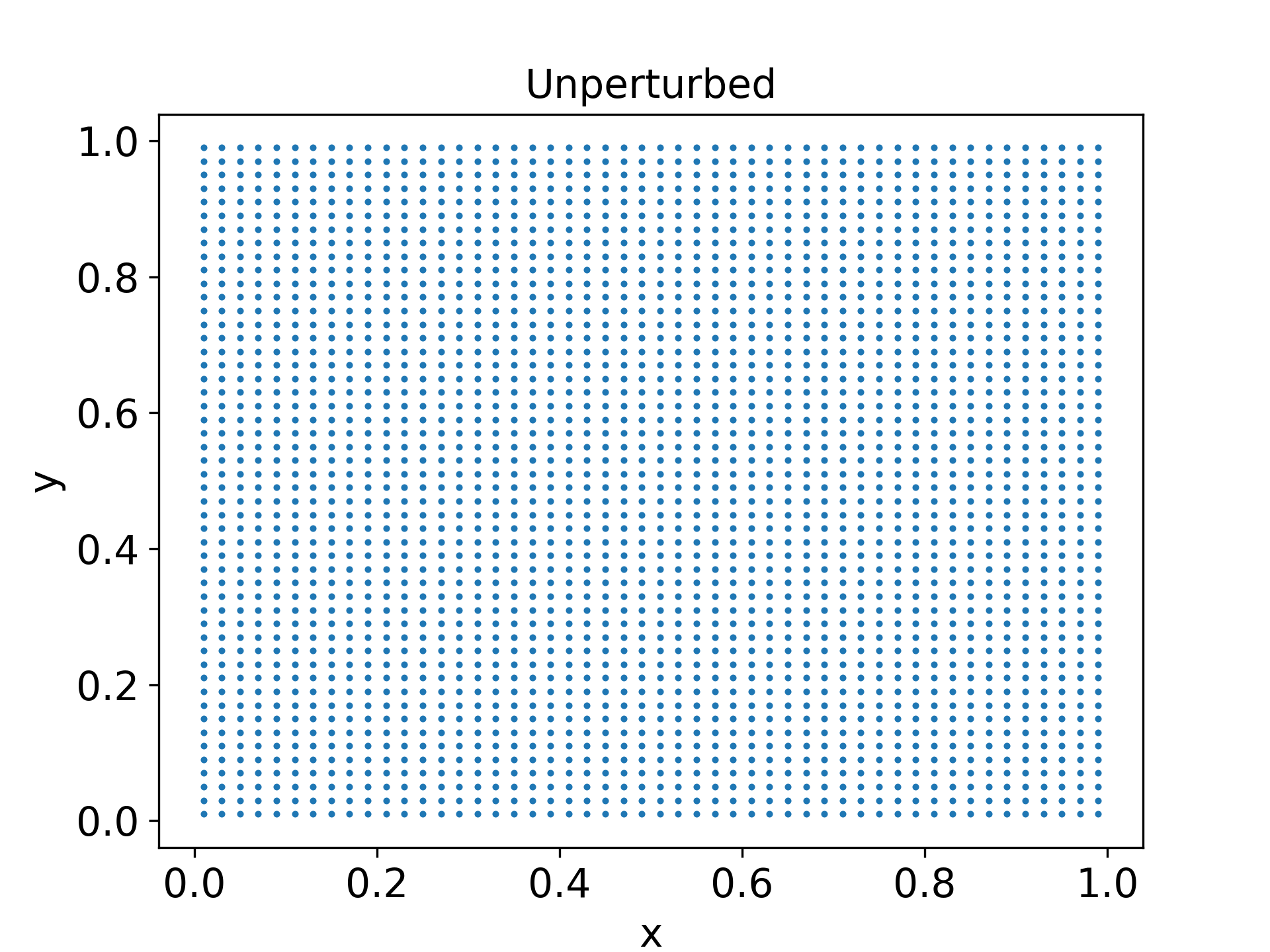}
  \includegraphics[width=\linewidth]{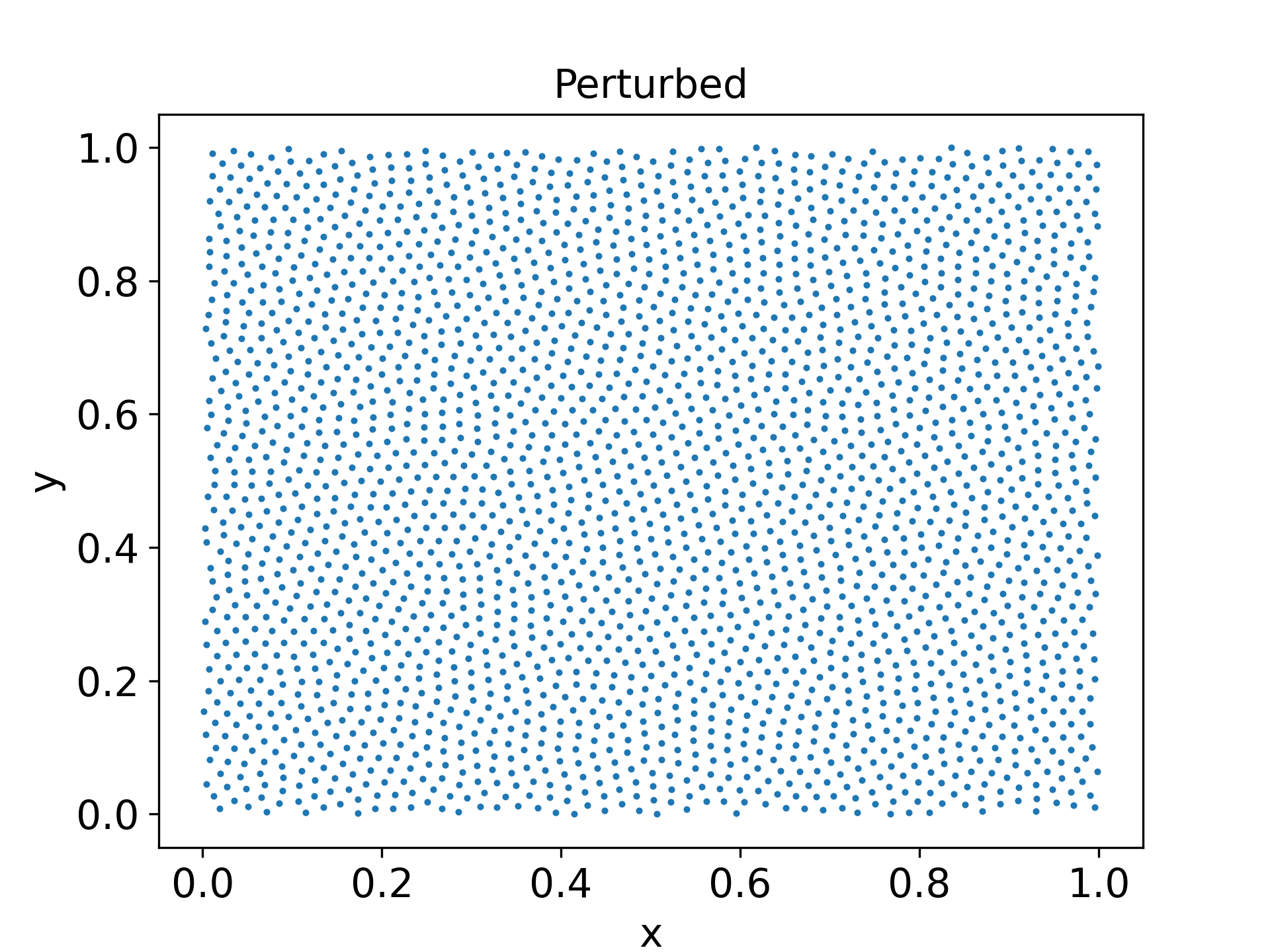}
  \caption{The unperturbed periodic particle and perturbed periodic particles.}
  \label{fig:par_plot}
\end{figure}

In this section, we compare various SPH kernels for their accuracy and order
of convergence in a discrete domain. We evaluate the error in function and
derivative approximation in a two-dimensional periodic domain. We simulate
periodicity by copying the appropriate particles and their properties near the
boundary such that the boundary particles have full support
\cite{randles1996smoothed}. The particles are either placed in a uniform mesh
(unperturbed) or in a \emph{packed} arrangement referred as unperturbed
periodic (UP) or perturbed periodic (PP), respectively. In order to obtain the
packed configuration, the particles are slightly perturbed from a uniform mesh
and their positions are moved and allowed to settle into a distribution with a
nearly constant density using a particle packing
algorithm~\cite{negi2019improved}. The algorithm effectively ensures that the
particles are not clustered and have minimal density variations. This mimics
the effect of many recent particle shifting
algorithms~\cite{colagrossi2012particle,lind2012incompressible,huang2019kernel}.
In \cref{fig:par_plot}, we show both the domains.

We compare the $L_1$ error in function and derivative approximation using
various approximating kernels commonly used in SPH with different support
radii. We present a detailed analysis in \cref{apn:kern}. The following is a
summary of the analysis:
\begin{itemize}
\item Errors for function approximation on an UP domain are
not affected by the choice of kernel.
\item The error increases with the increase in the $h_{\Delta s}$.
\item The error in a PP domain is dominated by the discretization error at
higher resolutions.
\end{itemize}

\begin{figure}[ht!]
  \centering
  \includegraphics[width=\linewidth]{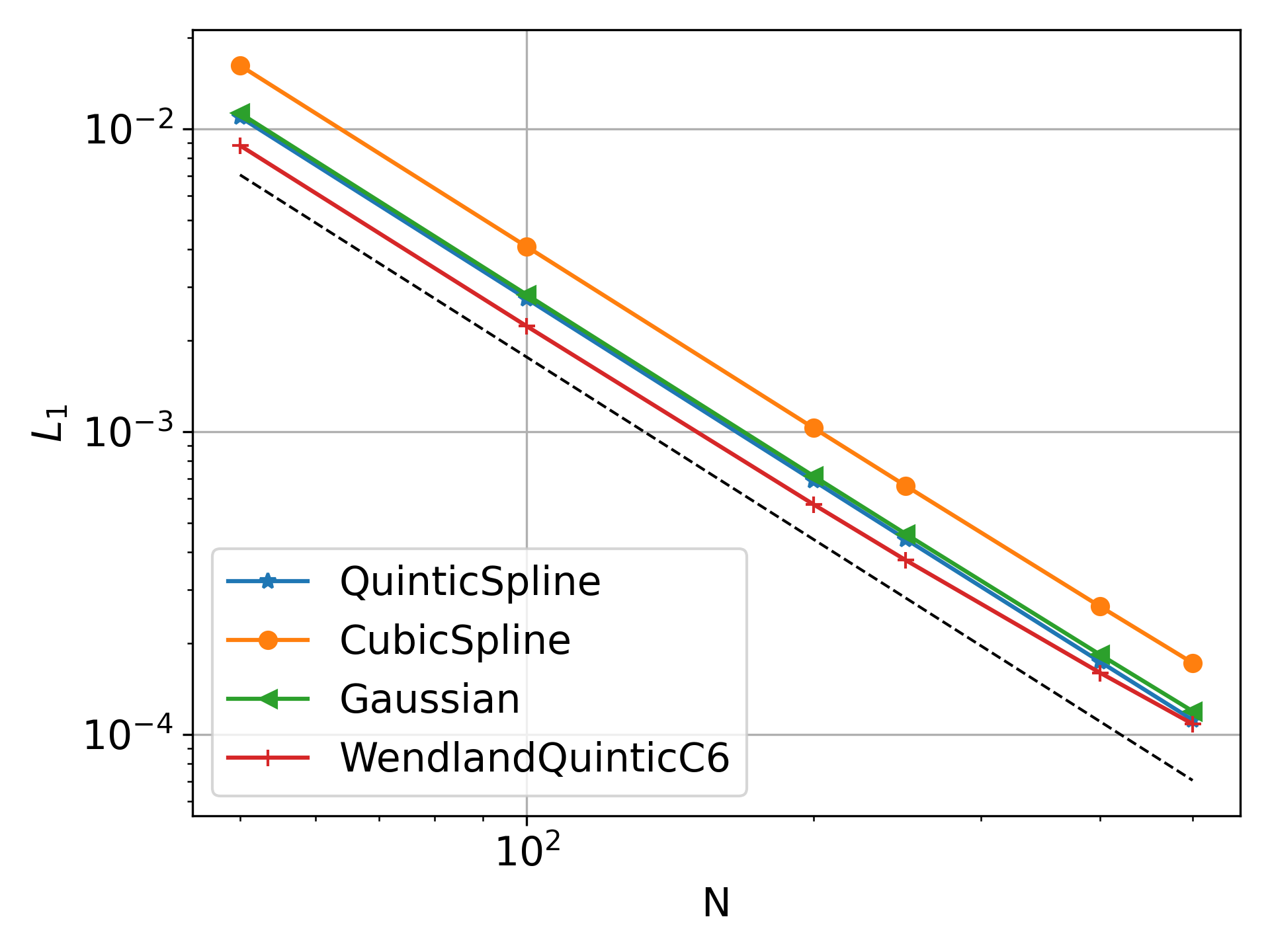}
  \caption{The convergence of the derivative approximation with different
    kernels when the gradient correction of \citet{bonet_lok:cmame:1999} is
    used on a PP domain. Here we use $h_{\Delta s}=1.2$. }
  \label{fig:d1p1_df_corr}
\end{figure}

In order to study the effect of kernel gradient correction, we apply the
correction proposed by \citet{bonet_lok:cmame:1999} for all the selected
kernels. The derivative approximation with correction is given by
\begin{equation}
  \left< \nabla f(\ten{x}_i) \right> = \sum_j (f(\ten{x}_j) - f(\ten{x}_i))
  B_i \nabla W_{ij} \omega_j.
  \label{eq:sph_dfij_corr}
\end{equation}
In \cref{fig:d1p1_df_corr}, we plot the $L_1$ error in the derivative
approximation as a function of resolution.  Clearly, all the kernels
Gaussiani ($G$), Wendland quintic 6\textsuperscript{th} order ($WQ_6$),
cubic spline ($CS$), and quintic spline ($QS$) show more or less the same behavior. Thus,
we can choose any of these kernels for our convergence study of the WCSPH
schemes. In the figure, it can be seen that the $WQ_6$ and $G$ kernels do
not sustain the second-order behavior. Therefore in this work, we choose
the $QS$ kernel with $h_{\Delta s} = 1.2$ for all the test cases
henceforth.

\subsection{Considerations while applying kernel gradient correction}
\label{sec:kern_grad}

The SPH method is widely used to solve fluid flow problems. In this work, we
focus on weakly-compressible SPH schemes that are used to simulate
incompressible fluid flows. We write the Navier-Stokes equation for a weakly
compressible flow along with the equation of state (EOS) as,
\begin{subequations}
  \begin{align}
  \frac{d \Varrho}{dt} &= - \Varrho \nabla \cdot \ten{u}, \label{eq:ns_wc:cont}\\
  \frac{d \ten{u}}{dt} &= - \frac{\nabla p}{\Varrho} + \nu \nabla^2 \ten{u},\label{eq:ns_wc:mom}\\
  p &= p(\Varrho, \Varrho_o, c_o) \label{eq:ns_wc:eos},
  \end{align}
  \label{eq:ns_wc}
\end{subequations}
where $\Varrho$, $\ten{u}$, and $p$ are the fluid density, velocity, and
pressure, respectively, $\nu$ is the kinematic viscosity, $\Varrho_o$ is the
reference fluid density, and $c_o$ is the artificial speed of sound in the
fluid.
  We note that the fluid density, $\Varrho$ is independent of the summation
  density $\rho$ (\cref{eq:num_den}). Normally, in SPH simulations these two are
  treated as the same and we discuss the reasons behind this choice in
  \cref{sec:comp_scheme}. The property, $\Varrho$ does not depend upon
  particle configuration and should be prescribed as an initial condition.

\begin{table}
  \centering
  \renewcommand{\arraystretch}{2}
  \begin{tabular}{l|l|l}
    \text{Name} & \text{Expression} & \text{Used in} \\
    \hline
    \text{\nam{sym1}} & \text{$\frac{(p_j + p_i)}{\rho_i \rho_j} \nabla W_{ij} m_j$} & \text{$\delta^+$ SPH \cite{sun2018numerical}} \\
    \text{\nam{sym2}} & \text{$m_j \left(\frac{p_j}{\rho_j^2} + \frac{p_i}{\rho_i^2}\right) \nabla W_{ij}$,} & \text{WCSPH \cite{wcsph-state-of-the-art-2010},} \\
    ~ & \text{$\frac{V_i^2 + V_j^2}{m_i} \tilde{p}_{ij} \nabla W_{ij} $} & \text{TVF \cite{Adami2013,zhang_hu_adams17,edac-sph:cf:2019}, ISPH \cite{isph:shao:lo:awr:2003}} \\
    \text{\nam{asym}} & \text{$\frac{(p_j - p_i)}{\Varrho_i} \nabla W_{ij} \omega_j$} & \text{WCSPH~\cite{monaghan-review:2005}} \\
  \end{tabular}
  \renewcommand{\arraystretch}{1}
  \caption{Different gradient approximations for $\frac{\nabla p}{\Varrho}$. The
  column ``expression'' is assumed to be summed over the index $j$ over all the
  neighbor particles inside the kernel support. The term $\Varrho=\rho$ for
  gradient comparison. $\tilde{p}_{ij}=\frac{p_i \rho_i + p_j \rho_j}{\rho_i +
  \rho_j}$ is the density averaged pressure.}
  \label{tab:grad}
\end{table}

There are many different ways to discretize \cref{eq:ns_wc} as can be seen
from
\cite{wcsph-state-of-the-art-2010,sun2018numerical,Adami2013,violeau2012fluid,crksph:jcp:2017,sun2017deltaplus}.
One of the key features of the discretization of the momentum equation (\cref{eq:ns_wc:mom}) is to ensure linear momentum conservation.
However, some researchers trade conservation for better
accuracy~\cite{sun2017deltaplus} and others use a conservative form that is
not as accurate~\cite{crksph:jcp:2017}. In view of this, we consider both
conservative and non-conservative discretizations in the present study.

In \cref{tab:grad}, we show the various pressure gradient approximations
employed in different SPH schemes. The conservative forms are usually
symmetric therefore are referred as \nam{sym}, the non-conservative forms are
asymmetric therefore referred as \nam{asym}. We compare the error in the
gradient approximation on both PP and UP domains with and without using a
corrected kernel gradient in \cref{apn:grad_comp}. We observe that only
\nam{asym} formulations can be corrected using all the corrections present in
the SPH literature. Whereas, \nam{sym1} can be corrected by the \textit{Liu
  correction} \cite{liu_restoring_2006} only.

In order to explain this behavior, we take first order Taylor series
approximation of a function, $f$ defined at $\ten{x}$ about $\ten{x}_i$ given by,
\begin{equation}
  f(\ten{x}) = f(\ten{x}_i) + (\ten{x} - \ten{x}_i) \cdot \nabla f(\ten{x}_i) + H.O.T
\end{equation}
integrating both sides with $\nabla W(\ten{x} - \ten{x}_i)$, we get
\begin{equation}
  \begin{split}
  \int f(\ten{x}) \nabla W d\ten{x} &= \int f(\ten{x}_i) \nabla W d\ten{x} + \int (\ten{x} - \ten{x}_i) \cdot \nabla f(\ten{x}_i) \nabla W d\ten{x} \\
  &= \int f(\ten{x}_i) \nabla W d\ten{x} + \int (\nabla W \otimes (\ten{x} - \ten{x}_i))  \nabla f(\ten{x}_i)  d\ten{x}. \\
  \end{split}
\end{equation}
Using one point quadrature approximation~\cite{dilts1999moving}, we get
\begin{equation}
  \begin{split}
\sum_j f_j \nabla W_{ij} \omega_j &= \sum_j f_i \nabla W_{ij} \omega_j + \sum_j \nabla W_{ij} \otimes (\ten{x}_j - \ten{x}_i) \nabla f(x_i) \omega_j  \\
\implies & \nabla f(x_i) = \sum_j (f_j - f_i) B_i \nabla W_{ij} \omega_j,\\
  \end{split}
  \label{eq:bonet_ts}
\end{equation}
where $B_i = \left(\sum_j \nabla W_{ij} \otimes (\ten{x}_j - \ten{x}_i)
\right)^{-1}$ is the correction matrix proposed by
\citet{bonet_lok:cmame:1999}. Clearly, the first order taylor-series
automatically suggests correction proposed in \onlinecite{bonet_lok:cmame:1999}
on an $\nam{asym}$ formulation. The \cref{eq:bonet_ts} is $O(h^2)$ accurate
\cite{fatehi_error_2011}.

On the other hand, the correction proposed by \citet{liu_restoring_2006},
originates by convolving the Taylor series with $W(\ten{x} - \ten{x}_j)$ and
$\nabla W(\ten{x} - \ten{x}_j)$, and solving all the equation
simultaneously. The matrix form is given by
\begin{widetext}
\begin{equation}
  \left[\begin{array} {cccc} {W_{kl} V_{l}} & {x_{lk}  W_{kl} V_{l}} &
  {y_{lk} W_{kl} V_{l}} & {z_{lk} W_{kl} V_{l}} \\ {W_{kl,x} V_{l} } &
  {x_{lk}W_{kl,x} V_{l}} & {y_{lk}W_{kl,x} V_{l}} & {z_{lk}W_{kl,x} V_{l}}
  \\ { W_{kl,y} V_{l}} & {x_{lk}W_{kl,y} V_{l}} & {y_{lk}W_{kl,y}V_{l}} &
  {z_{lk}W_{kl,y}V_{l} } \\ {W_{kl,z}V_{l}} & {x_{lk}W_{kl,z}V_{l}} &
  {y_{lk}W_{kl,z}V_{l} } & {z_{lk}W_{kl,z}V_{l}} \end{array} \right] \left[
  \begin{array} {c} {f_{k}} \\ {f_{k,x}} \\ {f_{k,y}} \\ {f_{k,z}}
\end{array} \right] = \left[\begin{array} {c} {f_{l}W_{kl}V_{l} } \\
  {f_{l}W_{kl,x}V_{l}} \\ {f_{l}W_{kl,y}V_{l} } \\ {f_{l}W_{kl,z}V_{l} }
\end{array} \right],
  \label{eq:liu_mat}
\end{equation}
\end{widetext}
where, $k$ is the destination particle index, $l$ is the neighbor particle
index, $W_{kl, \beta}$ for $\beta \in {x, y, z}$ is the kernel gradient
component in the $\beta$ direction. All the terms containing $l$ are summed
over all the neighbor particles. On solving the \cref{eq:liu_mat}, we
obtain a first order consistent gradient \cite{liu_restoring_2006}. For a
constant field, this method ensures that we satisfy $\sum \tilde{W}_{ij}
\omega_j=1$ and $\sum \nabla \tilde{W}_{ij} \omega_j = 0$, where
$\tilde{W}$ is the corrected kernel.  Therefore, with this correction in
both the \nam{sym1} and \nam{asym} forms the second term (i.e. $*_i$)
becomes zero, and we get the SOC approximation. Whereas, in \nam{sym2} the
term $p_i/\rho_i^2$ does not become zero, thus even this correction fails
to correct the approximation.

Using the Taylor series expansion, \citet{fatehi_error_2011} derived a
correction for the Laplacian operator. In \cref{apn:fatehi} we show the error
due to operators proposed by \citet{cleary1999conduction} and
\citet{fatehi_error_2011}. In \cref{apn:comp_disc}, we compare gradient,
divergence and Laplacian approximation with corrections proposed by various
authors in SPH literature. The comparison shows that the kernel gradient
correction must be used appropriately for second order convergence. However,
in case of divergence approximation the particle distribution plays a major
role that we discuss next.

\subsection{Considerations for the initial particle distribution}
\label{sec:div_error}

The particle distribution plays an important role in the error estimation
of divergence approximation. In this section, we use first order Taylor
series approximation to obtain error in divergence approximation as done in
previous section. We consider a two-dimensional velocity field. We write the
error $Er$, in the divergence evaluation as
\begin{equation}
  Er_i = \nabla \cdot \ten{u}_i - \sum_j (\ten{u}_j - \ten{u}_i) \cdot \nabla W_{ij} \omega_j,
\end{equation}
Using first order Taylor-series expansion of $\ten{u}_j$ about the point $\ten{x}_i$,
\begin{equation}
  \ten{u}_j = \ten{u}_i - (\ten{x}_{ij}\cdot \nabla) \ten{u}_i,
\end{equation}
we write,
\begin{equation}
  \begin{split}
  Er_i = & \left(1 - \sum_j x_{ij} \frac{\partial W_{ij}}{\partial x} \omega_j \right) \frac{\partial u_i}{\partial x} + \left(1 - \sum_j y_{ij} \frac{\partial W_{ij}}{\partial y} \omega_j \right) \frac{\partial v_i}{\partial y} \\
  ~ & - \sum_j y_{ij} \frac{\partial u_{i}}{\partial y} \frac{\partial W_{ij}}{\partial x} \omega_j - \sum_j x_{ij} \frac{\partial v_{i}}{\partial x} \frac{\partial W_{ij}}{\partial y} \omega_j.
  \end{split}
  \label{eq:div_err}
\end{equation}
In the case of a UP domain, in \cref{eq:div_err}, the last two terms are
exactly zero and the coefficient of the first two terms are of equal
magnitude. Furthermore, since for a divergence-free velocity field,
$\frac{\partial u_i}{\partial x} = - \frac{\partial v_i}{\partial y}$, the
overall error becomes zero. On the other hand, in a PP domain, the last two
terms are of equal magnitude thus cancel, and the first two terms are
different to the order $10^{-4}$ (see \cref{apn:comp_div}) which
becomes the leading error term. Thus, we always get an error of the order of
$10^{-4}$ even after applying the Bonet correction. As far as we are aware
there are no known SPH discretizations which can resolve this issue using a
simple correction as done in case of gradients. This is a possible avenue
for future research.

\subsection{Minimal requirements for a SOC scheme}
\label{sec:comp_scheme}

In this section, we discuss strategies to obtain a SOC scheme for weakly
compressible fluid flows. We consider the fluid density as a property,
$\Varrho$ carried by a particle. The numerical density, $\rho$, and volume,
$\omega$ are a function of the surrounding particle distribution. The mass,
$m$ of the particles satisfies $m_i = \Varrho_i V_i = \rho_i \omega_i$ where
$V_i$ is the physical volume occupied by the particle and $\omega_i$ is the
numerical volume used for integration. Thus, we can approximate the fluid
density using the standard SPH approximation given by
\begin{equation}
  \Varrho_i = \sum_j \Varrho_j W_{ij} \omega_j.
\end{equation}

In case of weakly compressible SPH, the requirement of linear momentum
conservation condition may be relaxed and is only satisfied
approximately~\cite{oger_improved_2007,bonet_lok:cmame:1999,sun2017deltaplus}.
Therefore, we use the SOC approximations that are non-conservative
as discussed in \cref{apn:comp_disc}. In \cref{tab:ope_disc}, we list all the
discretizations that we can employ to obtain a SOC WCSPH scheme.

\begin{table}[h!]
  \centering
  \renewcommand{\arraystretch}{2}
  \begin{tabular}{l|l}
    \text{Operators} & \text{Possible discretization for SOC}\\
    \hline
    \text{Gradient} & \text{\nam{asym\_c}, \nam{sym1\_l}}\\
    \text{Divergence} & \text{\nam{div\_c}} \\
    \text{Laplacian} & \text{\nam{coupled\_c}, \nam{Fatehi\_c}, \nam{Korzilius}}\\
  \end{tabular}
  \renewcommand{\arraystretch}{1}
\caption{The operators and their discretization suitable for a SOC scheme
(For details refer \cref{apn:comp_disc}).}
  \label{tab:ope_disc}
\end{table}

Furthermore, one can solve the fluid flow equations by using a Lagrangian
approach as well as an Eulerian approach. We discuss the scheme for both these
cases in the following sections.

\subsubsection{SOC for Lagrangian WCSPH}
\label{sec:lagrangian}

In the Lagrangian description, the continuity equation and the momentum equation are given by
\begin{equation}
  \begin{split}
    \frac{d \Varrho}{dt} &= -\Varrho \nabla \cdot \ten{u}\\
    \frac{d \ten{u}}{dt} &= - \frac{\nabla p}{\Varrho} + \nu \nabla^2 \ten{u},\\
  \end{split}
  \label{eq:wcsph_eq}
\end{equation}
In order to evaluate the RHS of the above equations, one may employ any method
listed in \cref{tab:ope_disc}. The pressure is evaluated using an equation of
state given by
\begin{equation}
  p = \frac{\Varrho_o c_o^2}{\gamma} \left( \left(\frac{\Varrho}{\Varrho_o}\right)^\gamma - 1\right),
  \label{eq:tait_eos}
\end{equation}
where $\gamma=7$, $\Varrho_o$ is the reference density, and $c_o$ is the reduced
speed of sound. We note that the linear equation of state where in
\cref{eq:tait_eos}, $\gamma=1$, works equally well. We integrate the particles
in time using a Runge-Kutta $2^{nd}$ order integrator.

Since we use an asymmetric form of the pressure gradient approximation,
particles tend to clump together due to absence of a redistributing background
pressure \cite{sun2018numerical}. We use the iterative particle shifting
proposed by \citet{huang2019kernel} after every few iterations to redistribute
the particles. We compute the shifting vector for the $m^{th}$ iteration (of
the shifting iterations) using
\begin{equation}
  \delta \ten{x}_i^{m} = h_i \sum_j \ten{n}_{ij} W_{ij} \omega_j,
\end{equation}
where $\ten{n}_{ij} = \ten{x}_{ij}/|\ten{x}_{ij}|$. The new particle position,
\begin{equation}
  \tilde{\ten{x}}_{i}^{m+1} = \ten{x}_{i}^{m} + \delta \ten{x}_{i}^{m}
\end{equation}
is computed. The particles are shifted until the criterion,
\begin{equation}
  |\max(\chi^m) - \chi_o| < \epsilon
\end{equation}
is satisfied up to a maximum of 10 iterations, where $\chi^m = h^2 \sum_j
W_{ij}$, $\chi_o$ is the value for uniform distribution, and $\epsilon$ is
an adjustable parameter. In order to keep the approximation of the particle
$O(h^2)$ accurate, we update the particle properties after shifting by,
\begin{equation}
  \phi(\tilde{\ten{x}}_i) = \phi({\ten{x}}_i) + (\tilde{\ten{x}}_i - \ten{x}_i) \cdot \nabla \phi(\ten{x}_i),
  \label{eq:update}
\end{equation}
where $\tilde{\ten{x}}_i$ is the final position after iterative shifting,
$\phi$ is the property to be updated, and $\nabla \phi(\ten{x}_i)$ is the
gradient of the property on the last position computed with the Bonet
correction. In a variation of the above scheme discussed in
\cref{sec:new_schemes}, we observe that usage of non-iterative PST
proposed by \citet{sun2019consistent} results in slightly higher errors but
still retains its SOC. We refer to the scheme discussed above as L-IPST-C
(Lagrangian with iterative PST and \nam{coupled\_c} viscosity formulation).
Similarly the method using \nam{Fatehi\_c} and \nam{Korzilius} formulation are
referred as L-IPST-F and L-IPST-K respectively. We note that we only perform
the IPST step every $10$ timesteps rather than at every timestep.

\subsubsection{SOC for Eulerian WCSPH}
\label{sec:eulerian}

In the Eulerian description, the continuity equation is written as
\begin{equation}
  \frac{\partial \Varrho}{\partial t} = - \Varrho \nabla \cdot \ten{u} - \ten{u}\cdot \nabla \Varrho.
\end{equation}
Since the fluid density $\Varrho$ is not same as the particle density $\rho$
we do not ignore this term, this is unlike what is done by \citet{nasar2019}.
The momentum equation is written as
\begin{equation}
  \frac{\partial \ten{u}}{\partial t} = -\frac{\nabla p}{\Varrho} + \nu \nabla^2 \ten{u} - \ten{u} \cdot \nabla \ten{u}.
\end{equation}
We discretize all the terms using SOC operators listed in
\cref{tab:ope_disc}. We perform time integration using the RK2 integrator;
however, we note that the positions of the particles are not updated.

\subsection{The effect of $c_o$ on convergence}
\label{sec:c0_effect}

In the schemes discussed in the previous section, we impose artificial
compressibility (AC) using the EOS, which is $O(M^2)$ accurate
\cite{chorin_numerical_1967,Clausen2013}, where $M=U_{max}/c_o$ is the Mach
number of the flow. \citet{chorin_numerical_1967} originally proposed this
method to obtain steady-state solutions of an incompressible flow. Some
authors have used artificial compressibility with dual-time stepping to
achieve truly incompressible time-accurate results
\cite{fatehi-2019,ramachandran_dual-time_2021}. We achieve the
incompressibility limit when $c_o \to \infty$. Therefore, in order to
increase the accuracy at higher resolution a higher speed of sound must be
used. We show the effect of the speed of sound on accuracy in
\cref{sec:var_co}.

\subsection{Variations of the SOC scheme}
\label{sec:new_schemes}

In this section, we show that the scheme presented in the
\cref{sec:lagrangian} (L-IPST-C) can be easily converted into other forms for
improved accuracy and ease of calculation. We note that regardless of the set
of governing equations employed, the discretizations from \cref{tab:ope_disc}
must be used to achieve SOC.

In order to remove high frequency oscillations, one could modify the
continuity equation given by,
\begin{equation}
  \frac{d \Varrho}{dt} = -\Varrho \nabla \cdot \ten{u} + D \nabla^2 \Varrho
  \label{eq:deltasph}
\end{equation}
where, $D=\delta h c_o^2$ is the damping constant, where $\delta=0.1$. This
corresponds to the $\delta$-SPH scheme~\cite{antuono-deltasph:cpc:2010}. In
this case we also use the linear equation of state to evaluate $p$ given by
\begin{equation}
  p = c_o^2 (\Varrho - \Varrho_o).
  \label{eq:linear_eos1}
\end{equation}

The following are different variations of the basic scheme:
\begin{enumerate}
\item Using different PST : One could use either IPST proposed by
  \citet{huang2019kernel} or the non-iterative PST proposed by
  \citet{sun_consistent_2019}. The properties like $u, v, p,$ and $\Varrho$
  need to be updated using first order Taylor expansions given by
  \begin{equation}
    \phi(\tilde{\ten{x}}_i) = \phi({\ten{x}}_i) + (\tilde{\ten{x}}_i - \ten{x}_i) \cdot \nabla \phi(\ten{x}_i)
  \end{equation}
  where $\phi$ is the desired property. We use the coupled formulation for the
  viscosity and non-iterative PST. We refer to this method as L-PST-C.

\item Using pressure evolution: On taking the derivative of EOS in
  \cref{eq:linear_eos1} w.r.t.\ time and using the \cref{eq:deltasph}, we get
  the pressure evolution equation given by
  \begin{equation}
    \frac{d p}{dt} = -\Varrho c_o^2 \nabla \cdot \ten{u} + D \nabla^2 p.
    \label{eq:edac}
  \end{equation}
  This is very similar to the EDAC pressure evolution~\cite{edac-sph:cf:2019}.
  The value of $\Varrho$ can be evaluated from the EOS in
  \cref{eq:linear_eos1} given by
  \begin{equation}
    \Varrho = \frac{p}{c_o^2} + \Varrho_o.
    \label{eq:eos_inv}
  \end{equation}
  We employ the coupled formulation for viscosity and use IPST for
  regularization. We refer to this method as PE-IPST-C.

\item Using remeshing for regularization: The regularization step performed
  using PST in the L-IPST-C method can be replaced with the remeshing
  procedure of \citet{hieberImmersedBoundaryMethod2008}. The remeshing is
  performed using the $M_4$ kernel given by,
  \begin{equation}
    M_4(q) = \begin{cases}
      1 - \frac{5 q^2}{2} + \frac{3 q^3}{2}& 0 \leq q < 1,\\
      \frac{(1 - q)(2-q)^2}{2}& 1  \leq q < 2, \\
      0 & q \geq 2,
    \end{cases}
  \end{equation}
  where $q = |\ten{x}|/\Delta s$, where $\Delta s$ is the initial particle
  spacing. The properties on the regular grid are computed using
  \begin{equation}
    \phi(\tilde{\ten{x}}_i) = \frac{\sum \phi(\ten{x}_j)
      M_4(|\tilde{\ten{x}}_i - \ten{x}_j|, h)}
    { \sum M_4(|\tilde{\ten{x}}_i - \ten{x}_j|, h)},
    \label{eq:m_4_remesh}
  \end{equation}
  where $\tilde{\ten{x}}$ are points on a regular Cartesian mesh. The
  remeshing procedure can be performed every few steps; however, we perform
  remeshing after every timestep. We use the coupled formulation for
  viscosity. We refer to this method as L-RR-C.

\item Including regularization in the form of shifting velocity: the methods
  of \onlinecite{Adami2013,oger_ale_sph_2016,sun_consistent_2019} use shifting by
  perturbing the velocity of the particles and adding corrections to the
  momentum equation. Thus the particles are advected using the transport
  velocity, $\tilde{\ten{u}} = \ten{u} + \delta \ten{u}$ and the displacement
  is given by,
  \begin{equation}
    \ten{x}_i^{n+1} = \ten{x}_i^{n} + \Delta t (\ten{u}_i + \delta \ten{u}_i).
    \label{eq:tv}
  \end{equation}
  The new continuity and momentum equations are given by
  \begin{equation}
    \begin{split}
      \frac{\tilde{d} \Varrho}{dt} &= -\Varrho \nabla \cdot \ten{u} + D \nabla^2 \Varrho + \delta \ten{u} \cdot \nabla \Varrho,\\
      \frac{\tilde{d} \ten{u}}{dt} &= - \frac{\nabla p}{\Varrho} + \nu \nabla^2 \ten{u} + \delta \ten{u} \cdot \nabla \ten{u},\\
    \end{split}
    \label{eq:tv_eq}
  \end{equation}
where $\frac{\tilde{d}(\cdot)}{dt} = \frac{\partial (\cdot) }{\partial t} +
\tilde{\ten{u}} \cdot \nabla (\cdot)$. In this method, we employ SOC
approximations mentioned in \cref{tab:ope_disc} along with correction
proposed in \cref{apn:delta_plus_corr}. We use the PST proposed in
\onlinecite{sun_consistent_2019} for this scheme along with the coupled
formulation for viscosity. We refer to this method as TV-C.

\item Eulerian method: The Eulerian method solves the equation of motion on a
  stationary grid. This can be derived using the TV-C method by setting
  $\delta \ten{u}_i = - \ten{u}_i$. This substitution makes the transport
  velocity in the \cref{eq:tv} equal to zero, thus the particle does not move.
  The modified equation on setting $\delta \ten{u}_i = - \ten{u}_i$ in
  \cref{eq:tv_eq}, we get
  \begin{equation}
    \begin{split}
      \frac{\partial \Varrho}{\partial t} &= -\Varrho \nabla \cdot \ten{u} + D \nabla^2 \Varrho - \ten{u} \cdot \nabla \Varrho,\\
      \frac{\partial \ten{u}}{\partial t} &= - \frac{\nabla p}{\Varrho} + \nu \nabla^2 \ten{u} - \ten{u} \cdot \nabla \ten{u}.\\
    \end{split}
    \label{eq:ewcsph}
  \end{equation}
Therefore, we recover the governing equation for the Eulerian method (See
\cref{sec:eulerian}). We note that unlike \cite{nasar2019}, we retain the
last term in the continuity equation. We use the coupled formulation to
discretize viscous term. We refer to this method as E-C.

\end{enumerate}

\section{Results and discussions}
\label{sec:result}

In this section, we compare the solution obtained from different schemes for
the Taylor-Green, Gresho-Chan vortex, and incompressible shear layer problems.
We first compare the $L_1$ error in velocity, pressure, and linear and angular
momentum conservation of the L-IPST-C with various existing schemes. In order
to observe the effect of $c_o$ on the convergence, we solve the Taylor-Green
problem with different speeds of sound using the L-IPST-C and L-IPST-F
schemes. For the highest value of $c_o=80 m/s$, we compare the results using
different variations of the SOC schemes. In order to observe the conservation
property, we compare the solutions for inviscid problems viz.\ incompressible
shear layer and Gresho-Chan vortex using existing schemes as well as the SOC
schemes. Furthermore, we compare the SOC scheme and existing schemes for long
time simulations for all the test cases. Finally, we compute the cost of
computation versus accuracy for all the schemes.

We implement the schemes using the open source PySPH~\cite{pysph2020}
framework and automate the generation of all the figures presented in this
manuscript using the \texttt{automan} framework~\cite{pr:automan:2018}. The
source code is available at \url{https://gitlab.com/pypr/convergence_sph}.

\subsection{Comparison with existing SPH schemes}
\label{sec:comp_other_scheme}

In this section, we compare the following schemes:
\begin{enumerate}
  \item TVF : Transport velocity formulation proposed by \citet{Adami2013}.
  \item $\delta^+$SPH : The improved $\delta$-SPH formulation proposed by
  \citet{sun_consistent_2019}.
  \item EDAC : Entropically Damped artificial compressibility SPH
  formulation proposed by \citet{edac-sph:cf:2019}.
  \item EWCSPH : The Eulerian SPH method proposed by \citet{nasar2019}.
\item L-IPST-C : The Lagrangian method with iterative PST and
\nam{coupled\_c} viscosity formulation discussed in \cref{sec:lagrangian}.
\end{enumerate}

In order to compare these schemes, we consider the Taylor-Green vortex
problem. We choose this problem, since it is periodic, has no solid
boundaries, and admits an exact solution \footnote{In this paper, we do not
consider solid boundaries since to our knowledge, no second-order
convergent boundary condition implementations exist in the SPH literature.
Therefore the error due to the boundary will dominate. In order to show the
convergence with solid boundary, we study one popular boundary
condition~\cite{maciaTheoreticalAnalysisNoSlip2011,Adami2012} in
\cref{apn:solid_bc}.}. The solution of the Taylor-Green problem is given by
\begin{equation}
  \begin{split}
    u &= - U e^{bt} \cos(2 \pi x) \sin(2  \pi y),\\
    v &= U e^{bt} \sin(2 \pi x) \cos(2  \pi y),\\
    p &= -0.25 U^2 e^{2bt} (\cos(4 \pi x) +\cos(4  \pi y)),\\
  \end{split}
  \label{eq:tg_exact}
\end{equation}
where $b = -8 \pi^2 / Re$, where $Re$ is the Reynolds number of the flow.
We consider $Re=100$ and $U = 1 m/s$. For the Lagrangian schemes, we
consider a perturbed periodic (PP) arrangement of particles shown in the
\cref{fig:par_plot} for different resolutions. At $t=0$ we initialize the
pressure $p$ and velocity $(u, v)$ using \cref{eq:tg_exact} for all the
schemes. Since the fluid density $\Varrho$ is a function of pressure, we
initialize density inverting \cref{eq:tait_eos}. In the case of the EWCSPH
scheme, we consider an unperturbed periodic (UP) arrangement of particles
and initialize the $\rho$ using the prescribed pressure. We compute the
$L_1$ error in pressure and velocity by
\begin{equation}
  L_1(f, h) = \sum_j \sum_i \frac{|f(\ten{x}_i, t_j) - f_o(\ten{x}_i, t_j)|}{N} \Delta t
  \label{eq:tg_l1}
\end{equation}
where, $h=h_{\Delta s} \Delta s$ is the smoothing length of the kernel,
$\Delta t$ is the timestep, $N$ is the total number of particles in the
domain, $f$ is either pressure or velocity, and $f_o$ is the exact value
obtained using \cref{eq:tg_exact}. The particle spacing, $\Delta s$ is set
according to the resolution. We consider resolutions of $50\times50$ to
$500\times500$ particles in a $1m\times1m$ periodic domain. In order to isolate
the effect of spatial approximations on the convergence, we set the timestep
$\Delta t = 0.3 h/(U_{max} + c_o)$, where $h=1.2/500 m$ is set corresponding to
highest resolution, $c_o=10U$, for all the simulations. We run all the
simulations for $1$ timestep and observe convergence. We choose one timestep
since most of the schemes considered diverge.
\begin{figure*}[ht!]
  \centering
  \includegraphics[width=0.8\linewidth]{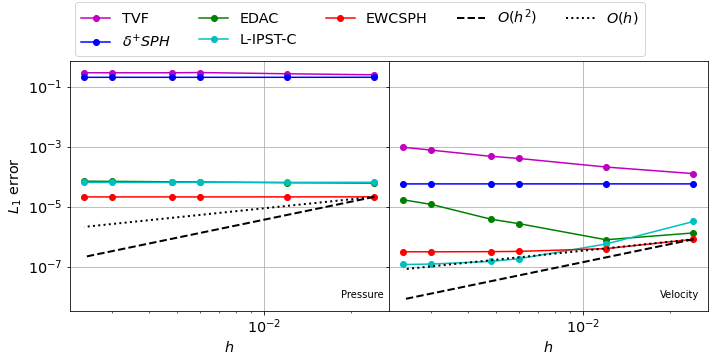}
  \caption{Convergence of $L_1$ error in pressure (left) and velocity (right)
    with the change in resolution. The $Re=100$, $c_o=10$, $\Delta
    t=6.54\times10^{-5}$ and only $1$ timestep taken.}
  \label{fig:tg_compare_scheme}
\end{figure*}

\begin{table}
  \centering
 \begin{tabular}{cllll}
\toprule
           Name & $\frac{F_T}{F_{max}}$ & $T_r$ & $L_1(|\ten{u}|)(O)$ &     $L_1(p)(O)$ \\
\midrule
$\delta^{+}SPH$ &              2.19e-05 &  1.49 &      5.69e-05(0.00) &  2.03e-01(0.00) \\
           EDAC &              1.88e-07 &  1.32 &     1.68e-05(-1.24) & 7.00e-05(-0.07) \\
         EWCSPH &              7.03e-15 &  2.25 &      3.13e-07(0.38) &  2.10e-05(0.00) \\
       L-IPST-C &              1.34e-05 &  3.36 &      1.18e-07(1.42) &  6.41e-05(0.00) \\
            TVF &              3.70e-16 &  1.00 &     9.45e-04(-0.88) & 2.88e-01(-0.07) \\
\bottomrule
\end{tabular}

  \caption{Table showing total force w.r.t.\ the maximum force in the domain
    and the time taken for $1$ iteration w.r.t.\ the TVF scheme for all the
    schemes.}
  \label{tab:tg_compare_time}
\end{table}

In the \cref{fig:tg_compare_scheme}, we plot the $L_1$ error evaluated using
\cref{eq:tg_l1} for pressure and velocity in the domain for different schemes.
Clearly, none of the schemes show convergence in pressure. This is because the
initial velocity is divergence-free, so there is no change in density and
thereby pressure. We observe that the EDAC, EWCSPH, and L-IPST-C schemes are
almost four orders more accurate than the TVF and $\delta^{+}$SPH schemes.  In
case of both the TVF and $\delta^{+}$SPH schemes, we link the pressure with
particle density $\rho$, which is a function of the particle configuration. The
particle positions are a result of the particle shifting, and therefore, the
pressure is incorrectly captured. On the other hand, the other schemes either
use a pressure evolution equation (EDAC) or a fluid density to evaluate
pressure. In the case of the EWCSPH scheme, we initialize density using the
pressure values in the \cref{eq:tait_eos} which results in better accuracy.

The $L_1$ error in velocity diverges in the case of the TVF and EDAC schemes
since these use a symmetric form of type \nam{sym2} in \cref{tab:grad} to
discretize the momentum equation. Whereas, in the case of the $\delta^{+}$SPH
scheme, \nam{sym1} type of discretization is employed leading to less errors.
Moreover, the $\delta^{+}$SPH scheme uses a consistent formulation and both TVF
and EDAC schemes are inconsistent when the shifting (transport) velocity is
added to the momentum equation \cite{sun_consistent_2019}. The EWCSPH, and
L-IPST-C formulations show convergence (not second-order) as expected. We
observe that in the velocity convergence a constant leading error term dominates
resulting in flattening at higher resolutions. Since, we use second-order
accurate formulations in L-IPST-C and EWCSPH \footnote{It is second-order
accurate since a uniform stationary grid is used.} formulations, the only
equation which is not converging with resolution is the equation of state (EOS).

In this section, we focus on highlighting the effect of using fluid density
$\Varrho$ different than the numerical density $\rho$. It allows for a
superior convergence rate and independence of density from particle positions.
In contrast to this, the use of numerical density as a function of particle
position is consistent with the volume used for the SPH approximation. In TVF,
$\delta^+$SPH, EDAC and, EWCSPH schemes, we make no such distinction, and use
the fluid density $\rho$ in the numerical volume $\omega_j=m_j/\rho_j$. The
poor convergence for these schemes show that it is important to treat the
fluid and numerical densities differently.

We also compare the linear momentum conservation and time taken to evaluate
the accelerations for the case with $500 \times 500$ particles. As shown in
\citet{bonet_lok:cmame:1999}, linear momentum is conserved when the total
force, $\sum_i F_i = 0$, where the sum is taken over all the particles and
$F_i= \frac{\nabla p_i}{\Varrho_i} + \nu \nabla^2 \ten{u}_i$. In the
\cref{tab:tg_compare_time}, we tabulate the total force and the time taken
by the scheme for one timestep with the errors and order of convergence in
pressure and velocity for the $500\times 500$ resolution case. It is clear
that the TVF and EWCSPH schemes conserve linear momentum and take the least
amount of time. The EDAC and the $\delta^+$SPH scheme do not conserve
linear momentum exactly. In the case of the EDAC scheme the use of average
pressure in the pressure gradient evaluation results in lack of
conservation. Whereas, in the case of $\delta^+$SPH the asymmetry of the
shifting velocity divergence causes lack of conservation. The L-IPST-C
scheme is known to be non-conservative; however, the value is comparable to
other schemes. The time taken by the L-IPST-C scheme is significantly
higher due to the evaluation of correction matrices.

\subsection{Convergence with varying speed of sound}
\label{sec:var_co}

\begin{figure*}[ht!]
  \centering
  \includegraphics[width=0.8\linewidth]{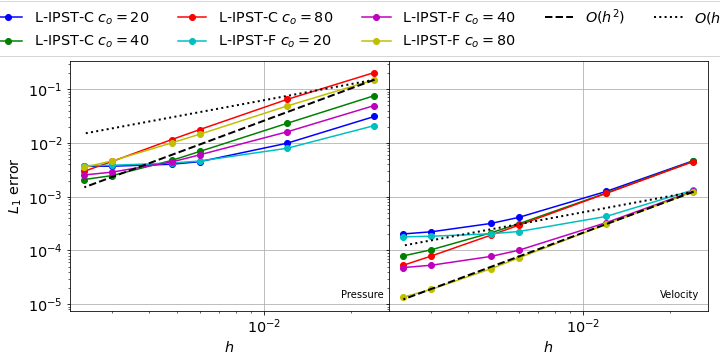}
  \caption{Convergence rates for pressure (left) and velocity (right). The
    L-IPST-C and L-IPST-F methods are compared for different values of $c_o$.}
  \label{fig:tg_conv_co}
\end{figure*}

\begin{table}[ht!]
  \centering
  \begin{tabular}{cllll}
\toprule
             Name & $\frac{F_T}{F_{max}}$ & $T_r$ & $L_1(|\ten{u}|)(O)$ &    $L_1(p)(O)$ \\
\midrule
L-IPST-C $c_o=20$ &              7.13e-05 &  1.00 &      2.03e-04(1.38) & 3.66e-03(0.93) \\
L-IPST-C $c_o=40$ &              9.31e-05 &  2.15 &      7.94e-05(1.78) & 2.09e-03(1.60) \\
L-IPST-C $c_o=80$ &              1.44e-04 &  3.76 &      5.32e-05(1.93) & 2.98e-03(1.85) \\
L-IPST-F $c_o=20$ &              7.11e-05 &  1.23 &      1.80e-04(0.85) & 3.77e-03(0.73) \\
L-IPST-F $c_o=40$ &              5.99e-05 &  2.84 &      4.81e-05(1.44) & 2.54e-03(1.31) \\
L-IPST-F $c_o=80$ &              1.66e-04 &  5.46 &      1.36e-05(1.98) & 3.52e-03(1.65) \\
\bottomrule
\end{tabular}

  \caption{Comparison of total force, time taken relative to the L-IPST-C with
    $c_o=20 m/s$, $L_1$ error in velocity and order for different values of
    $c_o$.}
  \label{tab:tg_compare_c0}
\end{table}

In this section, we compare the convergence of the L-IPST-F and L-IPST-C
scheme with change in AC parameter i.e.\ the speed of sound $c_o$. We
consider the Taylor-Green problem; however, we run the simulation for
$t=0.5s$. In \cref{fig:tg_conv_co}, we plot the $L_1$ error in the pressure
and velocity for both the schemes and different $c_o$. We observe that both
L-IPST-C and L-IPST-F methods are affected by the change in $c_o$ value
significantly as expected. In case of pressure, with the increase in the
$c_o$ value, the error in the lower resolutions increases; however, for
higher resolution it improves. Clearly, we attain the increase in the order
of convergence in case of the pressure due to increased error scales at
lower resolutions. The increase in error is attributed to the inability of
the SPH operators to correctly capture a divergence free velocity field as
discussed in \cref{sec:div_error}. However, on looking at the velocity
convergence, both the schemes attain SOC even at higher resolutions. We
observe, though, at lower resolutions, the error in the pressure increases
with the $c_o$ value.  As observed in the case of the Laplace operator
comparison in \ref{apn:comp_lap}, the use of \nam{Fatehi\_c}
discretization offers better accuracy.

In the \cref{tab:tg_compare_c0}, we tabulate the total force, relative time,
and the $L_1$ error in pressure and velocity with the order of convergence for
$500\times 500$ particles. We observe that at a higher $c_o$ value, the total
force is less compared to the simulation when $c_o$ values are lower. From the
table, we can see that the use of L-IPST-F scheme offers better accuracy at the
cost of the extra time taken. We also note that one can choose to use lower
values of $c_o$ at lower resolutions and increase the value as the resolution
increases to get better accuracy in pressure. Doing this does not affect the
convergence in velocity.

\subsection{Comparison of SOC variants}
\label{sec:comp_new}
\begin{figure*}[ht!]
  \centering
  \includegraphics[width=0.8\linewidth]{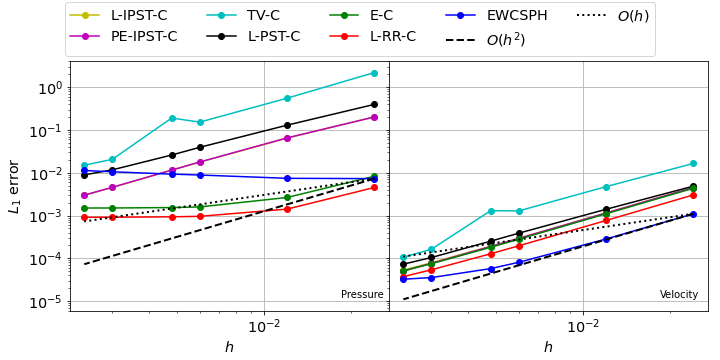}
  \caption{Convergence rates for pressure (left) and velocity (right) of
    different variants of SOC schemes.}
  \label{fig:tg_conv_var}
\end{figure*}

\begin{table}
  \centering
  \begin{tabular}{cllll}
\toprule
     Name & $\frac{F_T}{F_{max}}$ & $T_r$ & $L_1(|\ten{u}|)(O)$ &     $L_1(p)(O)$ \\
\midrule
      E-C &              8.13e-15 &  1.19 &      5.11e-05(1.94) &  1.50e-03(0.71) \\
   EWCSPH &              1.07e-14 &  1.00 &      3.25e-05(1.57) & 1.13e-02(-0.20) \\
 L-IPST-C &              1.44e-04 &  1.53 &      5.32e-05(1.93) &  2.98e-03(1.85) \\
  L-PST-C &              2.64e-06 &  2.02 &      7.31e-05(1.84) &  8.93e-03(1.68) \\
   L-RR-C &              1.51e-18 &  1.41 &      3.74e-05(1.92) &  9.11e-04(0.65) \\
PE-IPST-C &              6.17e-05 &  1.65 &      5.05e-05(1.96) &  3.01e-03(1.84) \\
     TV-C &             -5.01e-06 &  1.88 &      1.06e-04(2.18) &  1.51e-02(2.14) \\
\bottomrule
\end{tabular}

  \caption{Comparison of total force, time taken relative to L-IPST-C with
    $c_o=20$, $L_1$ error in velocity and pressure for variation of schemes
    with $c_o=80$.}
  \label{tab:tg_compare_var}
\end{table}

In this section, We simulate the Taylor-Green problem using $c_o=80 m/s$
for a duration of $0.5s$ with different resolutions discussed in
\cref{sec:new_schemes}. In addition, we study the performance of the EWCSPH
scheme since it is computationally efficient and accurate. In the
\cref{fig:tg_conv_var}, we plot the error in pressure and velocity for all
the schemes. In the \cref{tab:tg_compare_var}, we tabulate the total force,
relative time, $L_1$ error in pressure and velocity at $500\times500$
resolution, and the order of convergence.

The L-IPST-C and PE-IPST-C overlap in both the pressure and velocity convergence
plots, and these are both approximately second-order. Compared to L-IPST-C, the
L-PST-C shows lower convergence rate, and TV-C shows higher order of
convergence, whereas E-C and L-RR-C show very poor convergence rates in
pressure; however, the L-RR-C method shows very low errors in pressure. The
EWCSPH has a negative convergence rate in pressure.  While the TV-C shows a high
convergence rate, it has much larger errors than all the other schemes
considered for both pressure and velocity. The E-C, L-IPST-C, L-RR-C, and
PE-IPST-C show high convergence rates in velocity as expected. The L-PST-C shows
slightly high error and $1.84$ convergence rate. The EWCSPH shows a lower
convergence rate of $1.57$ but is the most accurate of all the schemes regarding
the velocity error.

The TV-C scheme shows low accuracy since we perform the shifting using an
additional term in the momentum equation compared to the PE-IPST-C and
L-IPST-C. This decrease in accuracy is also visible in the case of velocity.
The L-PST-C scheme show higher error suggesting that the non-iterative PST
does not perform the required amount of regularization. Both L-RR-C and E-C
are comparable and most accurate. These schemes have lower error since the
particles are fixed on a cartesian grid resulting in accurate computation of
divergence as discussed in \cref{sec:div_error}. The pressure convergence
flattens since it reaches the limit of accuracy possible with this value of
$c_o=80 m/s$ and further accuracy may be seen by increasing this further.

Clearly, the total force in case of E-C, L-RR-C, and EWCSPH scheme is zero
since we compute the acceleration on a uniform Cartesian grid of particles.
However, the total force in other schemes are accurate to order $10^{-5}$. The
times taken shows that L-PST-C is the highest since we apply the PST at
every timestep, the TV-C involves many terms in the equations and therefore
takes a lot of time. The E-C, and EWCSPH take the least time since they do not
use a PST~\footnote{These methods can be made even faster since the neighbors
  need not be updated, and the correction matrices can also be computed once
  and saved.}. The L-RR-C, L-IPST-C, and PE-IPST-C take a similar amount of
time.

\subsection{Comparison of conservation errors}
\label{sec:linear_conv}

Thus far, we have looked at the convergence of the various schemes. In this
section, we look at the schemes listed in \cref{tab:tg_compare_var} from
the perspective of conservation of linear and angular momentum. We solve
Gresho vortex and incompressible shear layer using all the schemes
discussed.

\subsubsection{The Grehso vortex}
\label{sec:grehso}

We consider the Gresho vortex problem \cite{gresho1990theory}, which is an
inviscid incompressible flow problem having the pressure and velocity fields
given by,
\begin{equation}
  p(r), u_\phi(r) =
  \begin{cases}
    12.5 r^2+5, 5r & 0 \leq r < 0.2,\\
    12.5r^2 - 20r + 4 \ln(5r) + 9, 2 - 5r & 0.2 \leq r < 0.4\\
    3 + 4 \ln(2), 0 & 0.4 \leq r\\
  \end{cases}
  \label{eq:gc_vortex}
\end{equation}

\begin{figure}[ht!]
  \centering
  \includegraphics[width=0.9\linewidth]{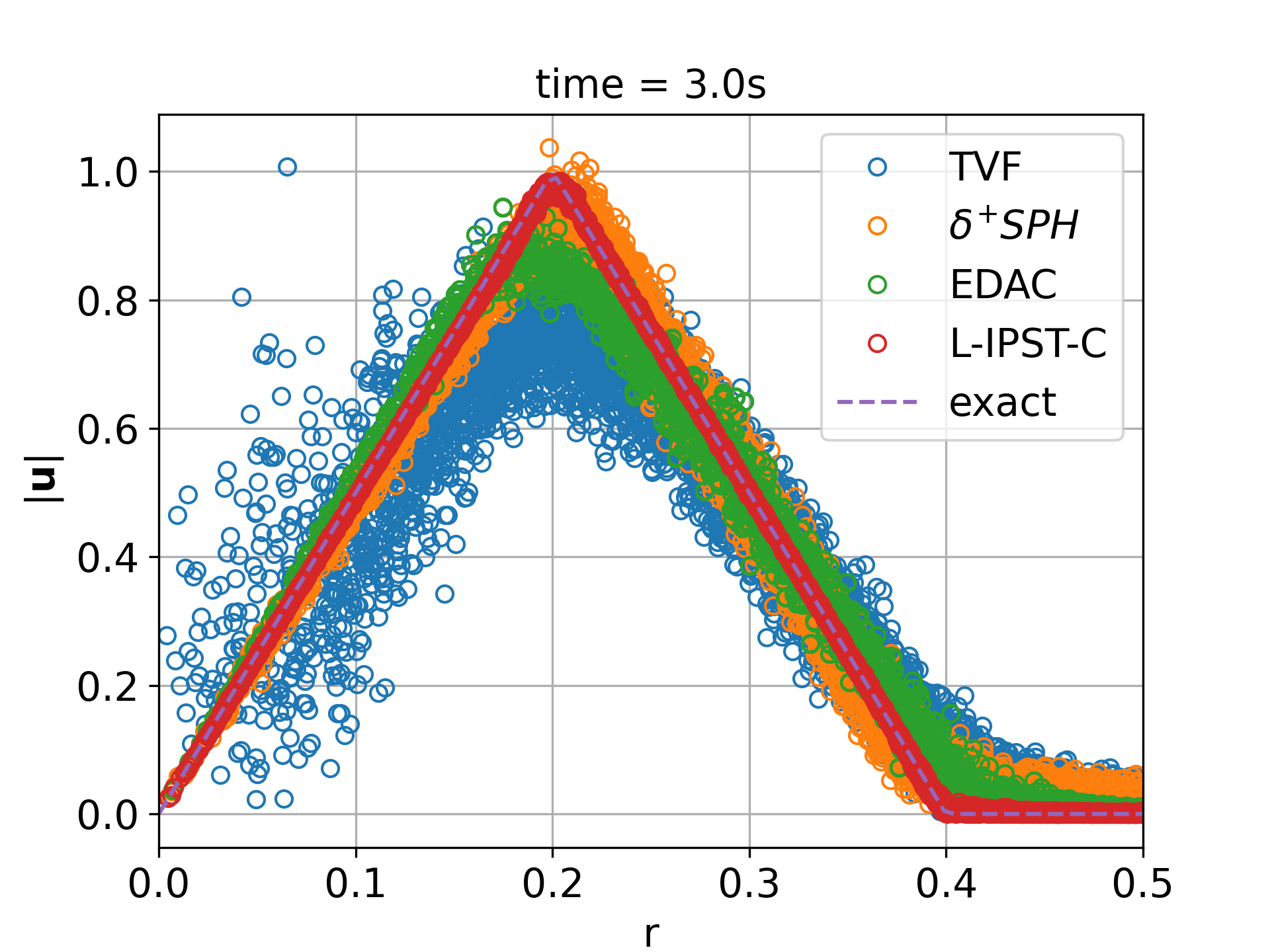}
  \includegraphics[width=0.9\linewidth]{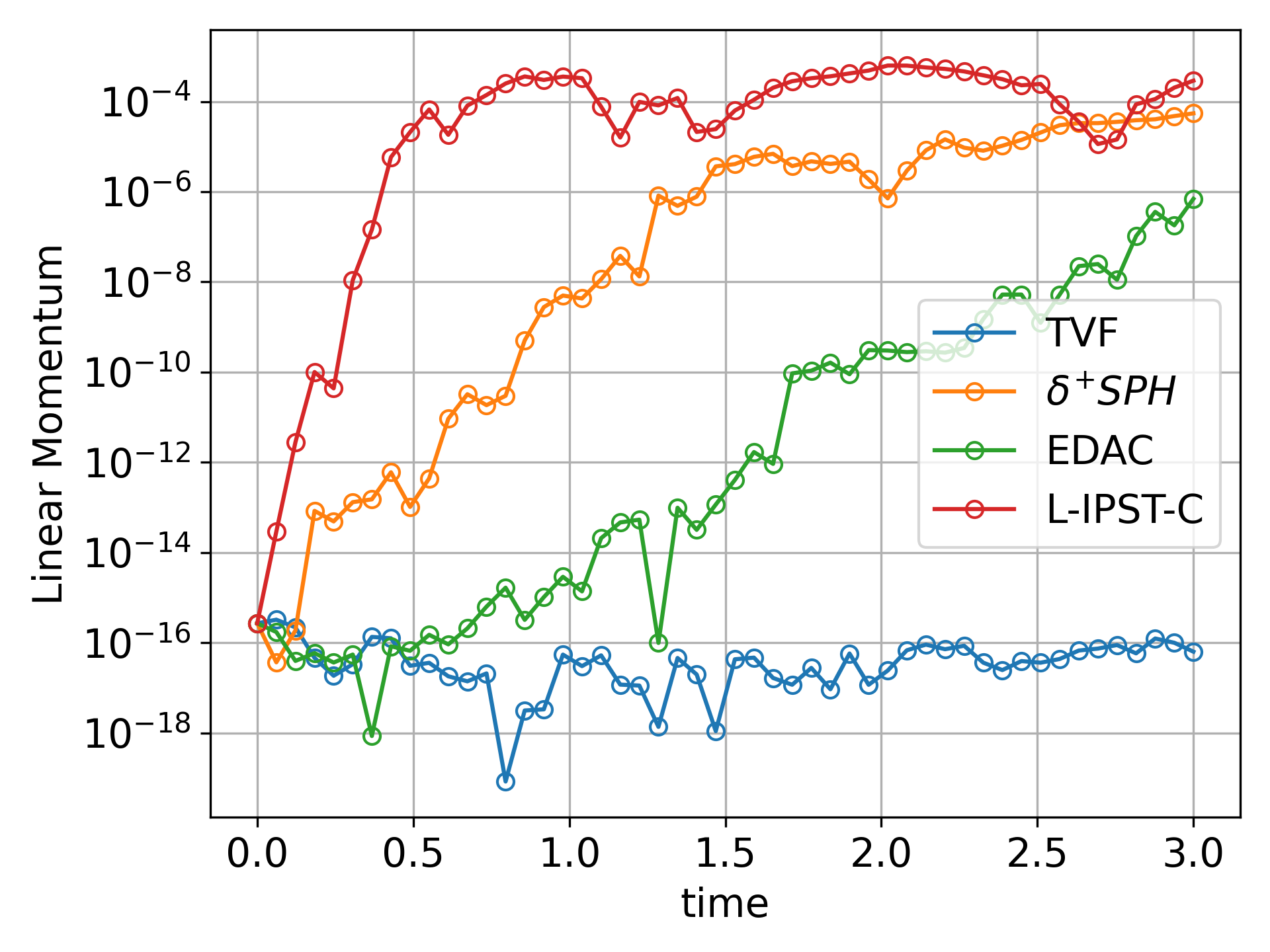}
  \caption{The velocity of particles with the distance from the center of the
    vortex (left) and the $x$-component of the total linear momentum (right)
    for all the schemes.}
  \label{fig:gc_vel_r}
\end{figure}

\begin{figure}[ht!]
  \centering
  \includegraphics[width=0.9\linewidth]{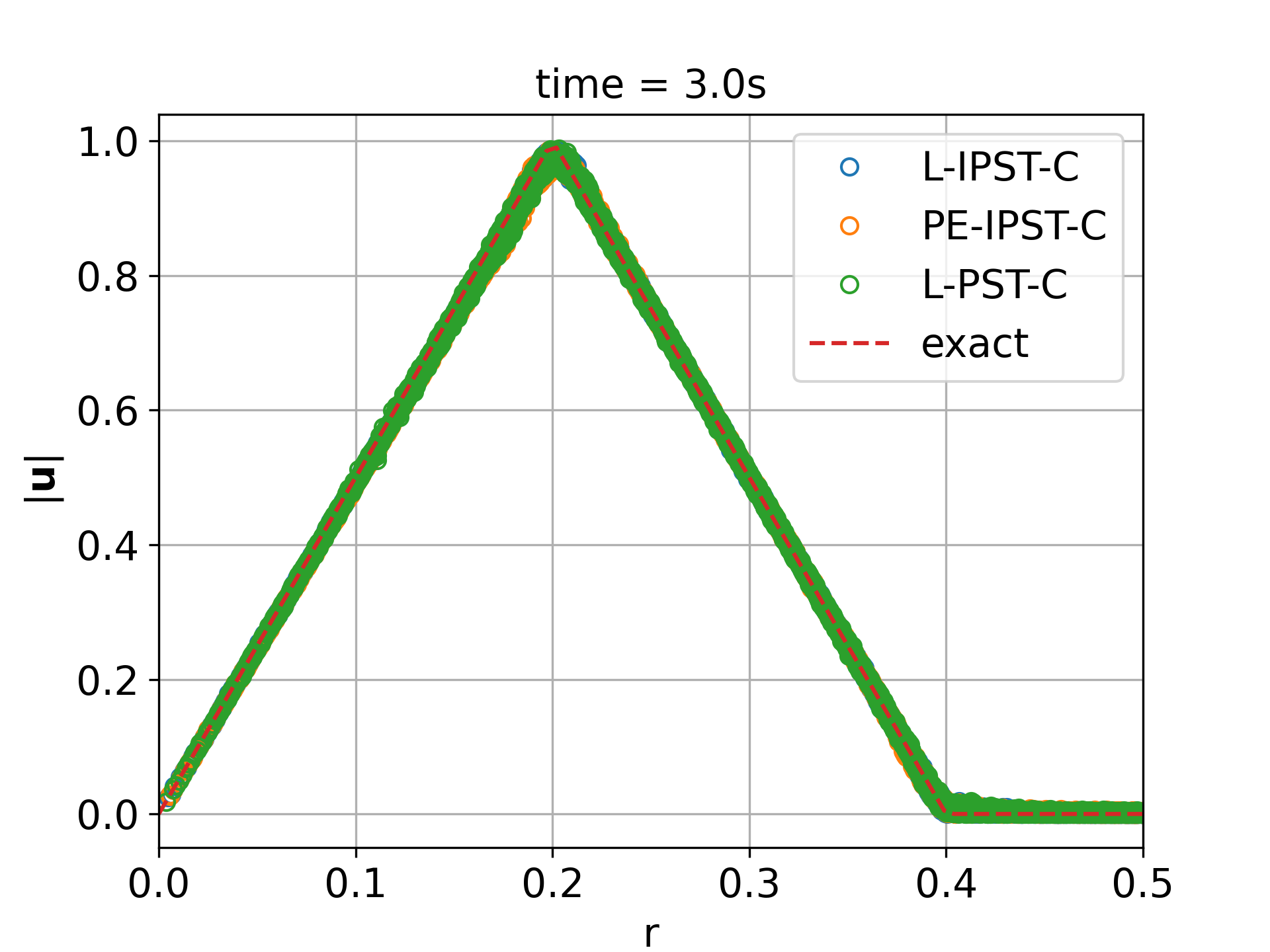}
  \includegraphics[width=0.9\linewidth]{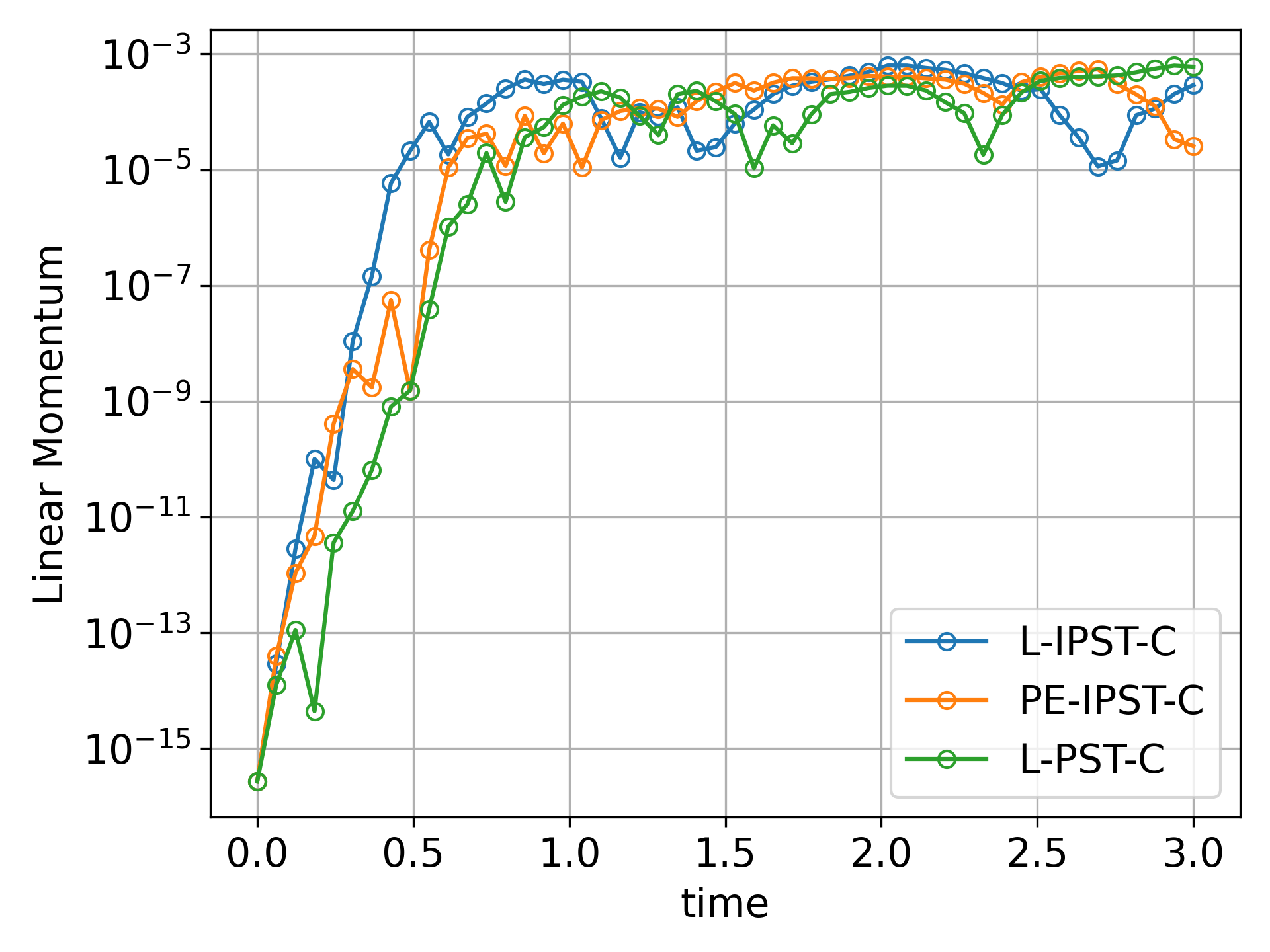}
  \caption{The velocity of particles with the distance from the center of the
    vortex (left) and the the $x$-component of the total linear momentum
    (right) for the variation of SOC scheme.}
  \label{fig:gc_scheme_var}
\end{figure}

We consider an unperturbed periodic domain of size $1\times1$ with the
center at $(0, 0)$. We set the kinematic viscosity, $\nu=0$, and the time
step and other properties as done in the Taylor-Green problem. The problem
is simulated until $t=3s$. Since the problem is inviscid, we expect the
scheme to retain the velocity and pressure field. We do not use artificial
viscosity in the simulations for any of the schemes. However, we use
density or pressure damping as given in \cref{eq:deltasph} or
\cref{eq:edac}, respectively in order to reduce the pressure oscillations.
Without this the solution becomes unstable in a short amount of time. We
perform the simulation of all the schemes listed in
\cref{tab:tg_compare_time} except the EWCSPH scheme \footnote{We discuss
the failed simulations are discussed in the \cref{apn:gc}.}. We note that
using an initial perturbed particle configuration results in very diffused
results for all schemes except the L-IPST-C.

In the \cref{fig:gc_vel_r}, we plot the velocity of the particles with the
distance, $r$ from the center (on left) and the $x$-component of the total
linear momentum with time for a $100\times100$ particle simulation. The L-IPST-C
scheme retains the velocity profile very well. The $\delta^+$SPH, EDAC, and
TVF schemes show diffusion due to inaccuracy in the pressure gradient
evaluation. Except for the TVF scheme, the rest show a finite increase in the
momentum bounded at $10^{-4}$. Clearly, approximate linear momentum conservation
is sufficient to obtain accurate results in the case of weakly compressible
flows.

We also perform the simulations with different versions of the SOC scheme listed
in \cref{tab:tg_compare_var}~\footnote{The L-RR-C, TV-C, and E-C scheme fail to
complete the simulation, and these are discussed in the \cref{apn:gc}}. In the
\cref{fig:gc_scheme_var}, we plot the velocity with the distance from the center
and the $x$-component of the linear momentum with time for $100\times100$
particle simulation. Clearly, all the schemes are accurate and approximately
conserve linear momentum as expected.

\begin{figure}[ht!]
  \centering
  \includegraphics[width=0.9\linewidth]{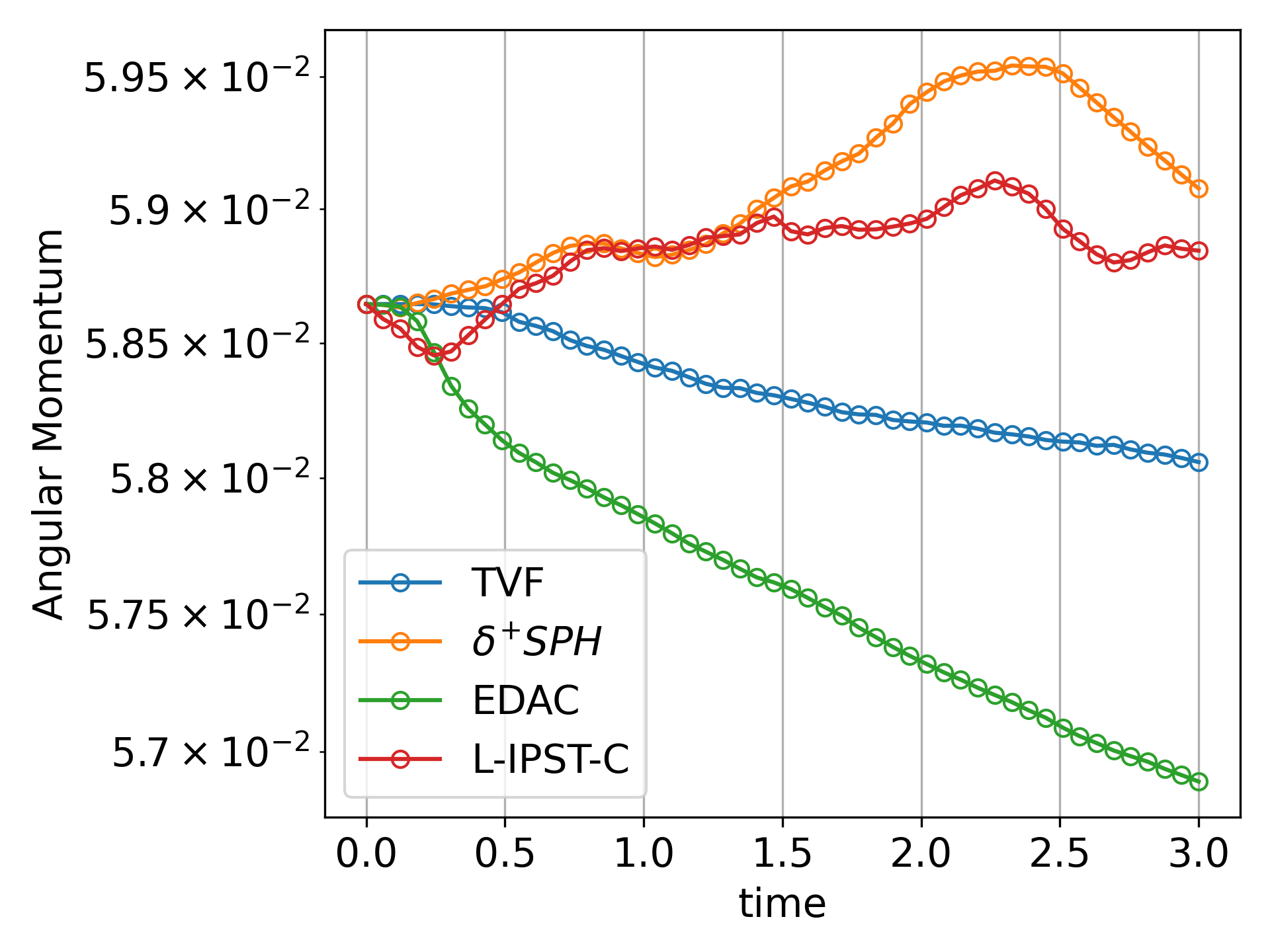}
  \includegraphics[width=0.9\linewidth]{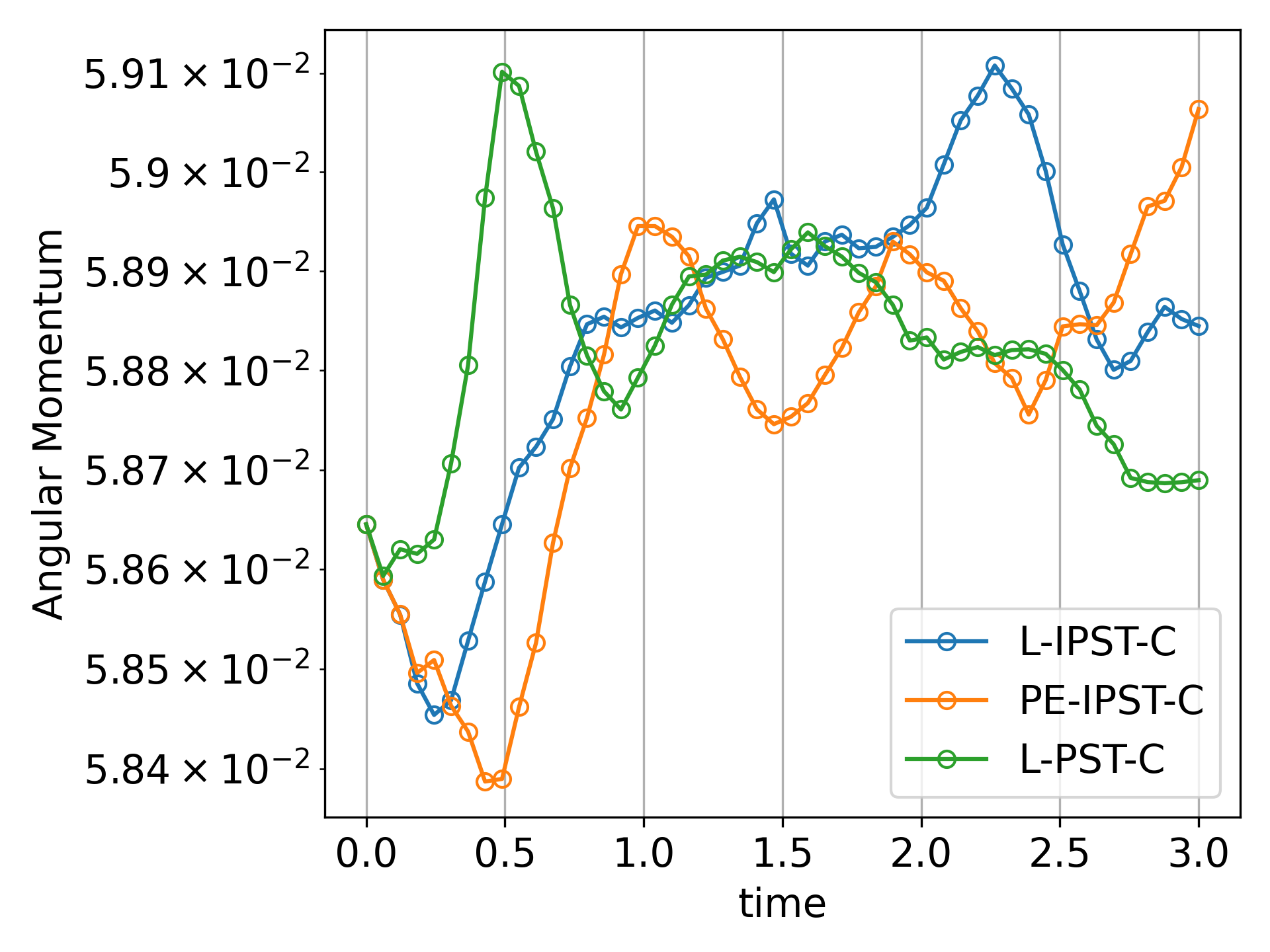}
  \caption{The angular momentum variation with time for Gresho-Chan vortex for different schemes.}
  \label{fig:gc_scheme_var_a_mom}
\end{figure}

In \cref{fig:gc_scheme_var_a_mom}, we show the angular momentum variation with
time for different schemes. None of the schemes conserve angular momentum, but
for the SOC schemes, the variations are very small and $O(5\times 10^{-4})$.

\subsubsection{The incompressible shear layer}
\label{sec:sl}

\begin{figure}[ht!]
  \centering
  \includegraphics[width=\linewidth]{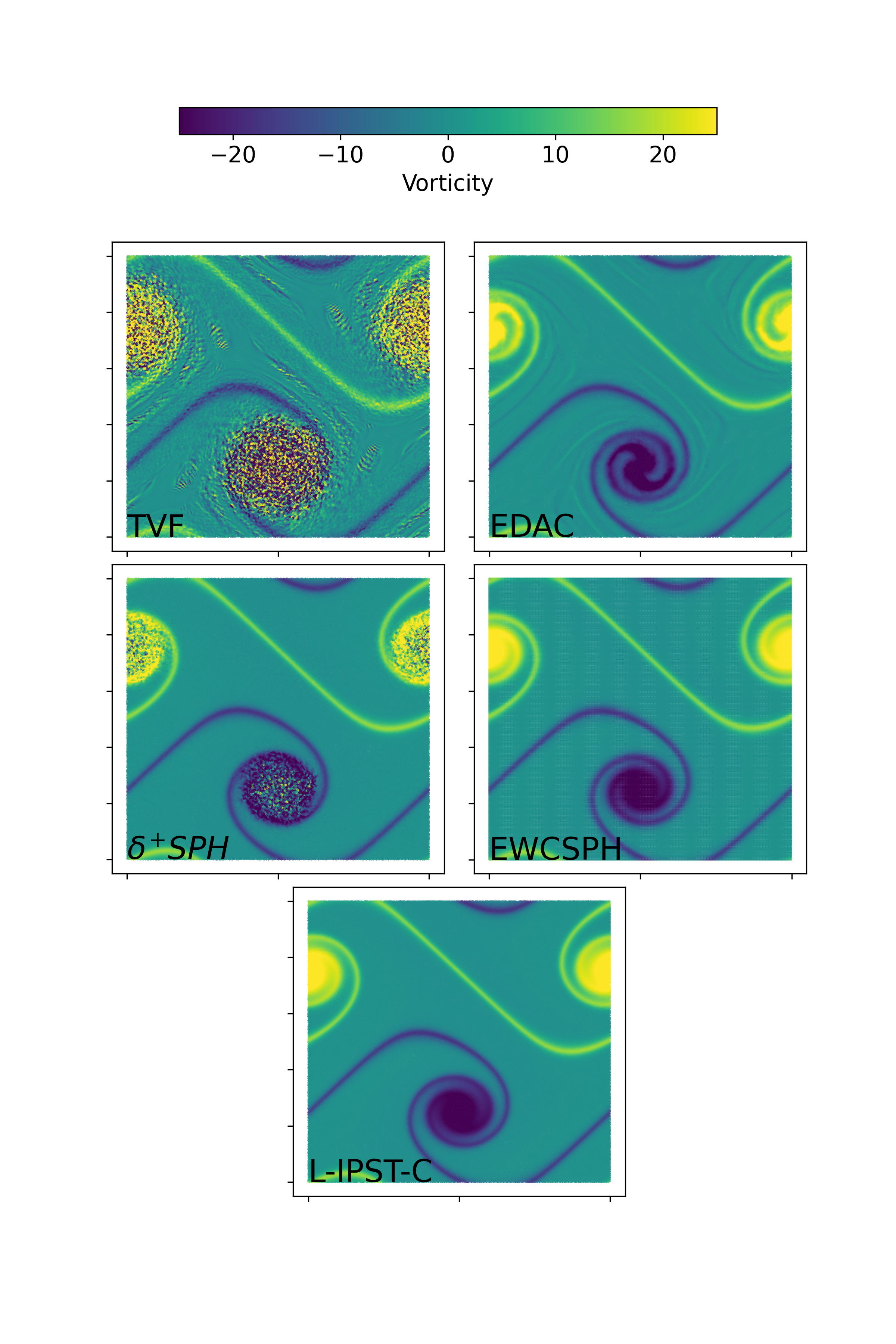}
  \caption{Vorticity contour plot for $500\times500$ resolution for all the schemes.}
  \label{fig:sl_vor}
\end{figure}

\begin{figure}[ht!]
  \centering
  \includegraphics[width=\linewidth]{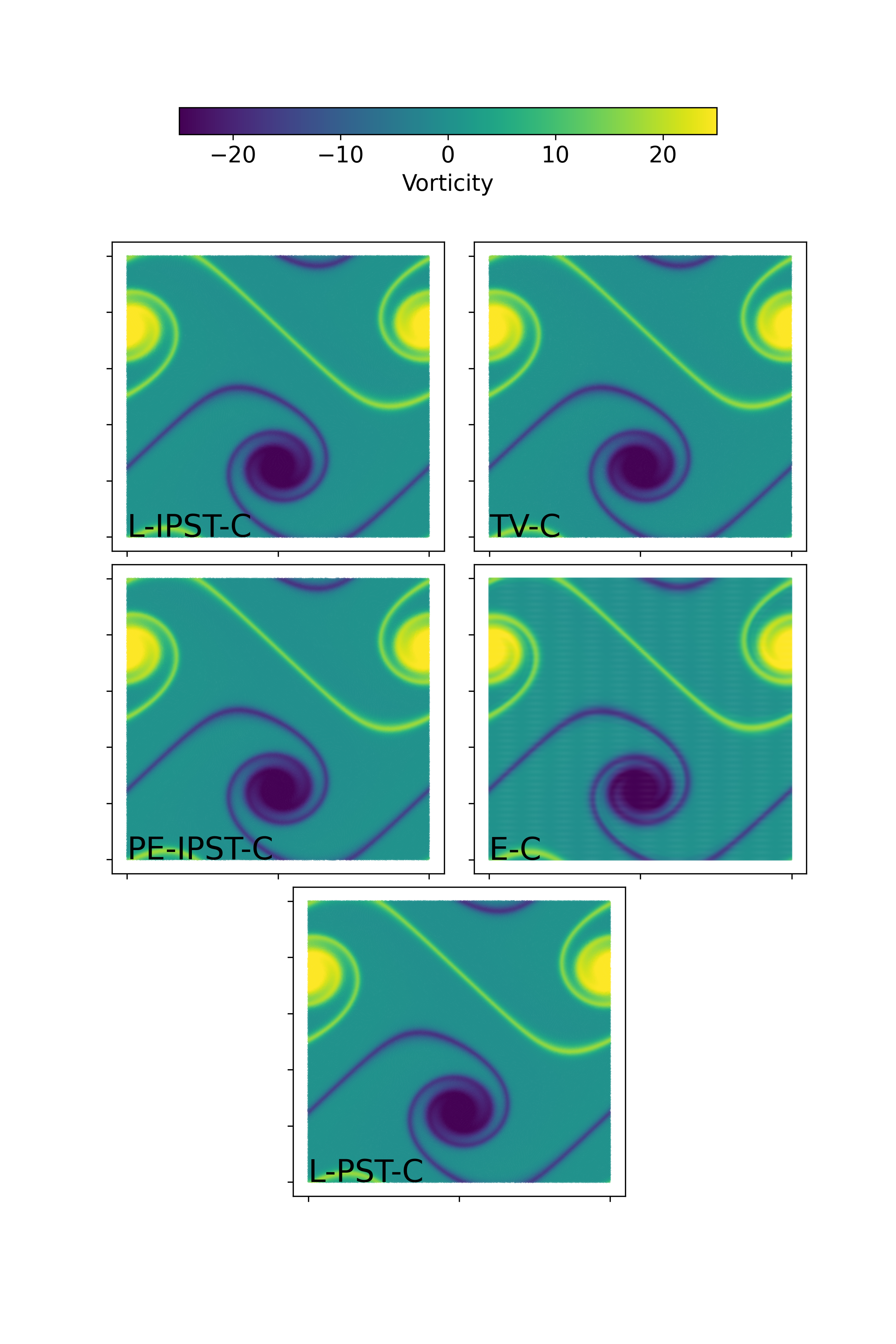}
  \caption{Vorticity contour plot for $500\times500$ resolution for all the schemes.}
  \label{fig:sl_vor_var}
\end{figure}

The incompressible shear layer simulates the Kelvin-Helmholtz instability in an
incompressible flow. This test case produces non-physical vortices for the
schemes where the operators are under resolved even when the scheme is
convergent \cite{diMovingMeshFinite2005}. The initial condition for the velocity
in x direction is given by
\begin{equation}
  u =
  \begin{cases}
    \tanh (\rho (y - 0.25))   & y \leq 0.5,\\
    \tanh (\rho (0.75 - y))   & y > 0.5,\\
  \end{cases}
\end{equation}
where $\rho=30$. In order to begin the instability, a small velocity is given in
y direction,
\begin{eqnarray}
  v = \delta sin(2 \pi x),
\end{eqnarray}
where $\delta=0.05$. We consider a small viscosity $\nu=1/10000$. We simulate
the problem using all the schemes listed in \cref{tab:tg_compare_time}. In
\cref{fig:sl_vor} and \cref{fig:sl_vor_var}, we plot the vorticity field for the
schemes~\footnote{The L-RR-C method failed to run due to discontinuity in the
initial velocity field.} discussed in this paper. We note that unlike the
inviscid problem of Gresho-Chan vortex, the scheme EWCSPH, TV-C and E-C shows
results matching other SOC schemes. In \cref{fig:sl_vor}, we observe that TVF
scheme and $\delta^+$SPH scheme show high frequency oscillations and while the
EDAC scheme is much better; However, it shows some undesired artifacts.

\subsection{Long time simulations}

In this section, we study the conservation for long time simulations using
the EDAC, TVF, and L-IPST-C schemes. We consider the Taylor-Green,
Gresho-Chan, and Poiseuille flow problem with the same condition as before.
We consider a UP particle distribution for all the schemes.

\begin{figure}[htbp]
  \centering
  \includegraphics[width=0.9\linewidth]{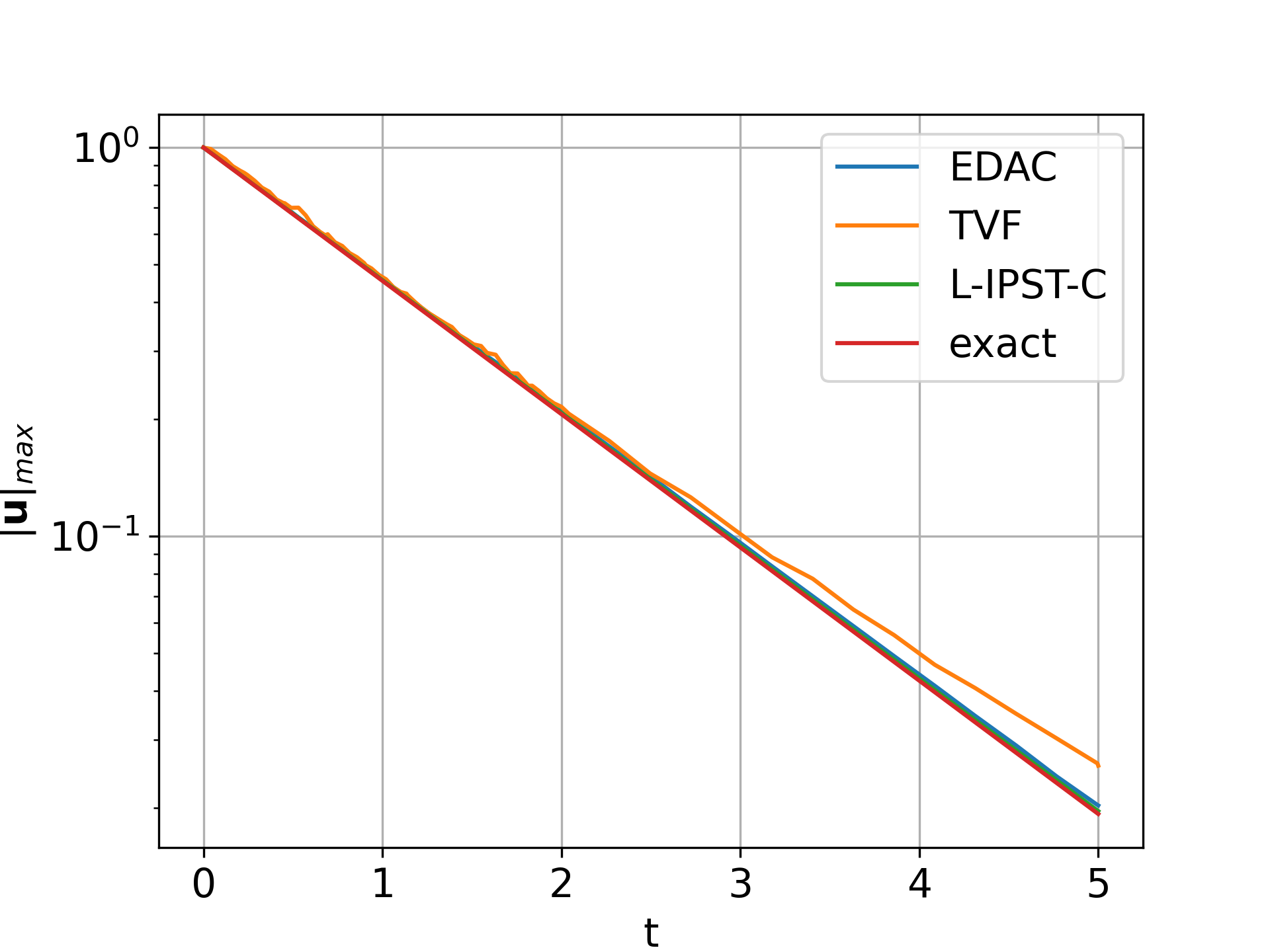}
  \includegraphics[width=0.9\linewidth]{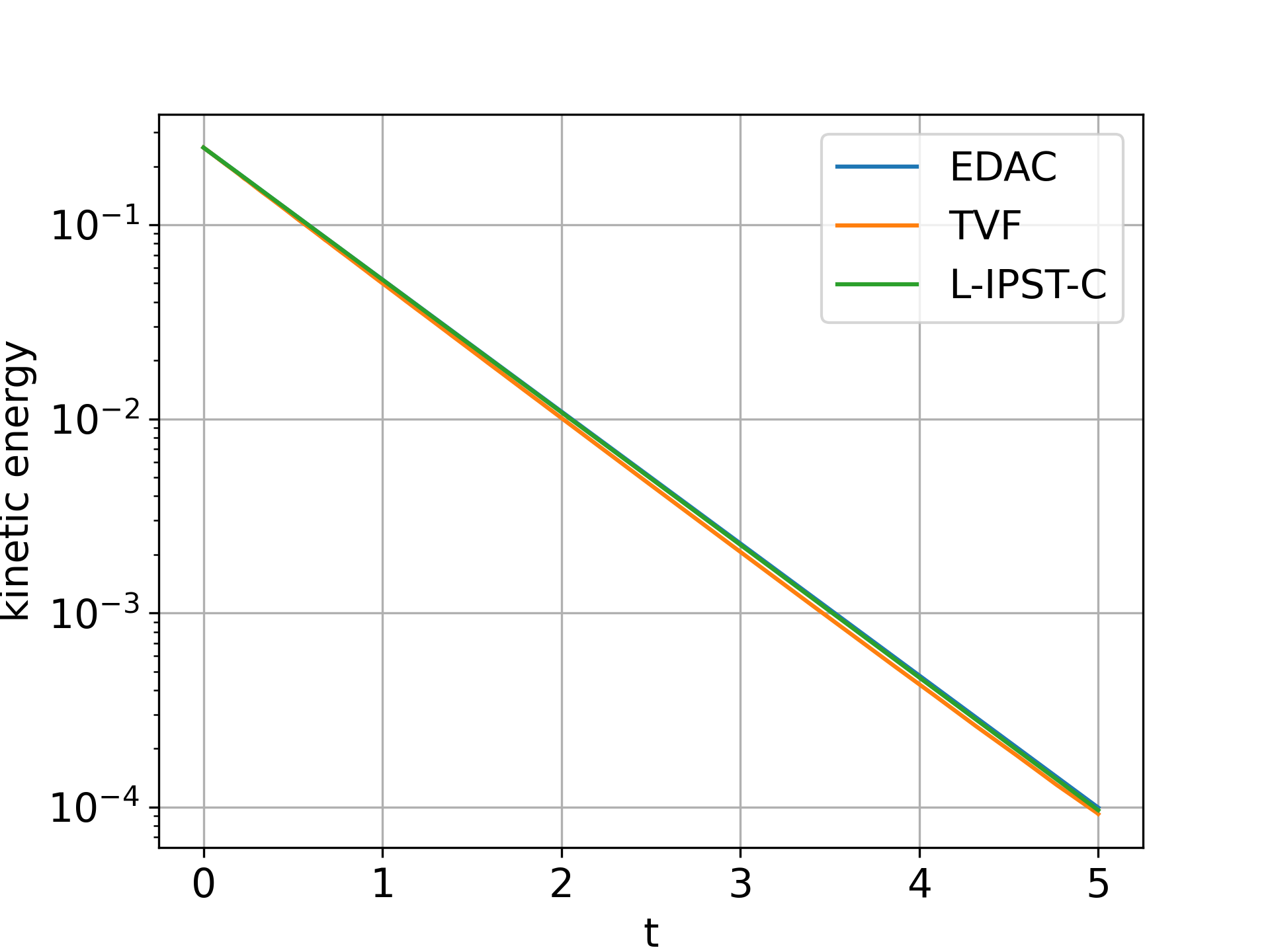}
  \caption{The maximum velocity decay and kinetic energy of the flow with
  respect time for Taylor-Green problem.}
  \label{fig:tg_lt}
\end{figure}

We simulate the Taylor-Green problem for $5 s$ at $Re=100$ compared to
the final time of $0.5 s$ in the previous simulations for all the
schemes. In \cref{fig:tg_lt}, we plot the velocity damping and the kinetic
energy of the flow as a function of time. The TVF scheme shows a significant
deviation from the exact result whereas the kinetic energy remains close
to other scheme solutions. We note that TVF scheme conserves linear
momentum exactly.

\begin{figure}[htbp]
  \centering
  \includegraphics[width=0.9\linewidth]{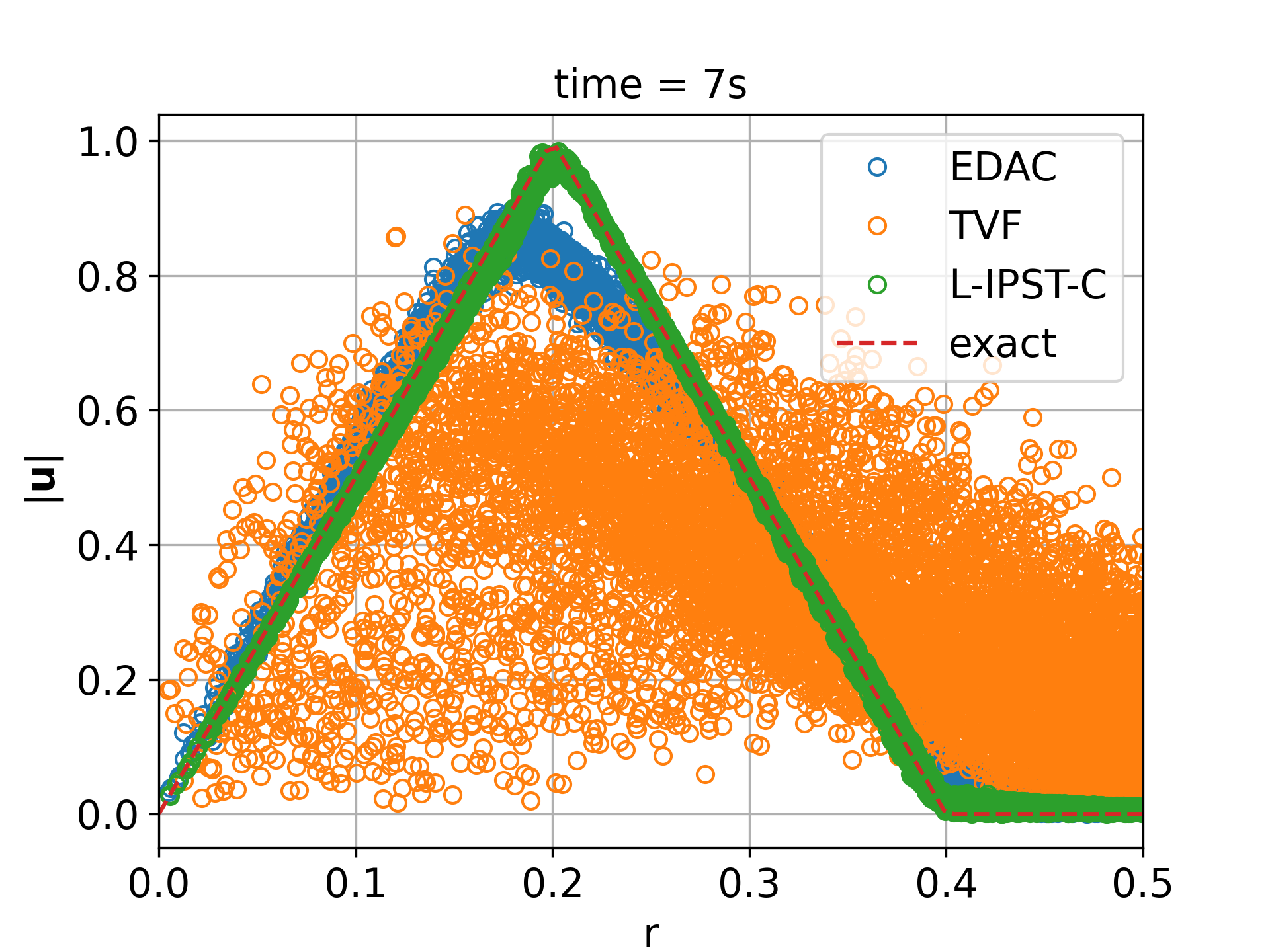}
  \includegraphics[width=0.9\linewidth]{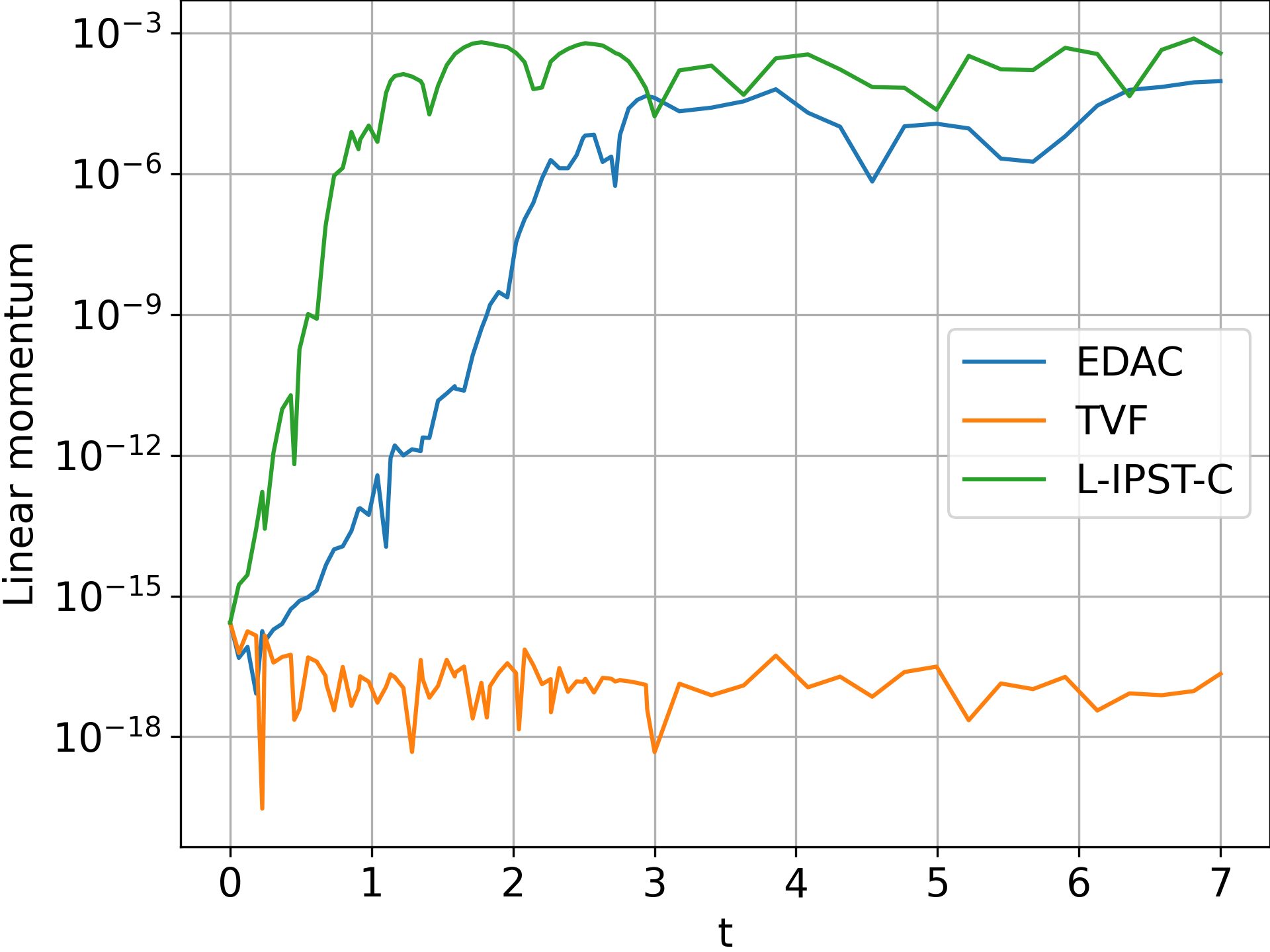}
  \includegraphics[width=0.9\linewidth]{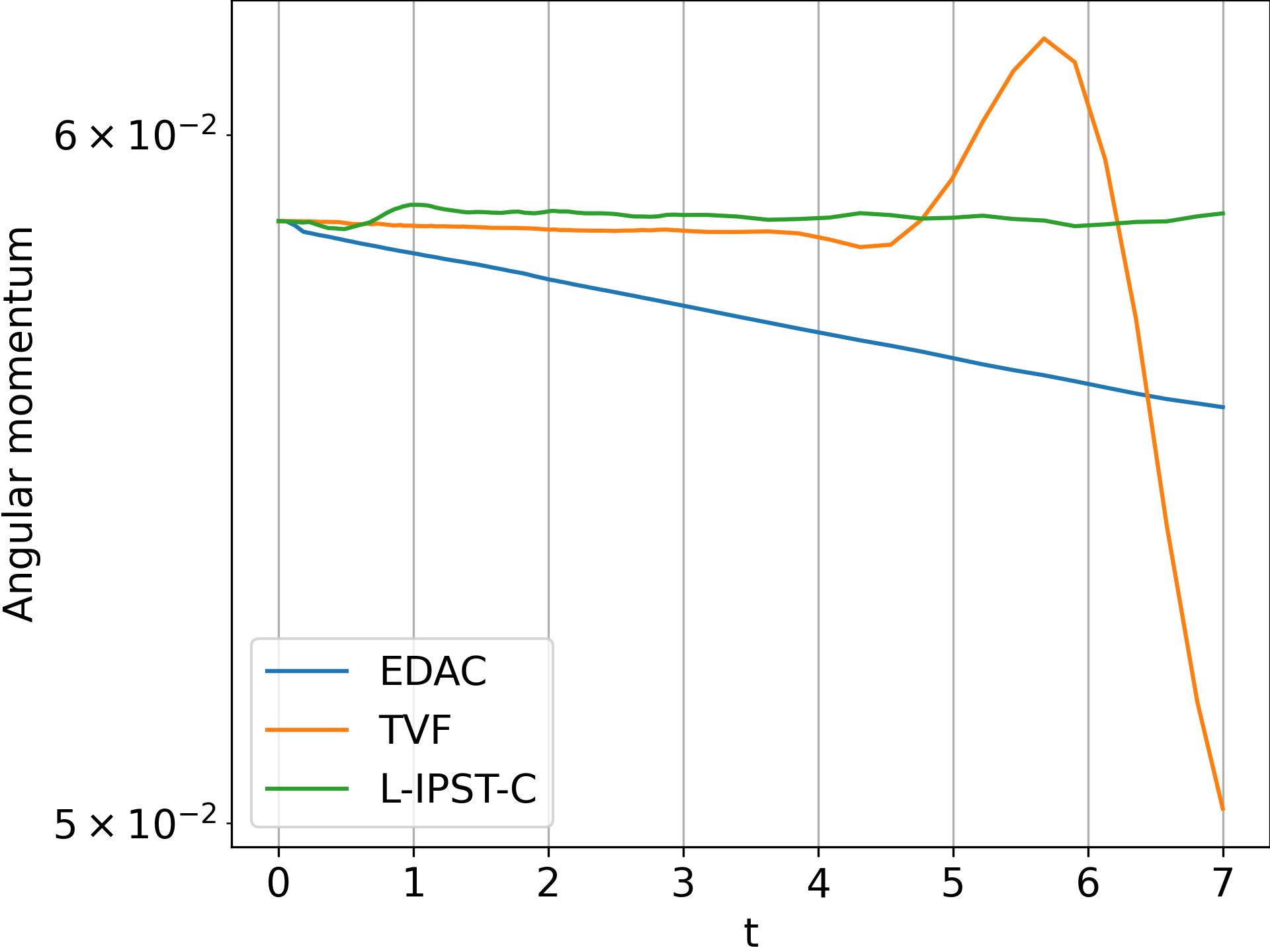}
  \caption{The velocity of particle with distance from the center, linear
  and angular momentum with respect to time for Gresho-Chan vortex
  problem.}
  \label{fig:gc_lt}
\end{figure}

We next simulate the Gresho-Chan vortex problem for $7 s$ compared to the
final time of $3 s$ in \cref{sec:grehso}. In \cref{fig:gc_lt}, we plot
the velocity as a function of $r$, and the linear and angular momentum with
respect to time for all the schemes. We observe that the TVF scheme does
not capture the physics of the problem however conserves linear momentum
but does not conserve angular momentum. In case of the EDAC scheme, the
physics is captured better. The linear momentum is not conserved, and the
solution looses angular momentum by a small amount, and the peak of the
velocity distribution is not captured accurately. The L-IPST-C schemes
retains the velocity field, and both the linear and angular momentum are
approximately conserved. After 7 seconds the L-IPST-C schemes is no longer
stable, and the velocity field is not captured accurately.

\begin{figure}[htbp]
  \centering
  \includegraphics[width=0.9\linewidth]{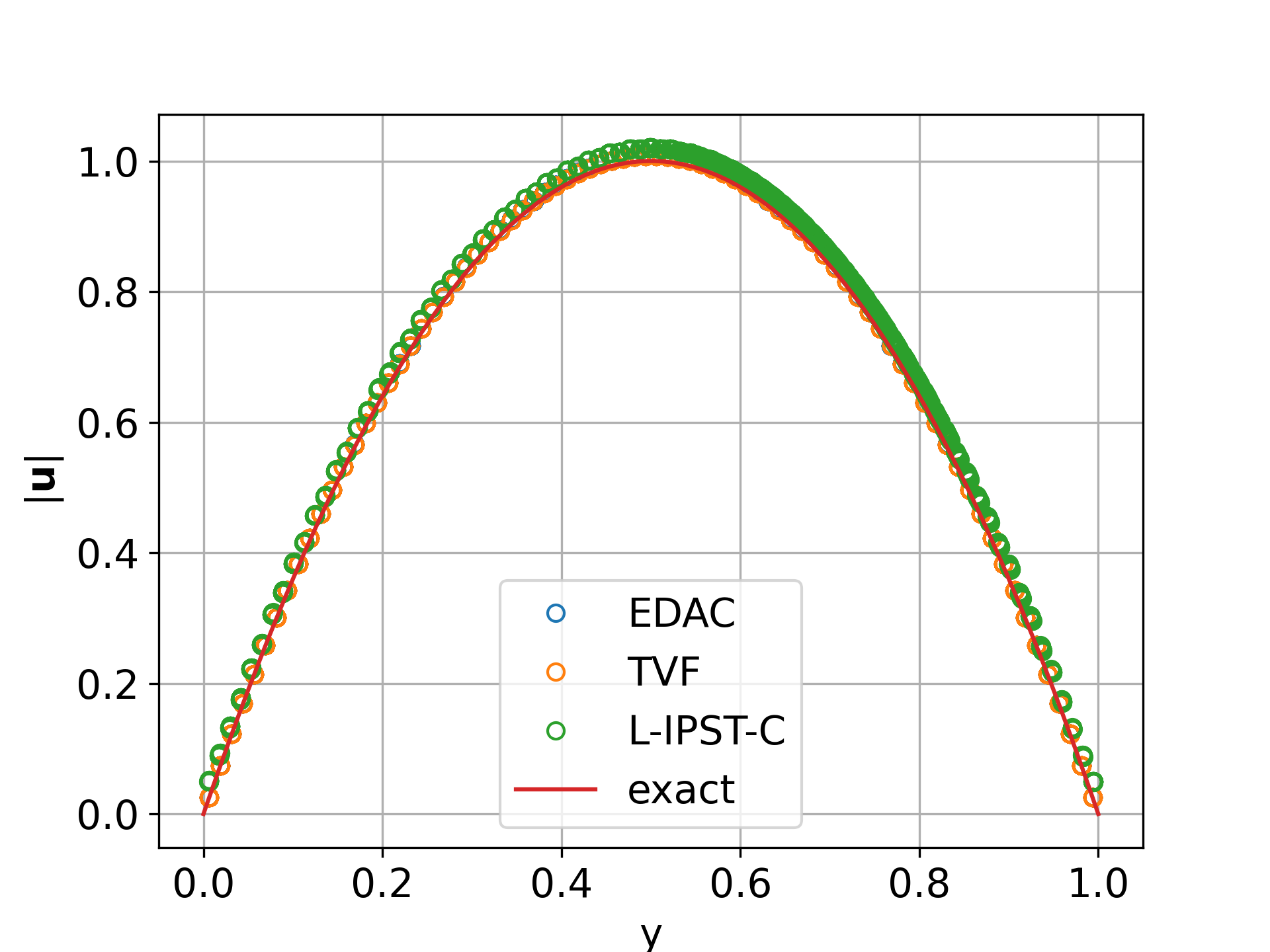}
  \includegraphics[width=0.9\linewidth]{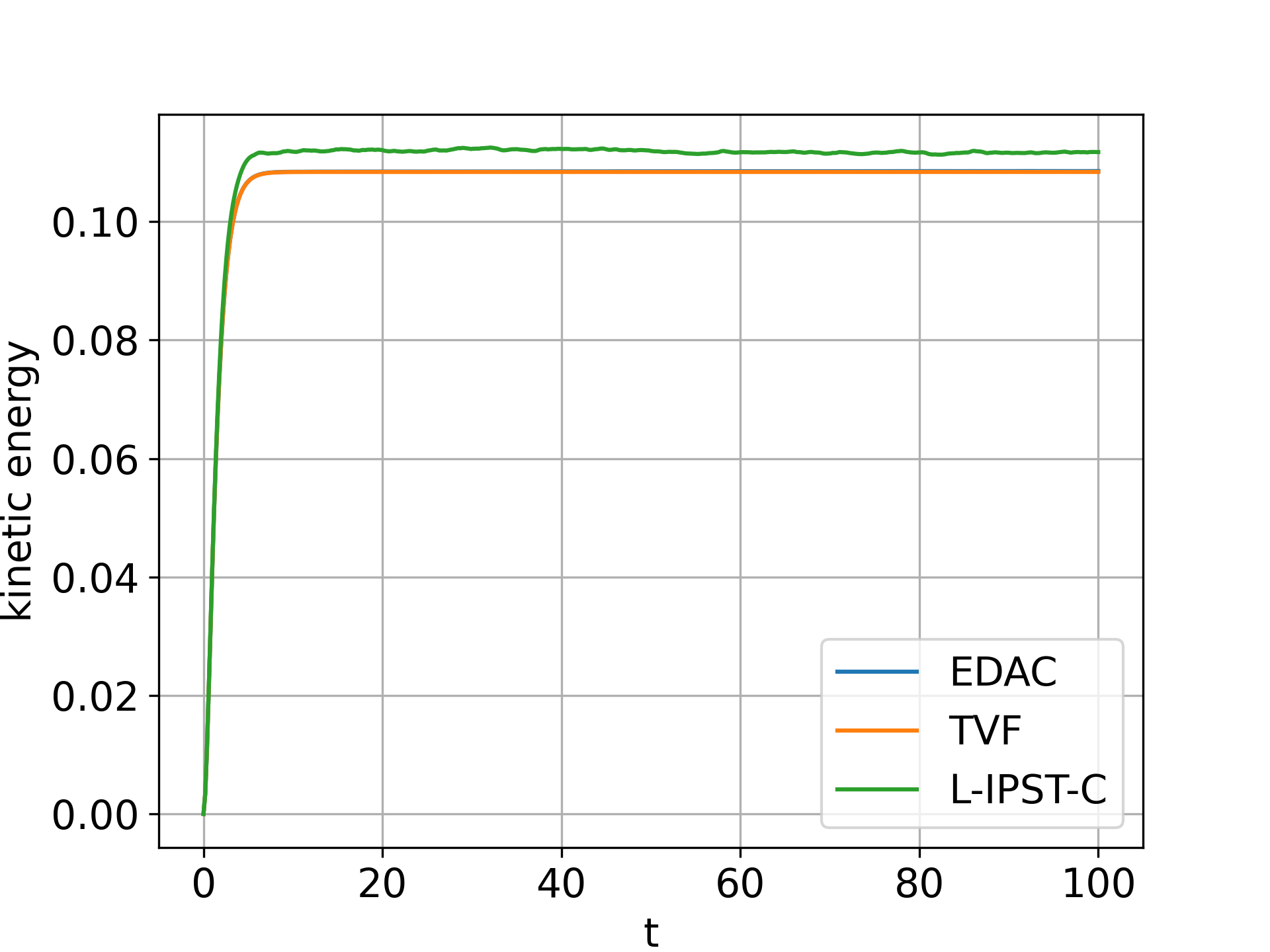}
  \caption{The $x$-component of the velocity across the plate and the kinetic
  energy of the flow with time for Poiseuille flow problem.}
  \label{fig:pois_lt}
\end{figure}

In the last test case, we simulate the Poiseuille flow problem for $100 s$
compared to $10 s$ in \cref{apn:solid_bc}. In the \cref{fig:pois_lt}, we plot the
$x$-component of the velocity and the kinetic energy of the flow with time.
Cleary, all the results are similar. In case of L-IPST-C a slight
deviation is observed near the wall due to approximation done near the wall
(See \cref{apn:solid_bc} for details).

These simulations suggests that even if a scheme is conservative like TVF
it may not produce accurate results. However, for a convergent scheme like
the L-IPST-C the results are accurate and despite there being no exact
conservation an approximate conservation is seen.

\subsection{Cost of computation}

\begin{figure*}[ht!]
  \centering
  \includegraphics[width=0.9\linewidth]{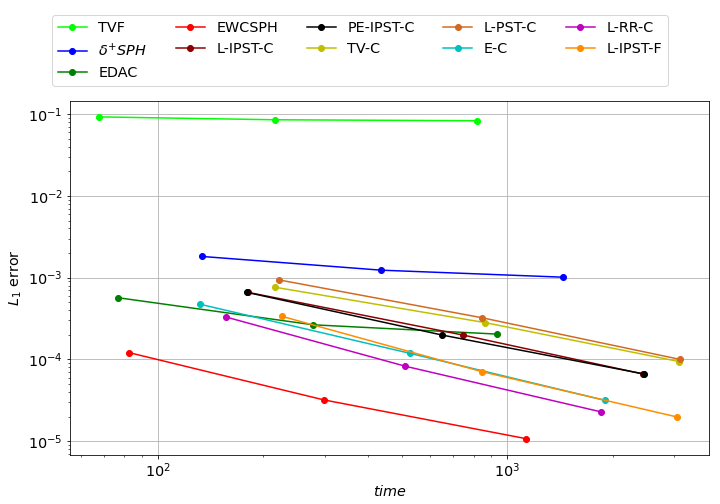}
  \caption{The $L_1$ error in velocity with respect to the time taken to
    evaluate 5000 timesteps for all the schemes discussed in the previous
    sections.}
  \label{fig:cost_comp}
\end{figure*}

In this section, we compare the cost of computation of all the schemes
considered in this study. We simulate the Taylor-Green problem for $5000$
timesteps with $50$, $100$, and $200$ resolutions for all the schemes. We
use an Intel(R) Xeon(R) CPU E5-2650 v3 processor and execute all the
simulations in serial. In \cref{fig:cost_comp}, we plot the $L_1$ error in
velocity computed using the \cref{eq:tg_l1} as a function of time taken for
the simulation. Clearly, all the SOC schemes are close to each other in
terms of errors. The E-C and EWCSPH scheme takes very less time and are
very accurate; however, EWCSPH is not convergent in pressure as shown in
\cref{sec:comp_other_scheme}. The EDAC scheme has lower error comparable to
the SOC schemes; however, its convergence rate reduces with increase in
resolution. We show that despite having higher time taken by the SOC
schemes, they achieve higher accuracy with a fewer number of particle. For
some schemes, these accuracy levels are not achievable at all.

\section{Conclusions}
\label{sec:disc}

In this paper, we have performed a numerical study of the accuracy and
convergence of a variety of SPH schemes in the context of weakly-compressible
fluids. Based on the numerical study performed in the previous sections, we
summarize the key findings below.

\subsection{Choice of smoothing kernel}

We first considered the SPH approximation of a function and its derivative
using different kernels. All the kernels considered here show second-order
convergence when the support radius is suitably chosen. The accuracy is
marginally effected by the change in type of kernel. The smoothing error of an
SPH approximation scales as $O(h^2)$ and this necessitates that the smoothing
length of the kernel be as small as possible. This implies that $h_{\Delta s}$
be small. As is well known, the discretization errors scale as
$O((\frac{\Delta s}{h})^{\beta + 2})$ and this necessitates that the smoothing
radius be larger. These two requirements are contradictory. We find that by
using a modest $h = 1.2 \Delta s$, along with the kernel corrections of
\citet{bonet_lok:cmame:1999} or \citet{liu_restoring_2006} we are able to
obtain close to second-order convergence for the kernels considered in this
work. It holds up to a resolution of $L/\Delta s = 500$, where $l=1m$ which
appears to be among the highest resolutions we have seen in the literature
concerning the convergence of SPH methods. In the literature, we find kernels
like the cubic spline to demonstrate pairing
instabilities~\cite{dehnen-aly-paring-instability-mnras-2012}. We can avoid
this instability by using a particle shifting technique (PST).

\subsection{Choice of suitable operators}

The SPH approximation of operators like the gradient, divergence, and
Laplacian must be chosen carefully. In this paper, we recommend two methods
for gradient approximation and three methods for viscous term approximation
that ensure second-order convergence. The approximations which ensure
pair-wise linear momentum conservation are always divergent. In the future,
one could explore pair-wise linear momentum conserving and second-order
convergent SPH approximation in a perturbed domain using SPH. Furthermore, the
widely used artificial compressibility assumption makes the scheme $O(M^2)$
accurate. We recommend using high $c_o$ values or a dual-time stepping
criteria to achieve convergence. Solving the pressure using the pressure
Poisson equation may also provide SOC, although those have not been studied in
this work.

\subsection{Particle density and fluid density}

We recommend that one employ the fluid density in the governing differential
equation as a property that convects with the particle. The approximation of
the SPH operators should not be a function of a property of the fluid i.e.\
density. We obtain the integration volume by \cref{eq:num_vol} where the
mass and kernel support radius of particles is kept constant. In the future,
it would be important to explore the convergence of SPH operators when the
mass as well as the support radius of the particles are varying as required by
an adaptive SPH algorithm.

\subsection{SOC scheme and variations}

We demonstrate Eulerian as well as Lagrangian SPH schemes that are second-order
convergent. We show that the Eulerian schemes captures the divergence accurately
due to symmetry in the particle distribution resulting in better accuracy in
pressure. We derive a pressure evolution equation using the continuity equation
that resembles the EDAC SPH scheme in literature. We show that the PST step in
the Lagrangian method can be replaced by a remeshing step which is another
moment conserving regularization. However, remeshing is not stable in the
presence of jumps in the properties as observed in the case of the Gresho vortex
(see \cref{apn:gc}) and incompressible shear layer. The PST step can be included
in the momentum equation resulting in the $\delta^+$SPH scheme. From the
$\delta^+$SPH method, one can obtain the Eulerian form of WCSPH method by
setting the shifting velocity to $-\ten{u}$. All these schemes are SOC when we
use a second-order convergent approximation for the operators. We show that even
though the schemes are non-conservative in the absolute sense, approximate
conservation also produces accurate results in the case of incompressible flows.

Thus, by a judicious choice of discretization, particle shifting, and a
separation of the fluid and particle densities we have shown that second-order
convergence is possible using the SPH method for weakly-compressible flows. We
do observe that the SPH discretization of the divergence operator introduces
errors for divergence-free fields which are noticeably absent in the case of
an Eulerian method due to the symmetry of the particle distribution. This
introduces significant errors into the pressure; it would be valuable to
develop more accurate divergence operators for the Lagrangian case.

Given that the proposed schemes are second-order, it would be important to
study the boundary conditions employed in the SPH to see how they affect
the accuracy and order of convergence. A preliminary analysis performed in
\cref{apn:solid_bc} suggests that a popular solid-wall boundary
condition~\cite{Adami2012,maciaTheoreticalAnalysisNoSlip2011} is not second
order convergent. The accuracy of the boundary conditions will be
investigated in the future.

A similar analysis in the context of variable smoothing length, and mass
would be very useful in light of many recent developments of adaptive SPH
methods~\cite{mutaadaptive2021} One concern of note is the increased
computational effort required to maintain second-order convergence and
future developments in this area would be important for practical
simulation using the SPH method.

\appendix

\section{Comparison of kernels}
\label{apn:kern}

\begin{table}[ht!]
  \centering
  \renewcommand{\arraystretch}{2}
  \begin{tabular}{l|l|l|l}
    \text{Name} & \text{Radius} & \text{$\beta$} & \text{Remark}\\
    \hline
    \text{$G$ - Gaussian \cite{monaghan-review:2005}} & 3 & 0 & \text{Truncated for low $N_{nbr}$}\\
    \text{$QS$ - Quintic spline \cite{edac-sph:cf:2019}} & 3 & 3 & \text{Tensile instability}\\
    \text{$CS$ - Cubic spline \cite{monaghan-review:2005}} & 2 & 5 & \text{Paring and tensile instability}\\
    \text{$WQ_2$ - Wendland $O(2)$ \cite{wendland_piecewise_1995}}& 2 & 5 & \text{No tensile or pairing instability}\\
    \text{$WQ_4$ - Wendland $O(4)$ \cite{wendland_piecewise_1995}} & 2 & 8 & \text{Produces higher accuracy}\\
    \text{$WQ_6$ - Wendland $O(6)$ \cite{wendland_piecewise_1995}} & 2 & 11 & \text{Produces higher accuracy}\\
\end{tabular}
  \renewcommand{\arraystretch}{1}
  \caption{Kernels and their properties}
  \label{tab:kernels}
\end{table}

We consider the set of kernels listed in \cref{tab:kernels}. It covers a
wide range (high order, kernels having tensile instability and pairing
instability
\cite{dehnen-aly-paring-instability-mnras-2012,sph:tensile-instab:monaghan:jcp2000}).
In order to assess the effect of $h_{\Delta s}$ for a kernel, we perform
the numerical experiment proposed by
\citet{dehnen-aly-paring-instability-mnras-2012}. We evaluate particle
density using \cref{eq:num_den} for increasing the number of neighbors
$N_{nbr}$, for each of the kernels. The increase in $N_{nbr}$ corresponds
to the scaling of the smoothing kernel using the $h_{\Delta s}$ parameter.
In this numerical experiment, we change both the resolution and $h_{\Delta
s}$.

\begin{figure}[ht!]
  \centering
  \includegraphics[width=0.9\linewidth]{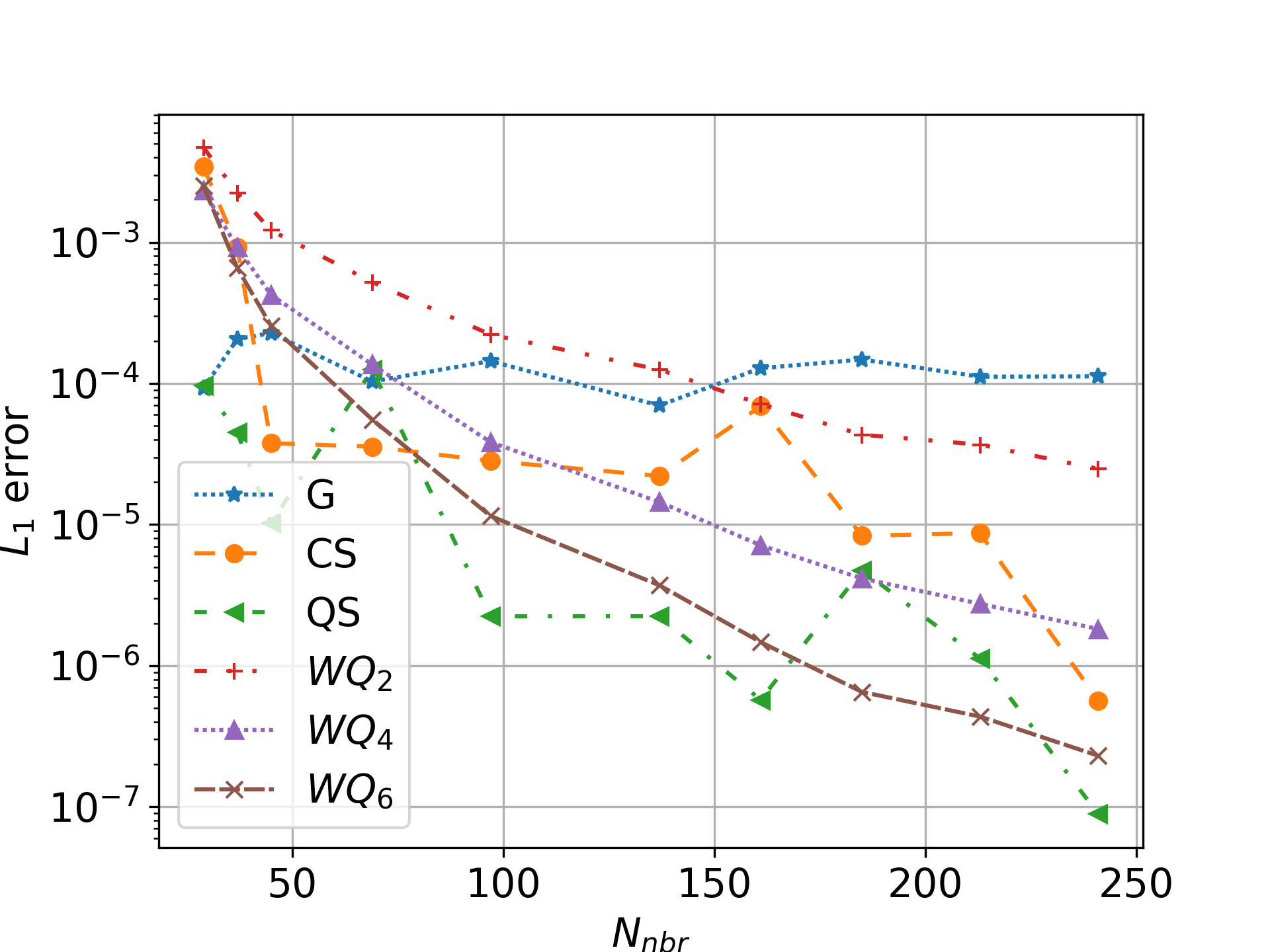}
  \caption{The particle density for different kernels with varying number of neighbors}
  \label{fig:rho_c}
\end{figure}

In the \cref{fig:rho_c}, we plot the absolute error in the particle density of
one particle in an UP domain for different kernels with the change in $N_{nbr}$
under the kernel support. Clearly, the Wendland class of kernel shows a
monotonic decrease in error with the increasing $N_{nbr}$. However, in the case
of the $G$ and $QS$ kernels, the errors are an order less at a lower $N_{nbr}$
compared to Wendland class of kernels. The error in the $G$ kernel does not
change significantly with the change in the $N_{nbr}$ compared to others. It
is because we truncate the $G$ kernel to have compact support. In the $QS$, the error is lower than the $WQ_4$ in the entire plot.
Therefore, we drop $WQ_2$ and $WQ_4$ in the subsequent investigations since it
reaches the order of accuracy of $QS$ when $N_{nbr}$ is approximately $60$. High
$N_{nbr}$ results in higher computational cost.

\begin{figure}[ht!]
  \centering
  \includegraphics[width=\linewidth]{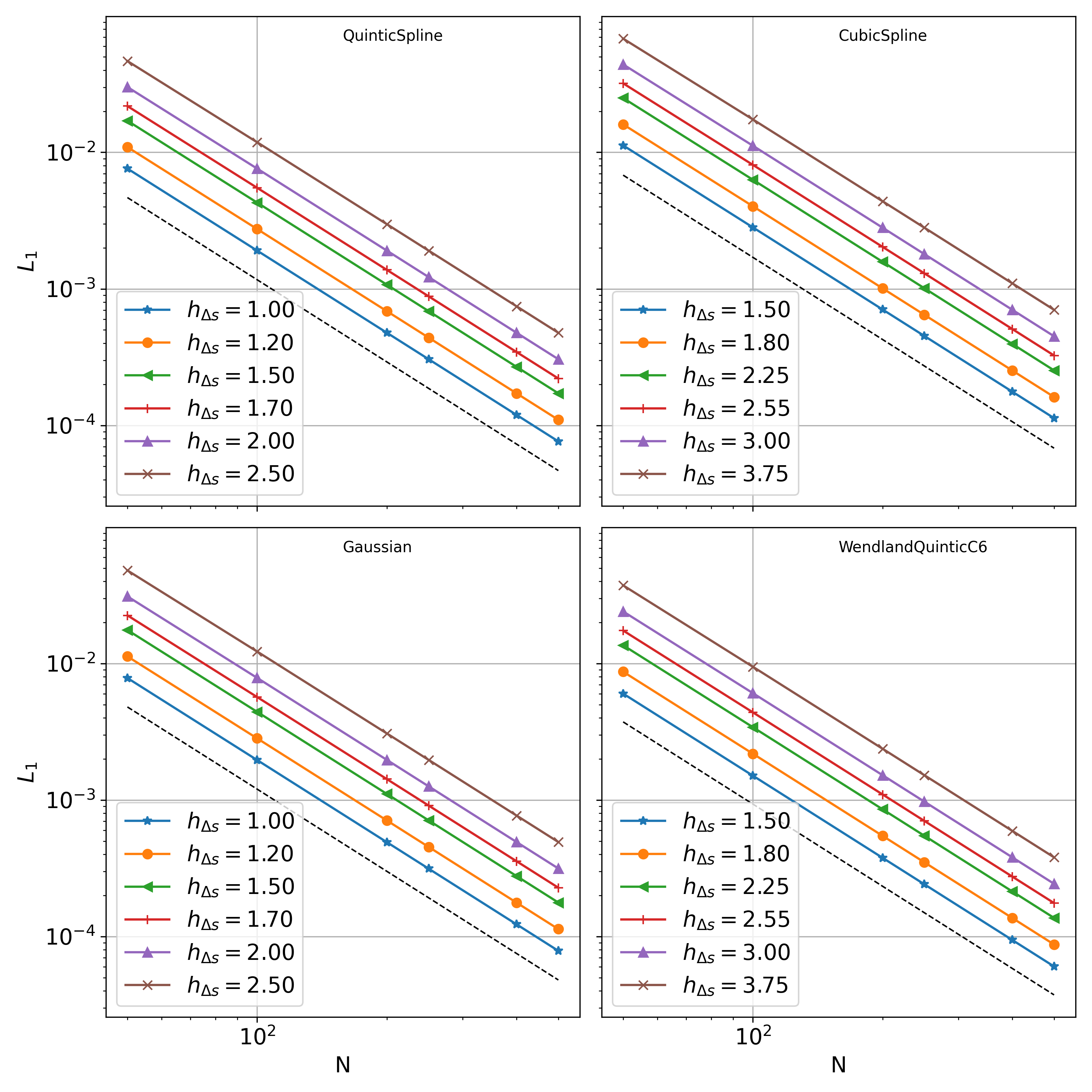}
  \caption{Order of convergence for the approximation of a function for
    different $h_{\Delta s}$ values in an UP domain. The dashed line shows the
    second order rate.}
  \label{fig:d0p0_f}
\end{figure}

\begin{figure}[ht!]
  \centering
  \includegraphics[width=\linewidth]{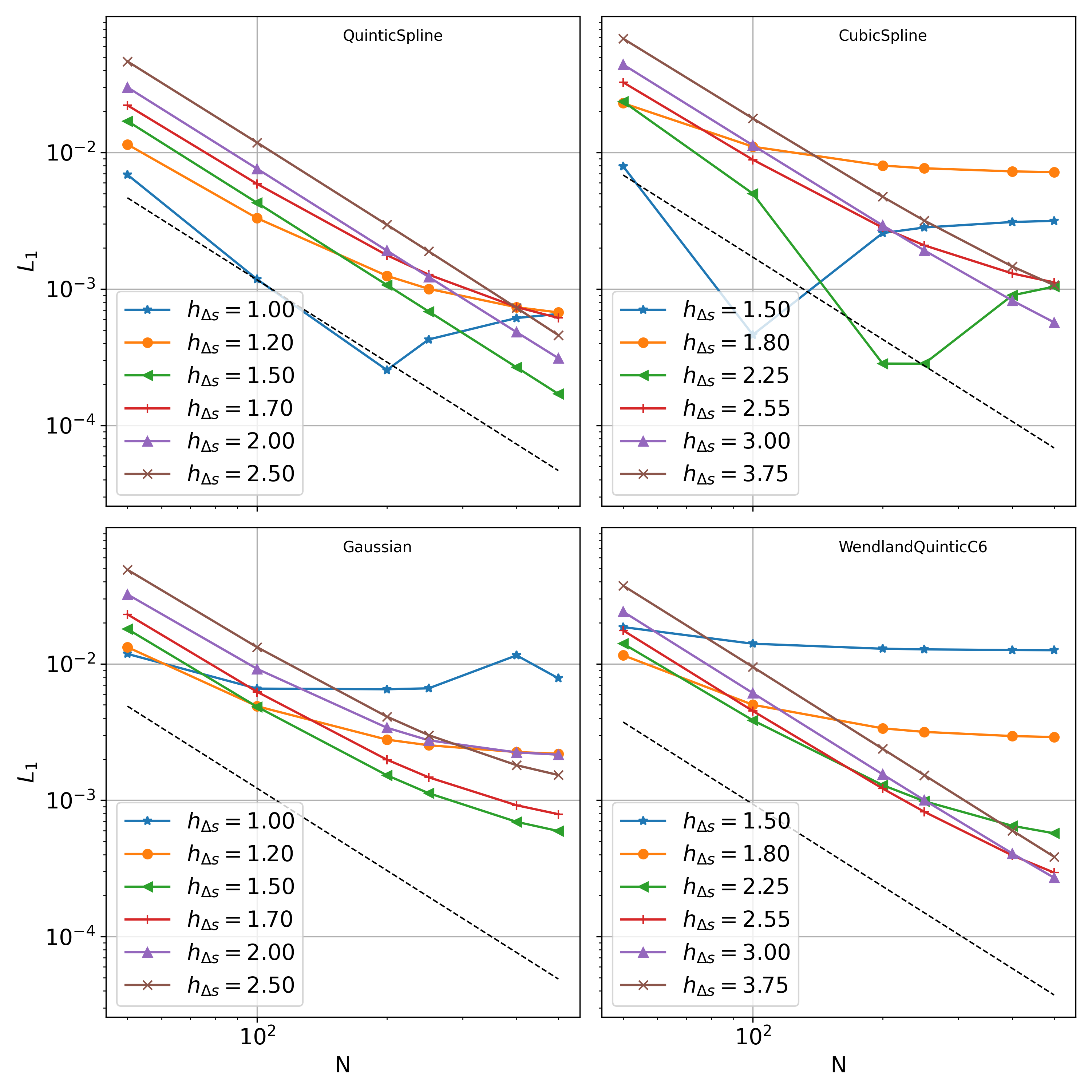}
  \caption{Order of convergence for derivative approximation for different
    $h_{\Delta s}$ values in an UP domain. The dashed line shows the second
    order rate.}
  \label{fig:d1p0_df}
\end{figure}

\begin{figure}[ht!]
  \centering
  \includegraphics[width=\linewidth]{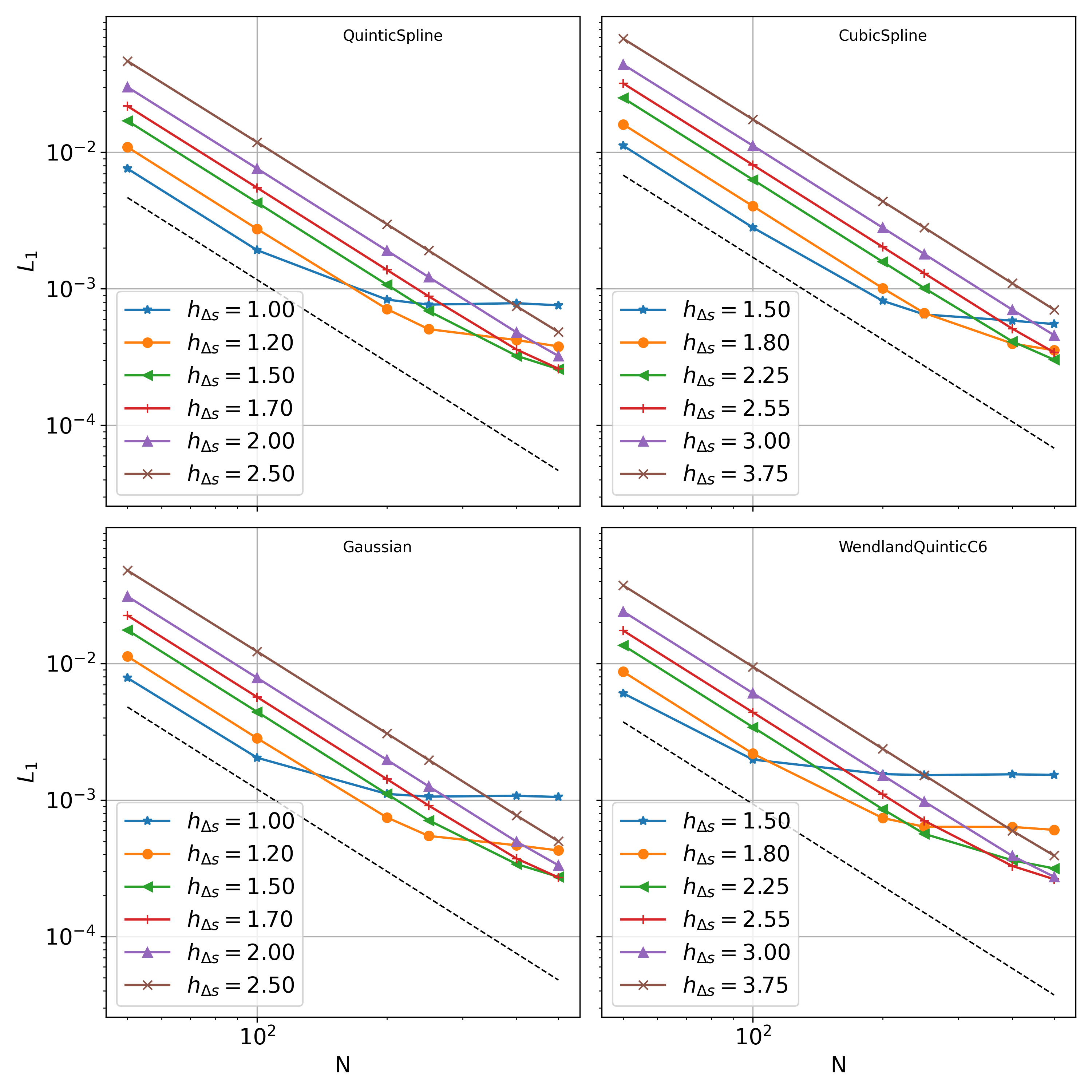}
  \caption{Order of convergence for function approximation for different
    $h_{\Delta s}$ values in a PP domain. The dashed line shows the second
    order rate.}
  \label{fig:d0p1_f}
\end{figure}

\begin{figure}[ht!]
  \centering
  \includegraphics[width=\linewidth]{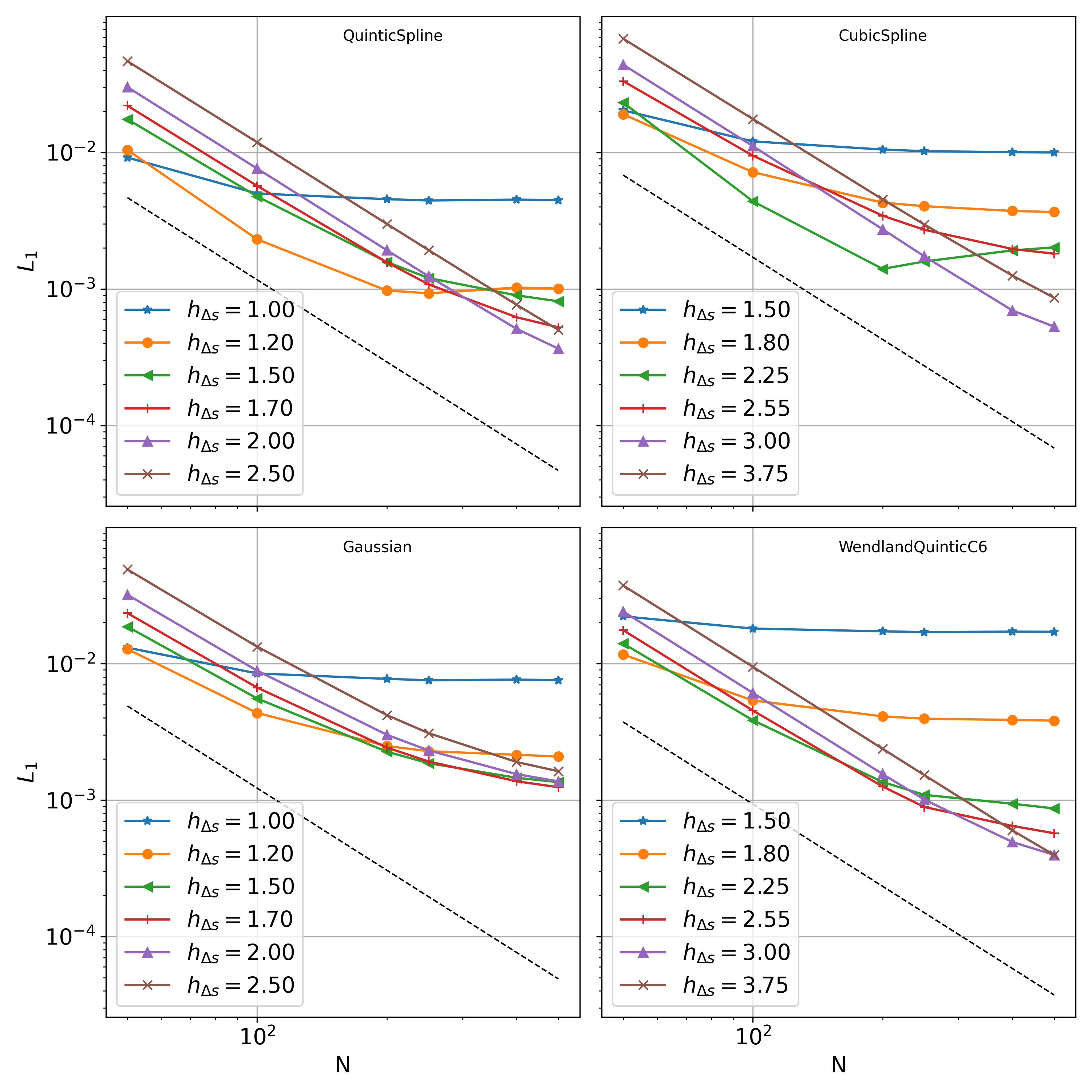}
  \caption{Order of convergence for the derivative approximation for different
    $h_{\Delta s}$ values in a PP domain. The dashed line shows the second
    order rate. The \cref{eq:sph_dfij} is used for the approximation.}
  \label{fig:d1p1_df}
\end{figure}

We compare the four kernels $G$, $CS$, $QS$ and $WQ_6$ for convergence of function and its gradient approximation. We consider the field,
\begin{equation}
  f = \sin(\pi (x + y)).
  \label{eq:function}
\end{equation}
Given a function $g_o$ and its approximation $g$, we evaluate the $L_1$ error
using,
\begin{equation}
  L_1 = \frac{\sum_i^N |g(\ten{x_i}) - g_o(\ten{x_i})|}{\sum_i^N |g_o(\ten{x}_i)|},
  \label{eq:l1_error}
\end{equation}
where $N$ is the total number of particles in the domain. Since the $CS$ and
$WQ_6$ kernels have support radius of $2$ whereas, the $G$ and $QS$ kernel have
support radius of $3$, we set the $h_{\Delta s}$ such that the $N_{nbr}$ is
same in an UP domain. Therefore, when $h_{\Delta s} = 1.0$ for $QS$ (or $G$),
we take $h_{\Delta s}=1.5$ for the $CS$ (or $WQ_6$). For the convergence
study, in this paper, we consider $50\times50$, $100\times100$,
$200\times200$, $250\times250$, $400\times400$, and $500\times500$ resolutions
for all the test cases unless stated otherwise.

\subsection{Unperturbed periodic domain}

In \cref{fig:d0p0_f}, we plot the $L_1$ in the function approximation as a
function of the resolution for different values of $h_{\Delta s}$ in a UP
domain. We observe similar error values for all the kernels except $CS$. We
obtain second-order convergence (SOC) in an UP domain upto a considerably high
resolution of $500 \times 500$ as expected, for all the kernels
\cite{kiara_sph_2013}.

In \cref{fig:d1p0_df}, we plot the $L_1$ error in the derivative approximation
of the function in \cref{eq:function} in a UP domain. The $G$ and $QS$ kernels
show a better convergence rate compared to $CS$ and $WQ_6$ for lower $h_{\Delta
s}$. The $G$ kernel does not show SOC even at $h_{\Delta s}=2.5$, since we use a
truncated Gaussian. The $CS$ and $WQ_6$ kernel shows SOC only when
$h_{\Delta s} \geq 3.0$. The $QS$ kernel shows SOC at $h_{\Delta s} > 1.5$
however, a reasonable convergence can be seen for $h_{\Delta s} = 1.2$ as well.

\subsection{Perturbed Periodic domain}

In an SPH simulation, the particles advect with different velocities, and thus
the distribution of particles is no longer uniform. Some particle shifting
techniques (PST) can be used to make the particle distribution uniform
\cite{acc_stab_xu:jcp:2009,lind2012incompressible}. Thus, it is essential to
observe the convergence rate in the PP domain as well.

In \cref{fig:d0p1_f}, we plot the $L_1$ error of the approximation of the field
given in \cref{eq:function} as a function of resolution in a PP domain for
different $h_{\Delta s}$ values. The convergence rates tend to zero for
higher resolution for low value of $h_{\Delta s}$ for all the kernels. The
$WQ_6$ kernel performs worse than the $CS$ kernel at lower $h_{\Delta s}$ values
however, the errors are significantly lower in $WQ_6$ when the $h_{\Delta s}$
value increase. On comparing $G$ and $QS$, the error plot looks exactly same
except when $h_{\Delta s} = 1.0$.

The SPH approximation of the gradient of a function is not even zero order
accurate in a perturbed domain \cite{nguyen_meshless_2008,kiara_sph_2013}. The
derivatives diverge when we evaluate it using \cref{eq:sph_df}. We use a zero
order consistent method proposed by \citet{sph:fsf:monaghan-jcp94} to compare
the kernels. We write this approximation as,
\begin{equation}
  \left< \nabla f(\ten{x}_i) \right> = \sum_j (f(\ten{x}_j) - f(\ten{x}_i)) \nabla W_{ij} \omega_j.
  \label{eq:sph_dfij}
\end{equation}

In \cref{fig:d1p1_df}, we plot the $L_1$ error in the function derivative
approximation as a function of resolution using \cref{eq:sph_dfij} in a PP
domain for different $h_{\Delta s}$ values and kernels. Clearly, the
approximation for all the kernels shows at least zero-order convergence. The $G$
kernel does not show SOC for high $h_{\Delta s}$ which is the same as observed
in the case of the UP domain. The accuracy in the case of $QS$ and $CS$
oscillates when going from lower $h_{\Delta s}$ to higher values.
\citet{zhu2015numerical} suggest that one should increase the $h_{\Delta s}$ as
one increases the resolution but given the inconsistent behavior of the $CS$ and
$QS$ kernels; these may not be suitable for that approach. The zero-order
convergence rate occurs due to dominance of discretization error (the term
$\left(\frac{\Delta s}{h} \right)^{\beta + 4}$ in \cref{eq:sph_df}) when the
resolution increases in the PP domain.

\section{Comparison of discretization operators}
\label{apn:comp_disc}

\subsection{Comparison of $\frac{\nabla p}{\Varrho}$ approximation}
\label{apn:grad_comp}

\begin{figure*}
  \centering
  \includegraphics[width=0.9\linewidth]{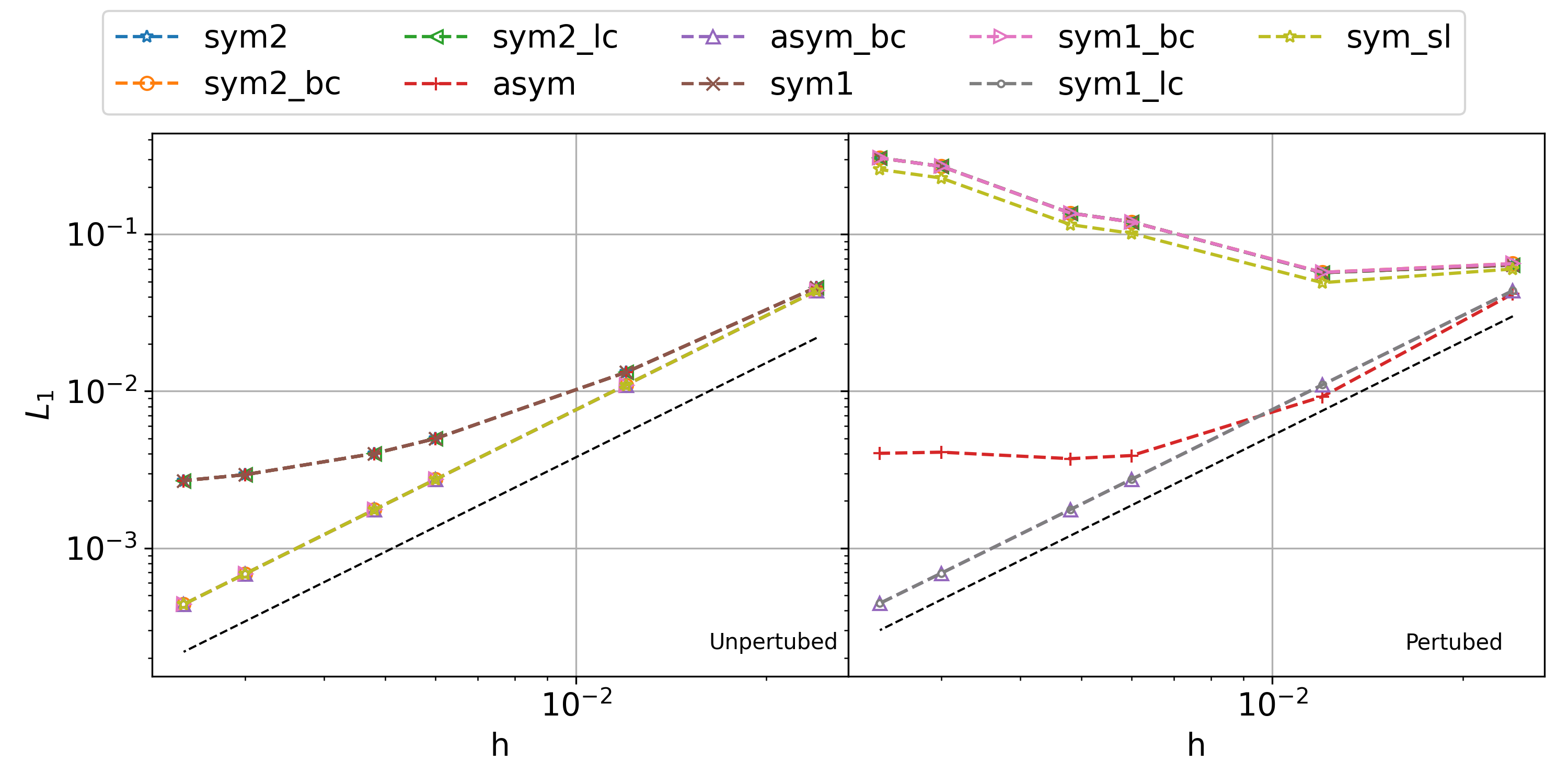}
  \caption{The rate of convergence in UP (left) and PP(right) domains for
    various pressure gradient listed in \cref{tab:grad}. The dashed line shows
    the SOC rate. The \nam{\_bc} and \nam{\_lc} suffixes represent the corresponding form
    with Bonet correction \cite{bonet_lok:cmame:1999} and Liu correction
    \cite{liu_restoring_2006}. The \nam{sym\_sl} is the \nam{sym1} formulation with
    symmetrization of kernel proposed in \citet{dilts1999moving}.}
  \label{fig:grad_p}
\end{figure*}

\begin{table}
  \centering
  \begin{tabular}{clll}
\toprule
          Name & $\frac{F_T}{F_{max}}$ & $T_{r}$ &        $L_1(O)$ \\
\midrule
\nam{asym\_bc} &                  0.01 &    1.97 &  4.45e-04(1.98) \\
    \nam{asym} &                  0.00 &    1.00 &  4.02e-03(0.07) \\
\nam{sym1\_bc} &                  0.01 &    1.87 & 3.06e-01(-0.55) \\
\nam{sym1\_lc} &                  0.01 &    2.35 &  4.45e-04(1.98) \\
    \nam{sym1} &                 -0.00 &    1.05 & 3.06e-01(-0.55) \\
\nam{sym2\_bc} &                  0.01 &    1.87 & 3.06e-01(-0.55) \\
\nam{sym2\_lc} &                 -0.00 &    2.43 & 3.06e-01(-0.55) \\
    \nam{sym2} &                 -0.00 &    1.01 & 3.06e-01(-0.55) \\
 \nam{sym\_sl} &                 -0.00 &    2.50 & 2.59e-01(-0.56) \\
\bottomrule
\end{tabular}

  \caption{The ratio $\frac{F_T}{F_{max}}$ showing the total force in the system
  due to lack of conservation in the approximation, the time taken $T_r$
  relative to the \nam{asym} formulation, the $L_1$ error for $500\times500$
  particle in a PP domain, and last column shows the order of convergence for all the kind of formulations considered in the first column.}
  \label{tab:grad_time}
\end{table}

In this section, we compare various pressure gradient approximations. In the
\cref{tab:grad}, we list the gradient approximations considered in this study.
The \nam{sym1} and \nam{sym2} are the symmetric, conservative form of the
gradient approximation. We note that conservative forms have $\Varrho=\rho$.
The \nam{asym} is the asymmetric form. Since the SPH kernel gradient does
not show SOC in a perturbed domain~\cite{quinlan_truncation_2006}, we also
consider the kernel correction employed to each of the approximation. In this
paper, we refer to the correction proposed by \citet{bonet_lok:cmame:1999} as
\emph{Bonet correction} and the one proposed by \citet{liu_restoring_2006} as
\emph{Liu correction}. We add the suffix \nam{\_bc}, and \nam{\_lc}
respectively in the plots and tables to indicate these corrections. The
application of corrections renders the symmetric forms non-conservative, we
use the method of symmetrization of the kernel proposed by
\citet{dilts1999moving} to again make it conservative. We refer to this
formulation as \nam{sym\_sl} which we write as
\begin{equation}
  \sum_j m_j \frac{p_j + p_i}{\rho_j \rho_i} (L_i \nabla W_{ij} - L_j \nabla W_{ji}),
  \label{eq:symm_corr_conv}
\end{equation}
where $L_i$ is the Liu correction applied to the kernel gradient
\footnote{We select the \nam{sym1} formulation over \nam{sym2} as the
latter does not perform well with the linear correction (see
\cref{sec:kern_grad} for details).}. This formulation is used in the scheme
proposed by \citet{crksph:jcp:2017}.

In order to compare the convergence, we consider a pressure field, $p=\sin(\pi
(x +y))$. We determine the $L_1$ error using \cref{eq:l1_error}, where
$g(\ten{x}_i)$ is the pressure gradient evaluated using the approximation and
$g_o({\ten{x}_i})$ is the exact pressure gradient. The exact pressure gradient,
$\nabla p = \pi \cos(\pi (x + y)) (\hat{\ten{i}} + \hat{\ten{j}})$. We compare
only the x-component of the results. In \cref{fig:grad_p}, we plot the error in
the various gradient approximations discussed above in both an UP and PP domain.
In the UP domain, barring the \nam{sym2\_lc}, all the corrected gradient
approximations behave the same, whereas the uncorrected gradients do not display
SOC. The corrected versions retains SOC even at high resolution since it reduces
the discretization error in the approximation
\cite{fatehi_error_2011,quinlan_truncation_2006}. We also observe that with the
correction the second term involving the $p_i$ term is zero in an UP domain
leading to the same expression.

In the case of the PP domain, we observe that both \nam{sym1} and
\nam{sym2} and their corresponding \nam{\_bc} versions overlap. The
symmetric formulations show an increase in the error in the approximation
with increasing resolution as suggested in \citet{fatehi_error_2011}.
Furthermore, as discussed in \cref{sec:disc}, the Bonet correction does not
correct the symmetric formulations. Clearly, the \nam{asym} formulation
shows better convergence, and the Bonet correction version shows SOC.
Therefore, the Bonet correction can be applied only when an asymmetric
formulation is employed.  Moreover, the Liu correction only corrects the
symmetric form \nam{sym1}, which suggests that the \nam{sym2} cannot be
corrected using traditional correction techniques. Finally, the
\nam{sym\_sl} method has a slightly lower error but looses SOC behavior due
to the symmetrization of the kernel gradient. \citet{crksph:jcp:2017}
reported the similar behavior.

We also compare the linear momentum conservation and time taken to evaluate the
gradient for the case with $500 \times 500$ particles. As shown in
\citet{bonet_lok:cmame:1999}, linear momentum is conserved when the total force,
$\sum_i F_i = 0$, where the sum is taken over all the particles and $F_i=
\frac{\nabla p_i}{\Varrho_i}$. In \cref{tab:grad_time}, we tabulate the ratio of
total force to the maximum force ($\max(F_i)$), the time taken to evaluate the
gradient scaled by the minimum time taken by all the methods, and the $L_1$
error with the order of convergence~\footnote{In this paper, we report order of
convergence by fitting a linear regression line and finding its slope.}, for all
the formulations plotted in \cref{fig:grad_p}. As expected, all the symmetric
forms of approximation have zero total force. The asymmetric formulation has a
very small total force. Clearly, the use of Bonet corrections increases the
total force and slows down the computation by a factor of $2$, whereas the Liu
correction makes it $2.4$ times slower. The \nam{sym\_sl} formulation shows zero
residual force as expected. Using the \cref{tab:grad_time}, we can see that
\nam{asym\_bc} and \nam{sym1\_lc} show SOC and have a very low total force which
makes them a suitable candidate for a scheme with SOC.

\subsection{Comparison of $\nabla \cdot \ten{u}$ approximation}
\label{apn:comp_div}

\begin{figure*}[ht!]
  \centering
  \includegraphics[width=0.9\linewidth]{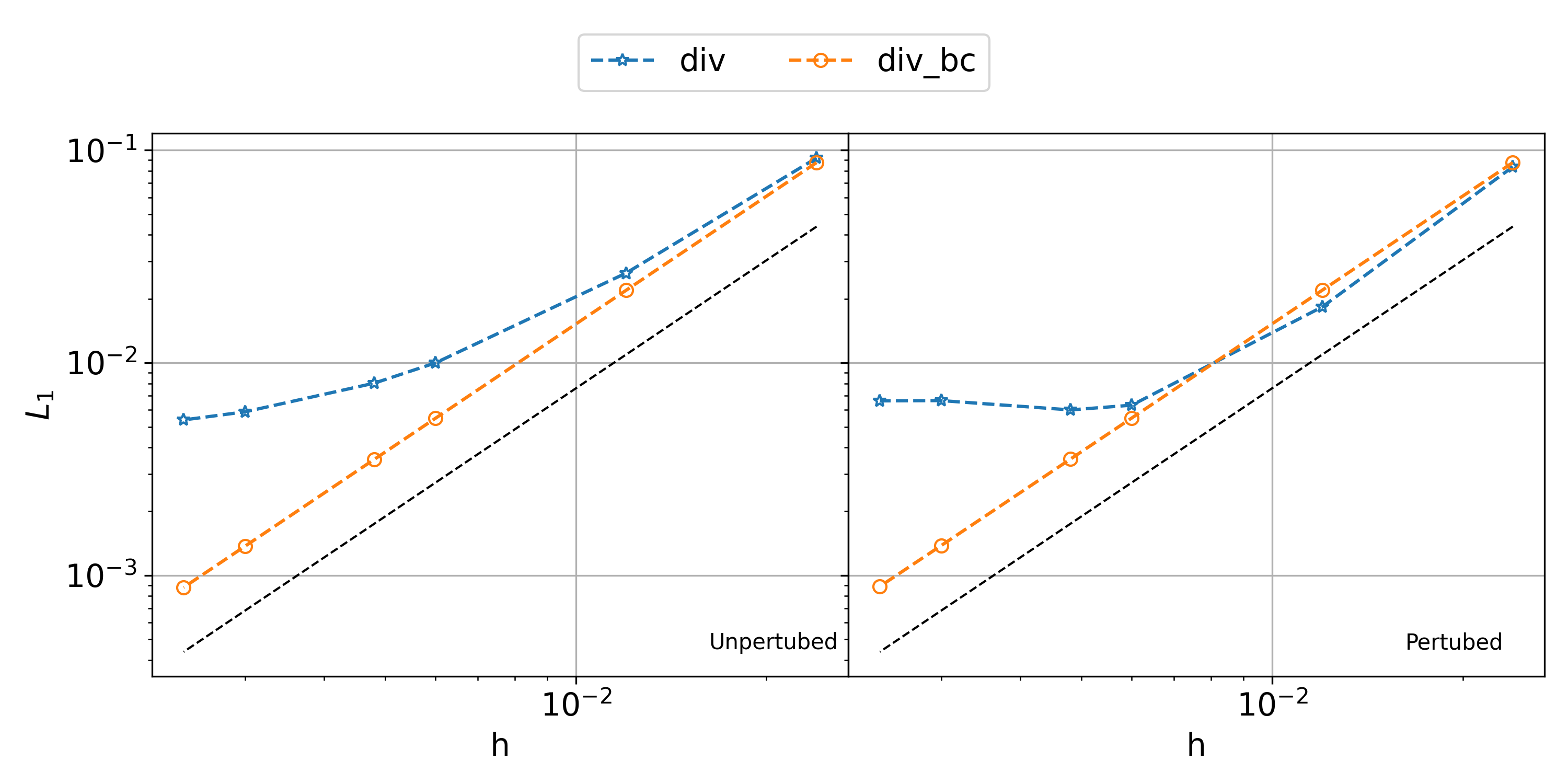}
  \caption{The rate of convergence UP (left) and PP(right) domains for
    velocity divergence in \cref{eq:div}. The dashed line shows the SOC rate.
    The suffix \nam{\_bc} represents the corresponding form with Bonet correction.}
  \label{fig:div_u}
\end{figure*}

A zero order consistent SPH approximation for the divergence
operator~\cite{monaghan-review:2005} is,
\begin{equation}
  \left< \nabla \cdot \ten{u}\right>_i = \sum_j (\ten{u}_j - \ten{u}_i) \cdot \nabla W_{ij} \omega_j.
  \label{eq:div}
\end{equation}
We refer to the approximation given in \cref{eq:div} as \nam{div}. We apply the
Bonet correction as done in the case of gradient approximation for a
first-order consistent approximation. We refer to the corrected form as
\nam{div\_bc}.

We consider the velocity field, $\ten{u} = \sin(\pi (x +y)) (\hat{\ten{i}} +
\hat{\ten{j}})$. The divergence of the velocity is given by, $\nabla \cdot
\ten{u} = 2 \pi \cos(\pi (x +y))$. We evaluate the $L_1$ error in the
approximation using \cref{eq:l1_error}. In \cref{fig:div_u}, we plot the $L_1$
error in the divergence approximation in an UP and PP domain. The uncorrected
approximation does not display SOC since the discretization error dominates as
we approach higher resolutions. Clearly, the corrected form shows SOC even in
the case of a PP domain.

In order to evaluate the accuracy of the approximation in a divergence-free
field, we consider the velocity field,
\begin{equation}
  \begin{split}
    u &= - \cos(2 \pi x) \sin(2  \pi y),\\
    v &= \sin(2 \pi x) \cos(2  \pi y).\\
  \end{split}
  \label{eq:div_free}
\end{equation}
\begin{figure*}[ht!]
  \centering
  \includegraphics[width=0.9\linewidth]{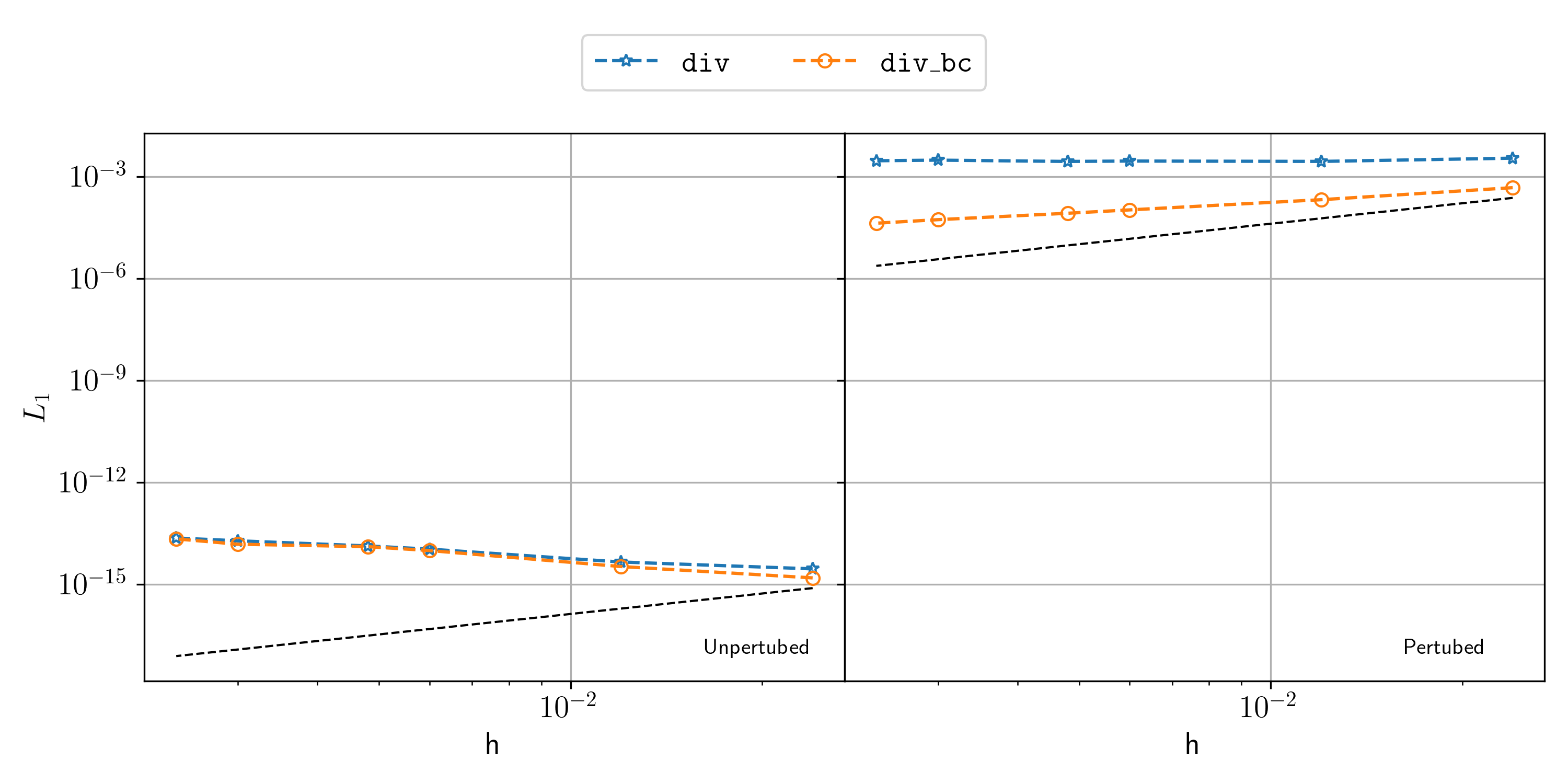}
  \caption{The rate of convergence UP (left) and PP(right) domains for velocity
  divergence for a divergence-free field. The dashed line shows the SOC rate.
  The suffix \nam{\_bc} represents the corresponding form with Bonet
  correction.}
  \label{fig:div_free_u}
\end{figure*}
In the \cref{fig:div_free_u}, we plot the $L_1$ error using
\cref{eq:l1_error} in the divergence computation as a function of
resolution for the UP and PP domain. Clearly, the divergence is zero in a
UP domain owing to the symmetry of the particles. However, the error in PP
domain remains about the same order as seen in the case of general field in
\cref{fig:div_u}~\footnote{See error in divergence approximation in
\cref{sec:div_error}}. Clearly, the Bonet correction does not correct this
issue. We observe the implication of this behavior when we compare the
schemes in \cref{sec:comp_new}.

The continuity equation corresponds to the mass conservation of
the system, since mass of each particle is kept constant, we satisfy the
global conservation of mass implicitly.

\subsection{Comparison of $\nabla^2 \ten{u}$ approximation}
\label{apn:comp_lap}

\begin{table}[ht!]
  \centering
  \renewcommand{\arraystretch}{2}
  \begin{tabular}{l|l|l}
    \text{Name} & \text{Expression} & \text{Used in} \\
    \hline
    \text{\nam{Cleary}} & \text{$  2 (\ten{u}_i - \ten{u}_j) \frac{\nabla W_{ij} \cdot \ten{x}_{ij}}{|\ten{x}_{ij}|^2} \omega_j$} & \text{WCSPH \cite{cleary1999conduction,wcsph-state-of-the-art-2010}} \\
    \text{\nam{Fatehi}} & \text{$  2 \omega_j \left( \frac{(\ten{u}_i - \ten{u}_j)}{|\ten{x}_{ij}|} - \frac{\ten{x}_{ij} \cdot (\nabla \ten{u})_i}{|\ten{x}_{ij}|} \right) \frac{\nabla W_{ij} \cdot \ten{x}_{ij}}{|\ten{x}_{ij}|}$} & \text{modified WCSPH \cite{fatehi_error_2011}} \\
    \text{\nam{tvf}} & \text{$  \frac{1}{m_i}\left( \omega_i^2 +  \omega_j^2 \right)(\ten{u}_i - \ten{u}_j) \frac{\nabla W_{ij} \cdot \ten{x}_{ij}}{|\ten{x}_{ij}|^2}$} & \text{TVF \cite{Adami2013}, EDAC \cite{edac-sph:cf:2019}} \\
    \text{\nam{coupled}} & \text{$((\nabla \ten{u})_j - (\nabla \ten{u})_i) \cdot \nabla W_{ij} \omega_j$} & \text{\citet{bonet_lok:cmame:1999}} \\
  \end{tabular}
  \renewcommand{\arraystretch}{1}
  \caption{The various approximations of $\nabla^2 \ten{u}$. The column ``expression'' is
    assumed to be summed over the index $j$ over all the neighbor particles
    inside the kernel support. The $\nabla u_i$ term are calculated using
    first-order consistent formulation i.e.~\nam{asym\_bc}}
  \label{tab:lap}
\end{table}

\begin{figure*}[ht!]
  \centering
  \includegraphics[width=0.9\linewidth]{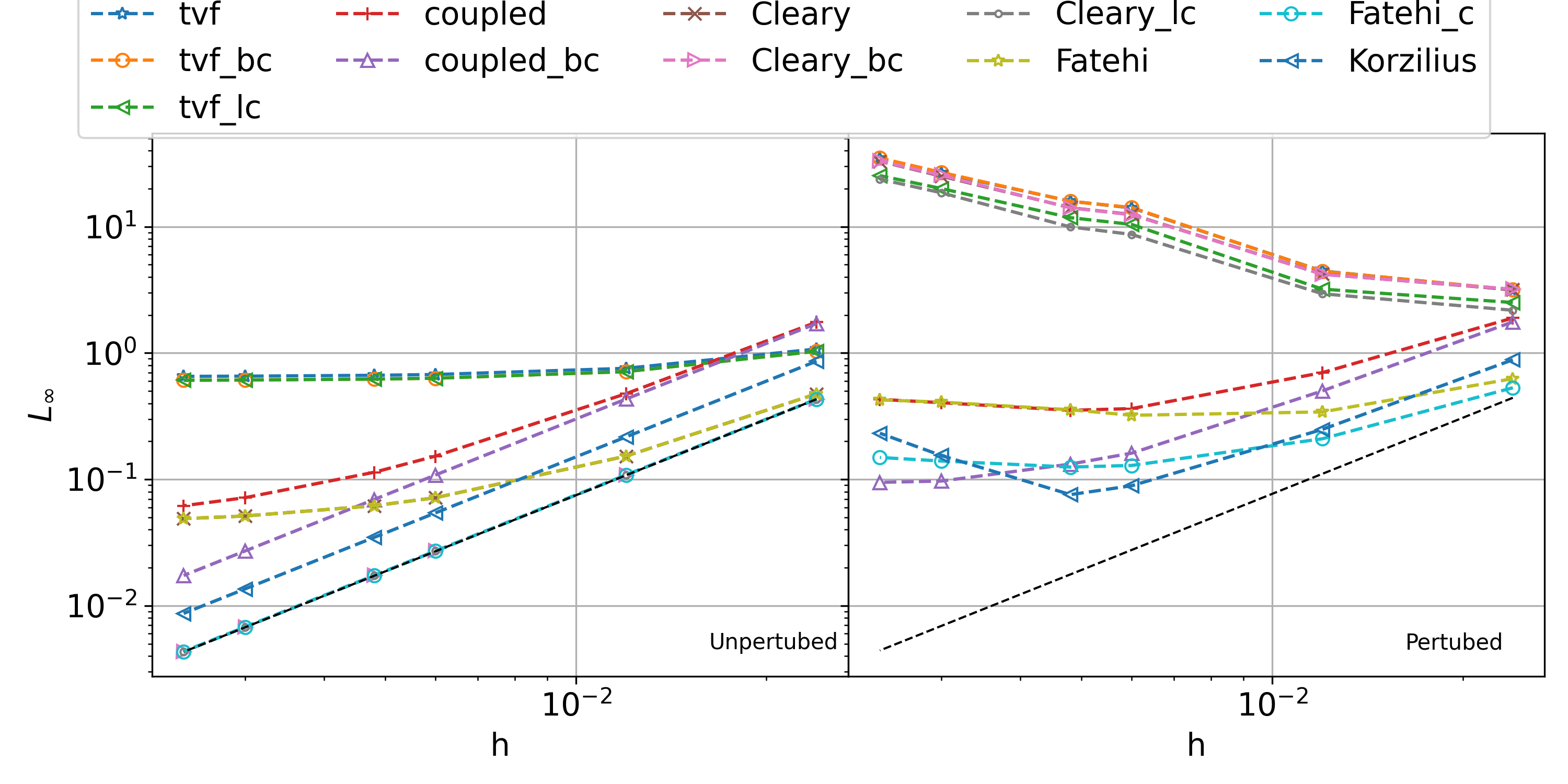}
  \caption{The rate of convergence UP (left) and PP(right) domains for various
  approximations of the Laplacian operator in \cref{tab:lap}. The dashed line
  shows the SOC rate. The suffixes \nam{\_bc} and \nam{\_lc} represent the
  corresponding form with the Bonet correction and Liu correction, respectively.
  The \nam{Fatehi\_c} refers to the \nam{fatehi} formulation with the correction
  proposed by \citet{fatehi_error_2011} (See \cref{apn:fatehi}).}
  \label{fig:lap_u}
\end{figure*}

\begin{table}[ht!]
  \centering
  \begin{tabular}{cllll}
\toprule
             Name & $\frac{F_T}{F_{max}}$ & $T_{r}$ & $L_1$ error & Order \\
\midrule
 \nam{Cleary\_bc} &             -1.28e+00 &    1.95 &    6.55e-02 & -0.83 \\
 \nam{Cleary\_lc} &             -2.15e-01 &    2.43 &    4.59e-02 & -0.79 \\
     \nam{Cleary} &             -1.08e-10 &    1.12 &    6.54e-02 & -0.84 \\
  \nam{Fatehi\_c} &              1.69e+00 &    4.10 &    2.88e-04 &  1.30 \\
     \nam{Fatehi} &              1.30e+00 &    2.70 &    8.43e-04 &  0.68 \\
  \nam{Korzilius} &              1.61e+00 &    4.66 &    2.49e-04 &  1.70 \\
\nam{coupled\_bc} &              1.61e+00 &    3.05 &    2.54e-04 &  1.95 \\
    \nam{coupled} &              1.30e+00 &    2.59 &    8.38e-04 &  1.46 \\
    \nam{tvf\_bc} &             -1.26e+00 &    2.07 &    6.80e-02 & -0.66 \\
    \nam{tvf\_lc} &             -2.16e-01 &    2.35 &    5.01e-02 & -0.56 \\
        \nam{tvf} &             -1.06e-10 &    1.00 &    6.77e-02 & -0.67 \\
\bottomrule
\end{tabular}

  \caption{The ratio $\frac{F_T}{F_{max}}$ showing the total force in the system
  due to lack of conservation of the approximation and the time taken, $T_r$
  relative to the \nam{tvf} formulation, and $L_1$ error for $500\times500$
  particle case in a PP domain. The last column shows order of convergence for
  all the methods listed in first column.}
  \label{tab:lap_time}
\end{table}

In this section, we compare various approximations for the Laplacian operator
listed in \cref{tab:lap}. We refer to the symmetric formulations of
\onlinecite{cleary1999conduction,brookshaw_method_1985} as \nam{Cleary}, and those
of \onlinecite{Adami2013} as \nam{tvf}. These ensure that linear momentum is
conserved. We also consider the coupled formulation used by
\onlinecite{bonet_lok:cmame:1999,nugent_liquid_2000} and refer to these as
\nam{coupled}. This formulation shows oscillations in the approximation when
the initial condition is discontinuous. However, to remedy this, one can
perform a first-order accurate approximation near the discontinuity and then
perform this approximation as shown in
\onlinecite{biriukovStableAnisotropicHeat2019a}. We consider the improved
formulation proposed by \citet{fatehi_error_2011} referred as \nam{Fatehi}.
Both the \nam{coupled} and \nam{Fatehi} formulations are asymmetric. These
formulations are performed in two steps where the first step involves
computation of velocity gradient for each particle. We also consider the
kernel correction applied to each of these formulations. In the case of
\nam{Cleary}, \nam{tvf}, and \nam{coupled} methods, we use the standard Bonet
and Liu corrections. However, in the case of \nam{Fatehi}, we use the
correction tensor proposed by \citet{fatehi_error_2011} given by
\begin{equation}
  \begin{split}
    \hat{B}_i^{\eta \mu} =& -\bigg(
  \sum_j  \omega_j
   \partial^{\mu} W_{ij} x_{ij}^{\eta}  +\\
   ~& \sum_j \omega_j  r_{ij}^2  \partial^{\theta} W_{ij} B_i^{T,\theta \alpha} \sum_j
   \omega_j e_{ij}^{\alpha} e_{ij}^{\eta}  \partial^{\mu}
   W_{ij} \bigg)^{-1 },
  \end{split}
  \label{eq:fatehi_corr}
\end{equation}
where the subscripts are SPH summation indices and superscripts are in tensor
notation \footnote{The term $*_{ij}$ is widely used in SPH literature and so
  we use tensor indices in the superscript to derive \cref{eq:fatehi_corr} in
  \cref{apn:fatehi}.}. We refer to this corrected formulation as
\nam{Fatehi\_c}. Additionally, the second derivative can also be obtained by
taking the double derivative of the kernel
\cite{chenGeneralizedSmoothedParticle2000a,zhangModifiedSmoothedParticle2004,korziliusImprovedCSPMApproach2017}.
Therefore, we also consider the method proposed by
\citet{korziliusImprovedCSPMApproach2017} that remedies the deficiencies in
earlier approaches where the second derivative was employed. We obtain the
second derivative of a scalar field using,
\begin{equation}
  \tilde{\nabla^2} u_i = (\tilde{\Gamma})^{-1} \left( \sum (u_j - u_i) \tilde{\nabla^2 W}_{ij} \omega_j  - \ten{x}_{ij} \left< \nabla u\right>_i \tilde{\nabla^2 W}_{ij} \omega_j   \right)
  \label{eq:korzilius}
\end{equation}
where $\tilde{\nabla^2} = \left[\frac{\partial^2}{\partial x^2}, \frac{\partial^2 }{\partial
x \partial y}, \frac{\partial^2 }{\partial y^2}\right]^{T}$ is the operator,
$\tilde{\nabla^2 W}_{ij} = \left[\frac{\partial^2 W_{ij}}{\partial x^2}, \frac{\partial^2 W_{ij}}{\partial
x \partial y}, \frac{\partial^2 W_{ij}}{\partial y^2}\right]^{T}$. The gradient
$\left< \nabla u\right>_i$ is approximated using the \nam{asym\_bc} formulation. The
correction $\tilde{\Gamma}$ is given by
\begin{equation}
  \tilde{\Gamma}_i = \sum \frac{1}{2} \tilde{\nabla^2 W}_{ij} \zeta_{ij}^{T} \omega_j - \sum \tilde{\nabla^2 W}_{ij} x_{ji}^{T} \omega_j B^{-1} \sum \frac{1}{2} \nabla W_{ij} \zeta_{ij}^{T} \omega_j
  \label{eq:kor_corr}
\end{equation}
where $\zeta_{ij} = [x_{ij}^2, x_{ij}y_{ij}, y_{ij}^2]$ and $B$ is the
Bonet correction matrix (see \eqref{eq:bonet_ts}). We refer to this
formulation as \nam{Korzilius}.

In the \cref{fig:lap_u}, we plot the rate of convergence for the various
formulations discussed above in both UP and PP domains. In the UP domain, all
the methods at least show zeroth-order convergence. All methods without
corrections suffer from high discretization error that dominates at higher
resolutions~\cite{fatehi_error_2011}. When either Bonet or Liu corrections are
employed, \nam{Cleary}, \nam{coupled}, \nam{Fatehi\_c}, and \nam{Korzilius}
methods show SOC. The \nam{coupled} method is approximately half an order less
accurate as compared to \nam{Cleary} and \nam{Fatehi}. The accuracy of
\nam{Korzilius} method is in between the \nam{coupled} and \nam{Fatehi}
method. The \nam{tvf} method is very inaccurate as the discretization error
increases due to the introduction of $\omega_i^2 + \omega_j^2$.

It is important to note that in the PP domain, the symmetric methods diverge
due to discretization error of $O(\frac{\tilde{d}}{h^2} \frac{\Delta s}{h})$,
where $\tilde{d}$ is the deviation from the regular particle arrangement
\cite{fatehi_error_2011}. Only the \nam{coupled}, \nam{Korzilius}, and
\nam{Fatehi} methods show a positive convergence rate. On applying the
corresponding correction, the \nam{coupled}, \nam{Korzilius}, and \nam{Fatehi}
methods improve. The accuracy for \nam{coupled}, \nam{Korzilius}, \nam{Fatehi}
is maintained as observed in the case of UP domain.

In the \cref{tab:lap_time}, we tabulate the total force as a result of the
approximation, the time taken for the approximation, and the error on a PP
domain consisting of $500\times 500 $ particles with the order of
convergence in the last column for each method plotted in \cref{fig:lap_u}.
We observe a similar increase in computational time due to the Bonet and
Liu corrections as seen in the case of gradient approximation.  The
\nam{coupled}, \nam{Korzilius}, and \nam{Fatehi} formulation have even
higher computational cost due to the additional step of velocity gradient
computation. The \nam{Korzilius} method requires additional time since the
double derivative of the kernel is involved. The \nam{Fatehi\_c} method has
an additional step where we compute the second-order tensor in
\cref{eq:fatehi_corr} for each particle resulting in a further increase in
computation time. We observe a similar increase in total force when an
asymmetric version of the formulation is employed, as seen in the case of
gradient approximation. Clearly, both the \nam{coupled}, \nam{Korzilius},
and \nam{Fatehi} formulation results in an equal amount of total force
resulting in a lack of conservation of linear momentum. In order to get a
SOC approximation, we can use either of \nam{coupled\_bc}, \nam{Korzilius},
or \nam{Fatehi\_c} formulations for viscous force estimation.

\section{The \nam{Cleary} and \nam{Fatehi} corrections}
\label{apn:fatehi}

In this section, we introduce the tensor notations for SPH that makes the
comprehension better. We use derivation for the error estimation from
\citet{fatehi_error_2011}. We write the Taylor series expansion of the
velocity component, $u_j$ defined at a point, $\ten{x}_j$ about a point $x_i$
as, expansion given by
\begin{equation}
  u_j = u_i - (\ten{x}_{ij}\cdot \nabla) u_i + \frac{1}{2} (\ten{x}_{ij} \cdot \nabla)^2 u_i - \frac{1}{6} (\ten{x}_{ij} \cdot \nabla)^3 u_i + \ \text{H.O.T}
  \label{eq:taylor_dot}
\end{equation}
where, $\ten{x}_{ij} = \ten{x}_i - \ten{x}_j$. Without loss of generality, we
consider only one component of velocity. We use tensor notation to represent
vector $\ten{x}_{ij}$ as $x_{ij}^{\alpha}$, where $i$ and $j$ are the particle
indices. We follow this notation since SPH approximation is performed using
sum over all its neighbors $j$. Thus, we write the \cref{eq:taylor_dot} in
this tensor notation as
\begin{equation}
  \begin{split}
  u_j =& u_i - x_{ij}^{\alpha} \partial^{\alpha} u_i + \frac{1}{2}
  x_{ij}^{\beta} x_{ij}^{\gamma}  \partial^{\beta} \partial^{\gamma}u_i -\\
  ~&\frac{1}{6} x_{ij}^{\alpha} x_{ij}^{\beta}  x_{ij}^{\gamma}
  \partial^{\alpha} \partial^{\beta} \partial^{\gamma} u_i + \text{H.O.T}.
  \end{split}
  \label{eq:taylor_tensor}
\end{equation}
We note that the subscripts are SPH notations and the superscripts are tensor notation indices.

We write the Laplacian of velocity, $\ten{u}$ using proposed by
\onlinecite{cleary1999conduction} as
\begin{equation}
  \left< \partial^{\eta} \partial^{\eta} u \right>_i = \sum_j 2 \omega_j (u_i - u_j) \frac{\partial^{\eta} W_{ij} x_{ij}^{\eta}}{r_{ij}^2}
  \label{eq:cleary}
\end{equation}
where $\left<*\right>$ is used to denote the approximation. We write the error, $E_i$
in the approximation as
 \begin{equation}
  E_i = \partial^{\eta} \partial^{\eta} u_i - \left< \partial^{\eta} \partial^{\eta} u \right>_i
\end{equation}
Using \cref{eq:taylor_tensor} and \cref{eq:cleary}, we obtain the error,
\begin{equation}
  \begin{split}
  E_i &= \partial^{\theta} \partial^{\theta}  u_i -  \sum_j 2 \omega_j \big[x_{ij}^{\alpha} \partial^{\alpha} u_i - \frac{1}{2} x_{ij}^{\beta} x_{ij}^{\gamma}  \partial^{\beta} \partial^{\gamma}u_i +\\
  & \frac{1}{6} x_{ij}^{\delta} x_{ij}^{\epsilon}  x_{ij}^{\zeta}  \partial^{\delta} \partial^{\epsilon} \partial^{\zeta} u_i + \text{H.O.T}\big] \frac{\partial^{\eta} W_{ij} x_{ij}^{\eta}}{r_{ij}^2}. \\
  \end{split}
\end{equation}
In the above equation, we can write $\partial^{\theta} \partial^{\theta}  u_i = \delta^{\theta \iota}\partial^{\theta} \partial^{\iota}  u_i$ and multiplying each term inside, we get,

\begin{equation}
  \begin{split}
  E_i &= -\partial^{\alpha} u_i \sum_j
  2 \omega_j e_{ij}^{\alpha} e_{ij}^{\eta}  \partial^{\eta}
  W_{ij}  +\\
  & \bigg(\delta^{\beta \gamma} + \sum_j \omega_j
  x_{ij}^{\beta} x_{ij}^{\gamma}
  \frac{\partial^{\eta} W_{ij} x_{ij}^{\eta}}{r_{ij}^2} \bigg) \partial^{\beta} \partial^{\gamma}u_i+ \text{H.O.T}\\
  \end{split}
  \label{eq:err_cleary}
\end{equation}

We can see that the first term is leading error term in the above equation. For a smoothing kernel, $W$ the term,
\begin{equation}
  \sum_j  \omega_j  (\ten{x}_{ij}\otimes \ten{x}_{ij})  \nabla W_{ij}
  \label{eq:m2_kern_grad}
\end{equation}
is the second moment of the kernel gradient. In a UP domain, the second moment
is zero. However, the leading term of the error is second moment scaled by
$\frac{1}{|x_{ij}|^2}$ which is still zero since it is a constant in a UP
domain. Whereas, the leading term is non-zero and causes the approximation to
deviate.

In the modified formulation proposed by \citet{fatehi_error_2011}, the
leading term is included in the approximation. We write the modified form as
\begin{equation}
  \left< \partial^{\theta} \partial^{\theta} u_i \right>_i = \sum_j 2 \omega_j ((u_i - u_j) - {x}_{ij}^{\alpha} \left< \partial^{\alpha} u \right>_i) \frac{\partial^{\eta} W_{ij} x_{ij}^{\eta}}{r_{ij}^2}
  \label{eq:fatehi}
\end{equation}
Using the similar algebraic manipulation, we write the error term as
\begin{equation}
  \begin{split}
  E_i &= \bigg(\sum_j \omega_j  x_{ij}^{\beta} x_{ij}^{\gamma}  \partial^{\theta} W_{ij} B_i^{T,\theta \alpha} \sum_j
  \omega_j e_{ij}^{\alpha} e_{ij}^{\eta}  \partial^{\eta}
  W_{ij}  +\\
  & \delta^{\beta \gamma} + \sum_j \omega_j
  x_{ij}^{\beta} x_{ij}^{\gamma}
  \frac{\partial^{\eta} W_{ij} x_{ij}^{\eta}}{r_{ij}^2} \bigg) \partial^{\beta} \partial^{\gamma}u_i+ \text{H.O.T}\\
  \end{split}
\end{equation}
where $B_i^T = \left(\sum_j \nabla W_{ij} \otimes (\ten{x}_j - \ten{x}_i)
\right)^{-T}$ is the correction matrix. \citet{fatehi_error_2011} also proposed a correction for the kernel gradient.
Let us assume the correction $\hat{B}_i^{\eta \mu}$ is applied to the kernel gradient. We write the modified equation as
\begin{equation}
  \left< \partial^{\theta} \partial^{\theta} u_i \right>_i = \sum_j 2 \omega_j ((u_i - u_j) - {x}_{ij}^{\alpha} \left< \partial^{\alpha} u \right>_i) \frac{\hat{B}_i^{\eta \mu} \partial^{\mu} W_{ij} x_{ij}^{\eta}}{r_{ij}^2}
  \label{eq:fatehi_mod}
\end{equation}
The Error in the above equation is given by
\begin{equation}
  \begin{split}
  E_i &= \bigg(\sum_j \omega_j  x_{ij}^{\beta} x_{ij}^{\gamma}  \partial^{\theta} W_{ij} B_i^{T,\theta \alpha} \sum_j
  \omega_j e_{ij}^{\alpha} e_{ij}^{\eta}  \hat{B}_i^{\eta \mu} \partial^{\mu}
  W_{ij}  +\\
  & \delta^{\beta \gamma} + \sum_j \omega_j
  x_{ij}^{\beta} x_{ij}^{\gamma}
  \frac{\hat{B}_i^{\eta \mu} \partial^{\mu} W_{ij} x_{ij}^{\eta}}{r_{ij}^2} \bigg) \partial^{\beta} \partial^{\gamma}u_i+ \text{H.O.T}\\
  \end{split}
\end{equation}
In order to make the approximation second order accurate, we must have the coefficient of $\partial^{\beta} \partial^{\gamma}u_i$ equal to zero. Thus we get,
\begin{equation}
  \begin{split}
  \sum_j \omega_j  x_{ij}^{\beta} x_{ij}^{\gamma}  \partial^{\theta} W_{ij} B_i^{T,\theta \alpha} \sum_j
  \omega_j e_{ij}^{\alpha} e_{ij}^{\eta}  \hat{B}_i^{\eta \mu} \partial^{\mu}
  W_{ij}  + &~\\
  \sum_j  \omega_j
  e_{ij}^{\beta} e_{ij}^{\gamma}
  \hat{B}_i^{\eta \mu} \partial^{\mu} W_{ij} x_{ij}^{\eta} &= -\delta^{\beta \gamma}\\
  \end{split}
\end{equation}
On inverting the system, we obtain,
\begin{equation}
  \begin{split}
    \hat{B}_i^{\eta \mu} =& -\bigg(
  \sum_j  \omega_j
   \partial^{\mu} W_{ij} x_{ij}^{\eta}  + \\
   ~&\sum_j \omega_j  r_{ij}^2  \partial^{\theta} W_{ij} B_i^{T,\theta \alpha} \sum_j
   \omega_j e_{ij}^{\alpha} e_{ij}^{\eta}  \partial^{\mu}
   W_{ij} \bigg)^{-1 }
  \end{split}
\end{equation}
The above equation is the correction matrix proposed by \onlinecite{fatehi_error_2011}
in a simple tensorial notation.

\section{The effect of solid-wall boundary conditions}
\label{apn:solid_bc}

There are many solid-wall boundary condition implementations in
SPH~\cite{Adami2012,marrone-deltasph:cmame:2011,randles1996smoothed,maciaTheoreticalAnalysisNoSlip2011}.
In this paper, we use the method due to
\citet{maciaTheoreticalAnalysisNoSlip2011} and \citet{Adami2012} that is
widely used in SPH. In order to apply the boundary condition, a few layers of
ghost particles are created outside the fluid domain such that the fluid
particles near the boundary have full support. The pressure and velocity on
the ghost particles are extrapolated from the fluid particles. The pressure
is determined using
\begin{equation}
  p_g = \frac{\sum_f p_f W_{gf}}{\sum_f W_{gf}},
\end{equation}
where $p_f$ is the pressure of the fluid particles, $W_{gf}$ is the kernel
weight between the ghost and fluid particle, and the sum is taken over all
the fluid particles near the ghost particle. The velocity on the ghost
particle is extrapolated using
\begin{equation}
  u_g = 2 u_s - \frac{\sum_f u_f W_{gf}}{\sum_f W_{gf}},
\end{equation}
where $u_s$ is the actual velocity of the solid, and $u_f$ is the velocity
of the fluid particles. In the L-IPST-C scheme, we use a slip boundary for
the continuity equation~\cite{muta2020efficient}. Additionally, since the
coupled formulation is prone to oscillations due to
discontinuity~\cite{quinlan_truncation_2006} we smooth-out oscillations
from the three fluid particle layers adjacent to the wall by setting
velocity values using a first order consistent interpolation.

We consider the Poiseuille flow problem. The exact solution of the
Poiseuille flow is given by
\begin{equation}
  u(y) = 0.5 \frac{F}{\nu} y (L - y)
\end{equation}
where $F=0.8N$ is the constant force applied on the flow, $\nu=0.1 m^2/s$
is the dynamic viscosity of the flow, $L=1m$ is the distance between the
parallel plates, and $y$ is the distance from the bottom plate. We consider
a domain of size $0.4 \times 1 m^2$ with maximum flow velocity $U=1 m/s$
and $Re = 10$. The domain is periodic in $x$-direction. We simulate the
problem for $10 s$ for each scheme for $50\times 50$, $100\times100$,
and $200\times200$ resolutions.

\begin{figure}[htbp]
  \centering
  \includegraphics[width=0.9\linewidth]{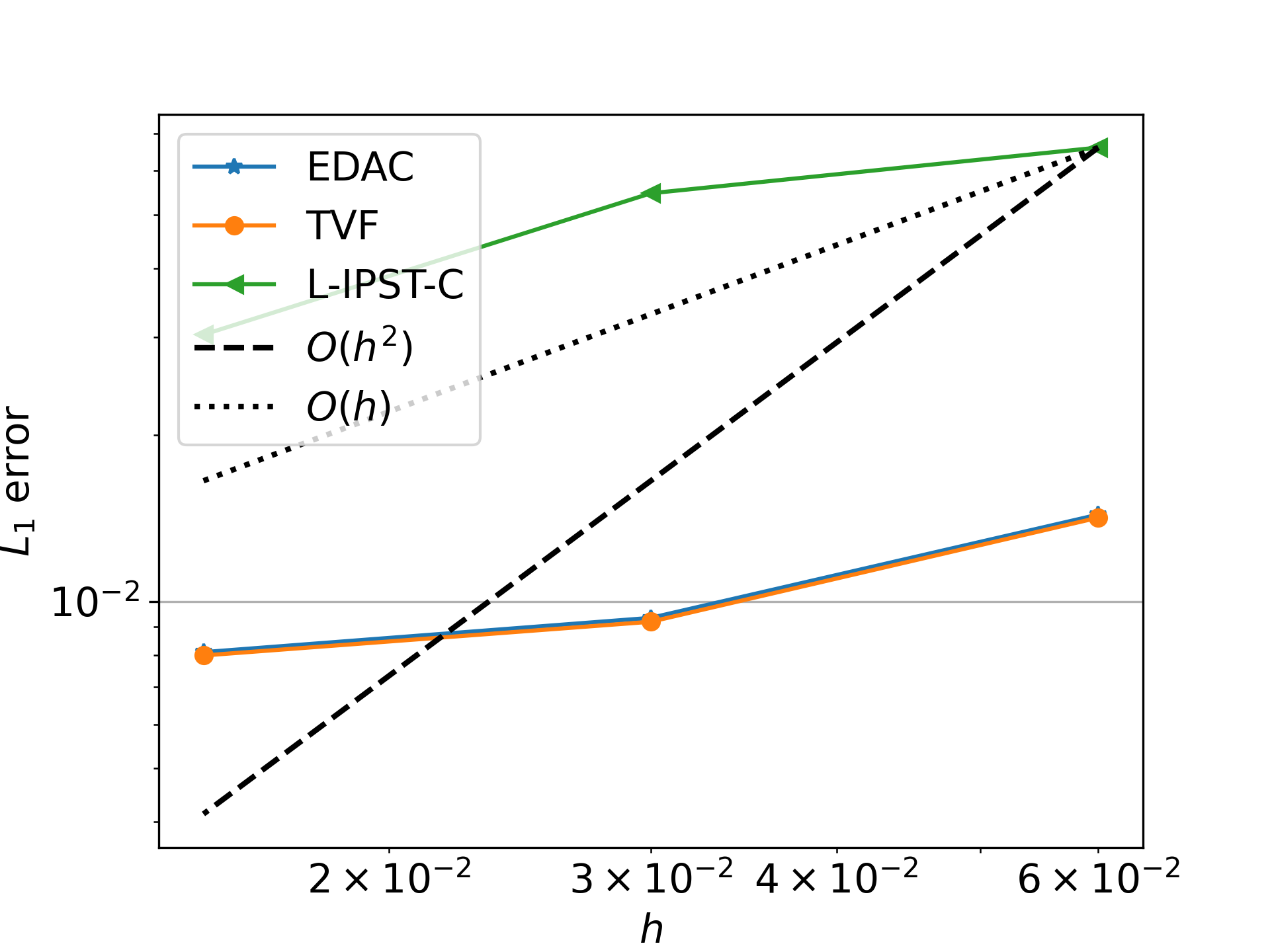}
  \caption{The convergence of error in the $x$-component of the velocity for Poiseuille flow problem.}
  \label{fig:pois}
\end{figure}

In \cref{fig:pois}, we plot the $L_1$ error in the $x$-component of the
velocity for all the schemes. Clearly, none of the schemes show convergence.
The L-IPST-C scheme shows a higher value of error due to the smoothing used
for the near wall fluid particles. It is clear from the above experiment
that we require a convergent boundary condition implementation such that it
does not dominate the error in the solution.

\section{$\delta^+$SPH formulation correction}
\label{apn:delta_plus_corr}

The the evolution equation of the $\delta^+$SPH equation has the form
\begin{equation}
  \frac{D f}{Dt} = \frac{df}{dt} + \nabla f \cdot \delta \ten{u},
\end{equation}
where $\frac{D f}{Dt} = \frac{\partial f}{\partial t} + (\ten{u} + \delta
\ten{u}) \cdot \nabla f$. The above equation can be written in terms of a
particle $i$ as,
\begin{equation}
  \frac{D f_i}{Dt} = \frac{df_i}{dt} + \nabla f_i \cdot \delta \ten{u}_i.
\end{equation}
We can use the vector identity for the last term,
\begin{equation}
  \begin{split}
  \nabla f \cdot \delta \ten{u} = \nabla \cdot (f \delta \ten{u}) - f \nabla \cdot (\delta \ten{u}).
  \end{split}
\end{equation}
On performing SPH approximation, we obtain
\begin{equation}
  \begin{split}
    \nabla f_i \cdot \delta \ten{u}_i &= \sum_j (f_j \delta \ten{u}_j-f_i \delta \ten{u}_i) \cdot \nabla W_{ij} \omega_j -\\
    ~& \sum_j f_i (\delta \ten{u_j} - \delta \ten{u}_i ) \cdot \nabla W_{ij} \omega_j\\
    & = \sum_j (f_j - f_i) \delta \ten{u}_j \cdot \nabla W_{ij} \omega_j.
  \end{split}
\end{equation}
Clearly, we cannot recover the LHS should we use the above discretization.
However, on using $f_j$ in place of $f_i$ in the second term, we get
\begin{equation}
  \begin{split}
    \nabla f_i \cdot \delta \ten{u}_i &= \sum_j (f_j \delta \ten{u}_j -f_i \delta \ten{u}_i) \cdot \nabla W_{ij} \omega_j -\\
    ~& \sum_j f_j (\delta \ten{u_j} - \delta \ten{u}_i ) \cdot \nabla W_{ij} \omega_j.\\
    & = \sum_j (f_j - f_i) \delta \ten{u}_i \cdot \nabla W_{ij} \omega_j.
  \end{split}
\end{equation}
Thus, in the $\delta^+$SPH we should use the above discretization.

\section{Schemes with issues solving the Gresho-Chan vortex}
\label{apn:gc}

\begin{figure}[ht!]
  \centering
  \includegraphics[width=0.9\linewidth]{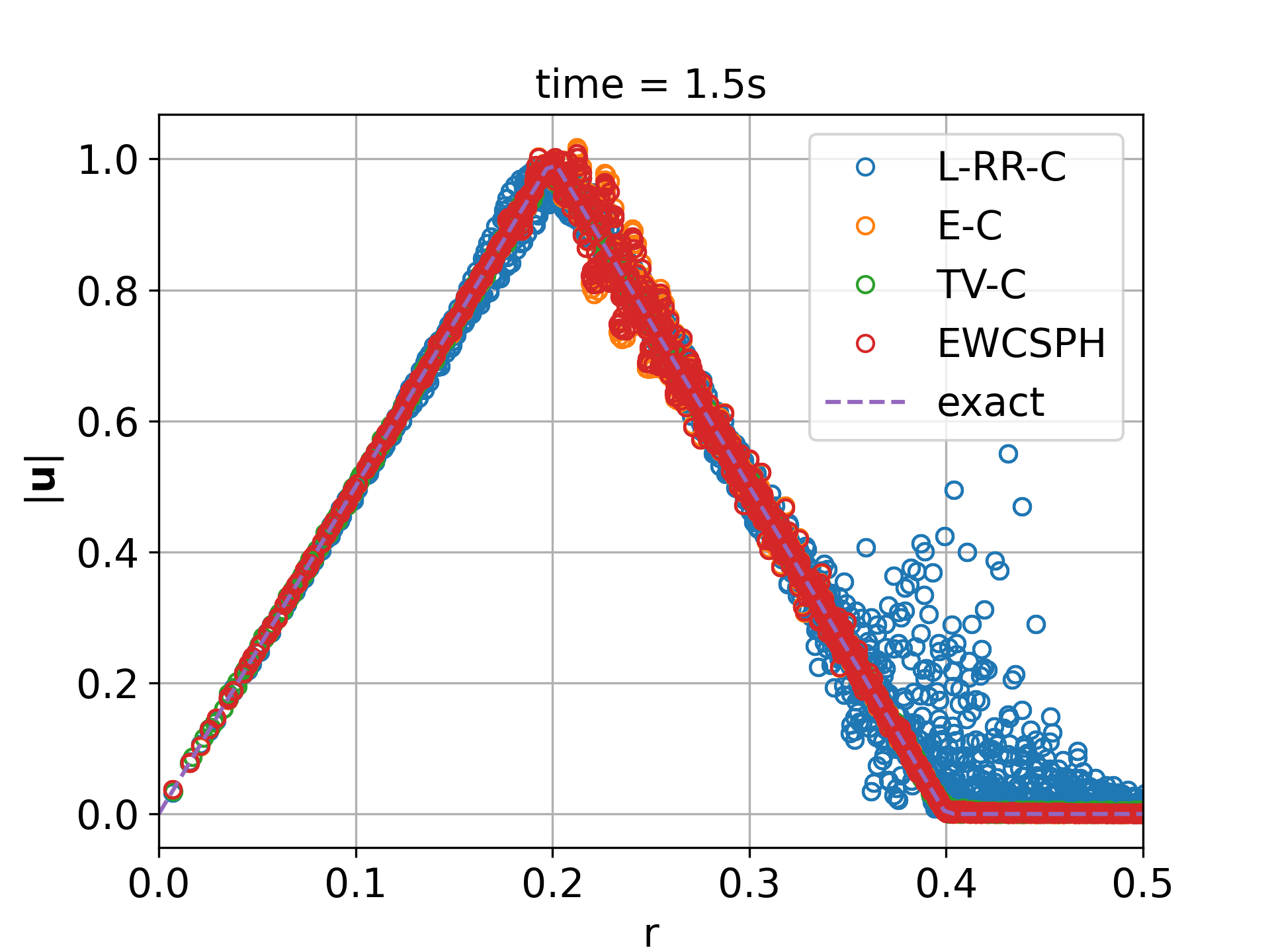}
  \includegraphics[width=0.9\linewidth]{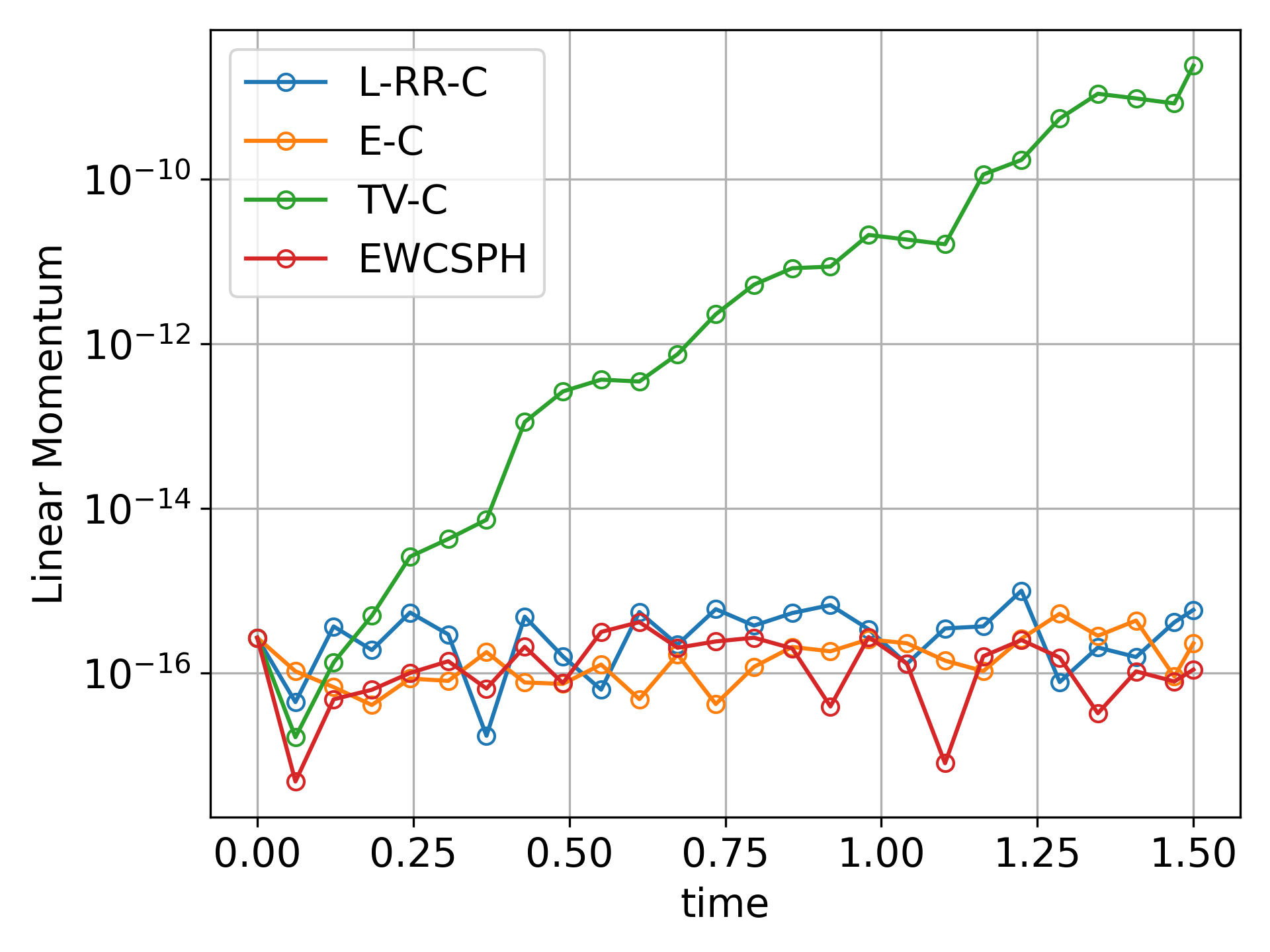}
  \caption{The velocity of particles with the distance from the center of the
    vortex (left) and the $x$-component of the total linear momentum (right)
    for the Gresho-Chan vortex problem.}
  \label{fig:gc_nw}
\end{figure}

In this section, we show the results for the scheme for which the Gresho-Chan
vortex problem failed to complete. In \cref{fig:gc_nw}, we plot the velocity
of the particles with the distance, $r$ from the center at $t=1.5s$, and the
linear momentum in the x-direction with time for a $100\times100$ simulation.
Clearly, all the schemes considered show better approximate conservation of
linear momentum compared to other scheme; however, they fail to complete.

In case of L-RR-C, due to the present of sharp change in the velocity field,
the remeshing procedure diverges \cite{remeshed_sph:jcp:2002}. In case of E-C,
TV-C and EWCSPH, we suspect that the advection term $\ten{u} \cdot \nabla
\ten{u}$ (or $\delta \ten{u} \cdot \nabla \ten{u}$ in case of TV-C) diverge in
the absence of viscosity. This opens possible avenues of research to obtain a
better discretization of the advection term.

\bibliography{references}

\begin{thebibliography}{67}%
\makeatletter
\providecommand \@ifxundefined [1]{%
 \@ifx{#1\undefined}
}%
\providecommand \@ifnum [1]{%
 \ifnum #1\expandafter \@firstoftwo
 \else \expandafter \@secondoftwo
 \fi
}%
\providecommand \@ifx [1]{%
 \ifx #1\expandafter \@firstoftwo
 \else \expandafter \@secondoftwo
 \fi
}%
\providecommand \natexlab [1]{#1}%
\providecommand \enquote  [1]{``#1''}%
\providecommand \bibnamefont  [1]{#1}%
\providecommand \bibfnamefont [1]{#1}%
\providecommand \citenamefont [1]{#1}%
\providecommand \href@noop [0]{\@secondoftwo}%
\providecommand \href [0]{\begingroup \@sanitize@url \@href}%
\providecommand \@href[1]{\@@startlink{#1}\@@href}%
\providecommand \@@href[1]{\endgroup#1\@@endlink}%
\providecommand \@sanitize@url [0]{\catcode `\\12\catcode `\$12\catcode
  `\&12\catcode `\#12\catcode `\^12\catcode `\_12\catcode `\%12\relax}%
\providecommand \@@startlink[1]{}%
\providecommand \@@endlink[0]{}%
\providecommand \url  [0]{\begingroup\@sanitize@url \@url }%
\providecommand \@url [1]{\endgroup\@href {#1}{\urlprefix }}%
\providecommand \urlprefix  [0]{URL }%
\providecommand \Eprint [0]{\href }%
\providecommand \doibase [0]{http://dx.doi.org/}%
\providecommand \selectlanguage [0]{\@gobble}%
\providecommand \bibinfo  [0]{\@secondoftwo}%
\providecommand \bibfield  [0]{\@secondoftwo}%
\providecommand \translation [1]{[#1]}%
\providecommand \BibitemOpen [0]{}%
\providecommand \bibitemStop [0]{}%
\providecommand \bibitemNoStop [0]{.\EOS\space}%
\providecommand \EOS [0]{\spacefactor3000\relax}%
\providecommand \BibitemShut  [1]{\csname bibitem#1\endcsname}%
\let\auto@bib@innerbib\@empty
\bibitem [{\citenamefont {Monaghan}(1994)}]{sph:fsf:monaghan-jcp94}%
  \BibitemOpen
  \bibfield  {author} {\bibinfo {author} {\bibfnamefont {J.~J.}\ \bibnamefont
  {Monaghan}},\ }\bibfield  {title} {\enquote {\bibinfo {title} {Simulating
  free surface flows with {SPH}},}\ }\href {\doibase 10.1006/jcph.1994.1034}
  {\bibfield  {journal} {\bibinfo  {journal} {Journal of Computational
  Physics}\ }\textbf {\bibinfo {volume} {110}},\ \bibinfo {pages} {399--406}
  (\bibinfo {year} {1994})}\BibitemShut {NoStop}%
\bibitem [{\citenamefont {Vacondio}\ \emph {et~al.}(2020)\citenamefont
  {Vacondio}, \citenamefont {Altomare}, \citenamefont {De~Leffe}, \citenamefont
  {Hu}, \citenamefont {Le~Touzé}, \citenamefont {Lind}, \citenamefont
  {Marongiu}, \citenamefont {Marrone}, \citenamefont {Rogers},\ and\
  \citenamefont {Souto-Iglesias}}]{vacondio_grand_2020}%
  \BibitemOpen
  \bibfield  {author} {\bibinfo {author} {\bibfnamefont {R.}~\bibnamefont
  {Vacondio}}, \bibinfo {author} {\bibfnamefont {C.}~\bibnamefont {Altomare}},
  \bibinfo {author} {\bibfnamefont {M.}~\bibnamefont {De~Leffe}}, \bibinfo
  {author} {\bibfnamefont {X.}~\bibnamefont {Hu}}, \bibinfo {author}
  {\bibfnamefont {D.}~\bibnamefont {Le~Touzé}}, \bibinfo {author}
  {\bibfnamefont {S.}~\bibnamefont {Lind}}, \bibinfo {author} {\bibfnamefont
  {J.-C.}\ \bibnamefont {Marongiu}}, \bibinfo {author} {\bibfnamefont
  {S.}~\bibnamefont {Marrone}}, \bibinfo {author} {\bibfnamefont {B.~D.}\
  \bibnamefont {Rogers}}, \ and\ \bibinfo {author} {\bibfnamefont
  {A.}~\bibnamefont {Souto-Iglesias}},\ }\bibfield  {title} {\enquote {\bibinfo
  {title} {Grand challenges for {Smoothed} {Particle} {Hydrodynamics} numerical
  schemes},}\ }\href {\doibase 10.1007/s40571-020-00354-1} {\bibfield
  {journal} {\bibinfo  {journal} {Computational Particle Mechanics}\ }
  (\bibinfo {year} {2020}),\ 10.1007/s40571-020-00354-1}\BibitemShut {NoStop}%
\bibitem [{\citenamefont {Hernquist}\ and\ \citenamefont
  {Katz}(1989)}]{hernquist1989treesph}%
  \BibitemOpen
  \bibfield  {author} {\bibinfo {author} {\bibfnamefont {L.}~\bibnamefont
  {Hernquist}}\ and\ \bibinfo {author} {\bibfnamefont {N.}~\bibnamefont
  {Katz}},\ }\bibfield  {title} {\enquote {\bibinfo {title} {{TREESPH}-{A}
  unification of {SPH} with the hierarchical tree method},}\ }\href {\doibase
  10.1086/191344} {\bibfield  {journal} {\bibinfo  {journal} {The Astrophysical
  Journal Supplement Series}\ }\textbf {\bibinfo {volume} {70}},\ \bibinfo
  {pages} {419--446} (\bibinfo {year} {1989})}\BibitemShut {NoStop}%
\bibitem [{\citenamefont {Quinlan}, \citenamefont {Basa},\ and\ \citenamefont
  {Lastiwka}(2006)}]{quinlan_truncation_2006}%
  \BibitemOpen
  \bibfield  {author} {\bibinfo {author} {\bibfnamefont {N.~J.}\ \bibnamefont
  {Quinlan}}, \bibinfo {author} {\bibfnamefont {M.}~\bibnamefont {Basa}}, \
  and\ \bibinfo {author} {\bibfnamefont {M.}~\bibnamefont {Lastiwka}},\
  }\bibfield  {title} {\enquote {\bibinfo {title} {Truncation error in
  mesh-free particle methods},}\ }\href {\doibase 10.1002/nme.1617} {\bibfield
  {journal} {\bibinfo  {journal} {International Journal for Numerical Methods
  in Engineering}\ }\textbf {\bibinfo {volume} {66}},\ \bibinfo {pages}
  {2064--2085} (\bibinfo {year} {2006})}\BibitemShut {NoStop}%
\bibitem [{\citenamefont {Zhu}, \citenamefont {Hernquist},\ and\ \citenamefont
  {Li}(2015)}]{zhu2015numerical}%
  \BibitemOpen
  \bibfield  {author} {\bibinfo {author} {\bibfnamefont {Q.}~\bibnamefont
  {Zhu}}, \bibinfo {author} {\bibfnamefont {L.}~\bibnamefont {Hernquist}}, \
  and\ \bibinfo {author} {\bibfnamefont {Y.}~\bibnamefont {Li}},\ }\bibfield
  {title} {\enquote {\bibinfo {title} {Numerical convergence in smoothed
  particle hydrodynamics},}\ }\href {\doibase 10.1088/0004-637X/800/1/6}
  {\bibfield  {journal} {\bibinfo  {journal} {The Astrophysical Journal}\
  }\textbf {\bibinfo {volume} {800}},\ \bibinfo {pages} {6} (\bibinfo {year}
  {2015})}\BibitemShut {NoStop}%
\bibitem [{\citenamefont {Kiara}, \citenamefont {Hendrickson},\ and\
  \citenamefont {Yue}(2013{\natexlab{a}})}]{kiara_sph_2013}%
  \BibitemOpen
  \bibfield  {author} {\bibinfo {author} {\bibfnamefont {A.}~\bibnamefont
  {Kiara}}, \bibinfo {author} {\bibfnamefont {K.}~\bibnamefont {Hendrickson}},
  \ and\ \bibinfo {author} {\bibfnamefont {D.~K.}\ \bibnamefont {Yue}},\
  }\bibfield  {title} {\enquote {\bibinfo {title} {{SPH} for incompressible
  free-surface flows. {Part} {I}: {Error} analysis of the basic assumptions},}\
  }\href {\doibase 10.1016/j.compfluid.2013.05.023} {\bibfield  {journal}
  {\bibinfo  {journal} {Computers \& Fluids}\ }\textbf {\bibinfo {volume}
  {86}},\ \bibinfo {pages} {611--624} (\bibinfo {year}
  {2013}{\natexlab{a}})}\BibitemShut {NoStop}%
\bibitem [{\citenamefont {Kiara}, \citenamefont {Hendrickson},\ and\
  \citenamefont {Yue}(2013{\natexlab{b}})}]{kiara_sph_2013-1}%
  \BibitemOpen
  \bibfield  {author} {\bibinfo {author} {\bibfnamefont {A.}~\bibnamefont
  {Kiara}}, \bibinfo {author} {\bibfnamefont {K.}~\bibnamefont {Hendrickson}},
  \ and\ \bibinfo {author} {\bibfnamefont {D.~K.~P.}\ \bibnamefont {Yue}},\
  }\bibfield  {title} {\enquote {\bibinfo {title} {{SPH} for incompressible
  free-surface flows. {Part} {II}: {Performance} of a modified {SPH} method},}\
  }\href {\doibase 10.1016/j.compfluid.2013.07.016} {\bibfield  {journal}
  {\bibinfo  {journal} {Computers \& Fluids}\ }\textbf {\bibinfo {volume}
  {86}},\ \bibinfo {pages} {510--536} (\bibinfo {year}
  {2013}{\natexlab{b}})}\BibitemShut {NoStop}%
\bibitem [{\citenamefont {Lind}\ and\ \citenamefont
  {Stansby}(2016)}]{lind_high-order_2016}%
  \BibitemOpen
  \bibfield  {author} {\bibinfo {author} {\bibfnamefont {S.~J.}\ \bibnamefont
  {Lind}}\ and\ \bibinfo {author} {\bibfnamefont {P.~K.}\ \bibnamefont
  {Stansby}},\ }\bibfield  {title} {\enquote {\bibinfo {title} {High-order
  {Eulerian} incompressible smoothed particle hydrodynamics with transition to
  {Lagrangian} free-surface motion},}\ }\href {\doibase
  10.1016/j.jcp.2016.08.047} {\bibfield  {journal} {\bibinfo  {journal}
  {Journal of Computational Physics}\ }\textbf {\bibinfo {volume} {326}},\
  \bibinfo {pages} {290--311} (\bibinfo {year} {2016})}\BibitemShut {NoStop}%
\bibitem [{\citenamefont {Bonet}\ and\ \citenamefont
  {Lok}(1999)}]{bonet_lok:cmame:1999}%
  \BibitemOpen
  \bibfield  {author} {\bibinfo {author} {\bibfnamefont {J.}~\bibnamefont
  {Bonet}}\ and\ \bibinfo {author} {\bibfnamefont {T.-S.}\ \bibnamefont
  {Lok}},\ }\bibfield  {title} {\enquote {\bibinfo {title} {Variational and
  momentum preservation aspects of smooth particle hydrodynamic
  formulations},}\ }\href {\doibase 10.1016/S0045-7825(99)00051-1} {\bibfield
  {journal} {\bibinfo  {journal} {Computer Methods in Applied Mechanics and
  Engineering}\ }\textbf {\bibinfo {volume} {180}},\ \bibinfo {pages} {97 --
  115} (\bibinfo {year} {1999})}\BibitemShut {NoStop}%
\bibitem [{\citenamefont {Liu}\ and\ \citenamefont
  {Liu}(2006)}]{liu_restoring_2006}%
  \BibitemOpen
  \bibfield  {author} {\bibinfo {author} {\bibfnamefont {M.}~\bibnamefont
  {Liu}}\ and\ \bibinfo {author} {\bibfnamefont {G.}~\bibnamefont {Liu}},\
  }\bibfield  {title} {\enquote {\bibinfo {title} {Restoring particle
  consistency in smoothed particle hydrodynamics},}\ }\href {\doibase
  10.1016/j.apnum.2005.02.012} {\bibfield  {journal} {\bibinfo  {journal}
  {Applied Numerical Mathematics}\ }\textbf {\bibinfo {volume} {56}},\ \bibinfo
  {pages} {19--36} (\bibinfo {year} {2006})}\BibitemShut {NoStop}%
\bibitem [{\citenamefont {Rosswog}(2015)}]{rosswog2015boosting}%
  \BibitemOpen
  \bibfield  {author} {\bibinfo {author} {\bibfnamefont {S.}~\bibnamefont
  {Rosswog}},\ }\bibfield  {title} {\enquote {\bibinfo {title} {Boosting the
  accuracy of {SPH} techniques: Newtonian and special-relativistic tests},}\
  }\href {\doibase 10.1093/mnras/stv225} {\bibfield  {journal} {\bibinfo
  {journal} {Monthly Notices of the Royal Astronomical Society}\ }\textbf
  {\bibinfo {volume} {448}},\ \bibinfo {pages} {3628--3664} (\bibinfo {year}
  {2015})}\BibitemShut {NoStop}%
\bibitem [{\citenamefont {Huang}\ \emph {et~al.}(2019)\citenamefont {Huang},
  \citenamefont {Long}, \citenamefont {Li},\ and\ \citenamefont
  {Liu}}]{huang2019kernel}%
  \BibitemOpen
  \bibfield  {author} {\bibinfo {author} {\bibfnamefont {C.}~\bibnamefont
  {Huang}}, \bibinfo {author} {\bibfnamefont {T.}~\bibnamefont {Long}},
  \bibinfo {author} {\bibfnamefont {S.}~\bibnamefont {Li}}, \ and\ \bibinfo
  {author} {\bibfnamefont {M.}~\bibnamefont {Liu}},\ }\bibfield  {title}
  {\enquote {\bibinfo {title} {A kernel gradient-free {SPH} method with
  iterative particle shifting technology for modeling low-reynolds flows around
  airfoils},}\ }\href {\doibase 10.1016/j.enganabound.2019.06.010} {\bibfield
  {journal} {\bibinfo  {journal} {Engineering Analysis with Boundary Elements}\
  }\textbf {\bibinfo {volume} {106}},\ \bibinfo {pages} {571--587} (\bibinfo
  {year} {2019})}\BibitemShut {NoStop}%
\bibitem [{\citenamefont {Adami}, \citenamefont {Hu},\ and\ \citenamefont
  {Adams}(2013)}]{Adami2013}%
  \BibitemOpen
  \bibfield  {author} {\bibinfo {author} {\bibfnamefont {S.}~\bibnamefont
  {Adami}}, \bibinfo {author} {\bibfnamefont {X.}~\bibnamefont {Hu}}, \ and\
  \bibinfo {author} {\bibfnamefont {N.}~\bibnamefont {Adams}},\ }\bibfield
  {title} {\enquote {\bibinfo {title} {{A transport-velocity formulation for
  smoothed particle hydrodynamics}},}\ }\href {\doibase
  10.1016/j.jcp.2013.01.043} {\bibfield  {journal} {\bibinfo  {journal}
  {Journal of Computational Physics}\ }\textbf {\bibinfo {volume} {241}},\
  \bibinfo {pages} {292--307} (\bibinfo {year} {2013})}\BibitemShut {NoStop}%
\bibitem [{\citenamefont {Monaghan}\ and\ \citenamefont
  {Gingold}(1983)}]{monaghan1983shock}%
  \BibitemOpen
  \bibfield  {author} {\bibinfo {author} {\bibfnamefont {J.~J.}\ \bibnamefont
  {Monaghan}}\ and\ \bibinfo {author} {\bibfnamefont {R.~A.}\ \bibnamefont
  {Gingold}},\ }\bibfield  {title} {\enquote {\bibinfo {title} {Shock
  simulation by the particle method {SPH}},}\ }\href {\doibase
  10.1016/0021-9991(83)90036-0} {\bibfield  {journal} {\bibinfo  {journal}
  {Journal of computational physics}\ }\textbf {\bibinfo {volume} {52}},\
  \bibinfo {pages} {374--389} (\bibinfo {year} {1983})}\BibitemShut {NoStop}%
\bibitem [{\citenamefont {Violeau}(2012)}]{violeau2012fluid}%
  \BibitemOpen
  \bibfield  {author} {\bibinfo {author} {\bibfnamefont {D.}~\bibnamefont
  {Violeau}},\ }\href {\doibase 10.1093/acprof:oso/9780199655526.001.0001}
  {\emph {\bibinfo {title} {Fluid mechanics and the SPH method: theory and
  applications}}}\ (\bibinfo  {publisher} {Oxford University Press},\ \bibinfo
  {year} {2012})\BibitemShut {NoStop}%
\bibitem [{\citenamefont {Frontiere}, \citenamefont {Raskin},\ and\
  \citenamefont {Owen}(2017)}]{crksph:jcp:2017}%
  \BibitemOpen
  \bibfield  {author} {\bibinfo {author} {\bibfnamefont {N.}~\bibnamefont
  {Frontiere}}, \bibinfo {author} {\bibfnamefont {C.~D.}\ \bibnamefont
  {Raskin}}, \ and\ \bibinfo {author} {\bibfnamefont {J.~M.}\ \bibnamefont
  {Owen}},\ }\bibfield  {title} {\enquote {\bibinfo {title} {{CRKSPH} - a
  conservative reproducing kernel smoothed particle hydrodynamics scheme},}\
  }\href {\doibase 10.1016/j.jcp.2016.12.004} {\bibfield  {journal} {\bibinfo
  {journal} {Journal of Computational Physics}\ }\textbf {\bibinfo {volume}
  {332}},\ \bibinfo {pages} {160--209} (\bibinfo {year} {2017})}\BibitemShut
  {NoStop}%
\bibitem [{\citenamefont {Dilts}(1999)}]{dilts1999moving}%
  \BibitemOpen
  \bibfield  {author} {\bibinfo {author} {\bibfnamefont {G.~A.}\ \bibnamefont
  {Dilts}},\ }\bibfield  {title} {\enquote {\bibinfo {title}
  {Moving-least-squares-particle hydrodynamics—{I}. consistency and
  stability},}\ }\href {\doibase
  10.1002/(SICI)1097-0207(19990320)44:8<1115::AID-NME547>3.0.CO;2-L} {\bibfield
   {journal} {\bibinfo  {journal} {International Journal for Numerical Methods
  in Engineering}\ }\textbf {\bibinfo {volume} {44}},\ \bibinfo {pages}
  {1115--1155} (\bibinfo {year} {1999})}\BibitemShut {NoStop}%
\bibitem [{\citenamefont {Dilts}(2000)}]{dilts2000moving}%
  \BibitemOpen
  \bibfield  {author} {\bibinfo {author} {\bibfnamefont {G.~A.}\ \bibnamefont
  {Dilts}},\ }\bibfield  {title} {\enquote {\bibinfo {title} {Moving
  least-squares particle hydrodynamics {II}: conservation and boundaries},}\
  }\href {\doibase
  10.1002/1097-0207(20000810)48:10<1503::AID-NME832>3.0.CO;2-D} {\bibfield
  {journal} {\bibinfo  {journal} {International Journal for numerical methods
  in engineering}\ }\textbf {\bibinfo {volume} {48}},\ \bibinfo {pages}
  {1503--1524} (\bibinfo {year} {2000})}\BibitemShut {NoStop}%
\bibitem [{\citenamefont {Chen}\ and\ \citenamefont
  {Beraun}(2000)}]{chenGeneralizedSmoothedParticle2000a}%
  \BibitemOpen
  \bibfield  {author} {\bibinfo {author} {\bibfnamefont {J.~K.}\ \bibnamefont
  {Chen}}\ and\ \bibinfo {author} {\bibfnamefont {J.~E.}\ \bibnamefont
  {Beraun}},\ }\bibfield  {title} {\enquote {\bibinfo {title} {A generalized
  smoothed particle hydrodynamics method for nonlinear dynamic problems},}\
  }\href {\doibase 10.1016/S0045-7825(99)00422-3} {\bibfield  {journal}
  {\bibinfo  {journal} {Computer Methods in Applied Mechanics and Engineering}\
  }\textbf {\bibinfo {volume} {190}},\ \bibinfo {pages} {225--239} (\bibinfo
  {year} {2000})}\BibitemShut {NoStop}%
\bibitem [{\citenamefont {Zhang}\ and\ \citenamefont
  {Batra}(2004)}]{zhangModifiedSmoothedParticle2004}%
  \BibitemOpen
  \bibfield  {author} {\bibinfo {author} {\bibfnamefont {G.~M.}\ \bibnamefont
  {Zhang}}\ and\ \bibinfo {author} {\bibfnamefont {R.~C.}\ \bibnamefont
  {Batra}},\ }\bibfield  {title} {\enquote {\bibinfo {title} {Modified smoothed
  particle hydrodynamics method and its application to transient problems},}\
  }\href {\doibase 10.1007/s00466-004-0561-5} {\bibfield  {journal} {\bibinfo
  {journal} {Computational Mechanics}\ }\textbf {\bibinfo {volume} {34}},\
  \bibinfo {pages} {137--146} (\bibinfo {year} {2004})}\BibitemShut {NoStop}%
\bibitem [{\citenamefont {Korzilius}, \citenamefont {Schilders},\ and\
  \citenamefont {Anthonissen}(2017)}]{korziliusImprovedCSPMApproach2017}%
  \BibitemOpen
  \bibfield  {author} {\bibinfo {author} {\bibfnamefont {S.~P.}\ \bibnamefont
  {Korzilius}}, \bibinfo {author} {\bibfnamefont {W.~H.~A.}\ \bibnamefont
  {Schilders}}, \ and\ \bibinfo {author} {\bibfnamefont {M.~J.~H.}\
  \bibnamefont {Anthonissen}},\ }\bibfield  {title} {\enquote {\bibinfo {title}
  {An {{Improved CSPM Approach}} for {{Accurate Second}}-{{Derivative
  Approximations}} with {{SPH}}},}\ }\href {\doibase 10.4236/jamp.2017.51017}
  {\bibfield  {journal} {\bibinfo  {journal} {Journal of Applied Mathematics
  and Physics}\ }\textbf {\bibinfo {volume} {05}},\ \bibinfo {pages} {168--184}
  (\bibinfo {year} {2017})}\BibitemShut {NoStop}%
\bibitem [{\citenamefont
  {Schwaiger}(2008)}]{schwaigerImplicitCorrectedSPH2008}%
  \BibitemOpen
  \bibfield  {author} {\bibinfo {author} {\bibfnamefont {H.~F.}\ \bibnamefont
  {Schwaiger}},\ }\bibfield  {title} {\enquote {\bibinfo {title} {An implicit
  corrected {{SPH}} formulation for thermal diffusion with linear free surface
  boundary conditions},}\ }\href {\doibase 10.1002/nme.2266} {\bibfield
  {journal} {\bibinfo  {journal} {International Journal for Numerical Methods
  in Engineering}\ }\textbf {\bibinfo {volume} {75}},\ \bibinfo {pages}
  {647--671} (\bibinfo {year} {2008})}\BibitemShut {NoStop}%
\bibitem [{\citenamefont {Maci{\`a}}\ \emph {et~al.}(2012)\citenamefont
  {Maci{\`a}}, \citenamefont {Gonz{\'a}lez}, \citenamefont {{Cercos-Pita}},\
  and\ \citenamefont {{Souto-Iglesias}}}]{maciaBoundaryIntegralSPH2012}%
  \BibitemOpen
  \bibfield  {author} {\bibinfo {author} {\bibfnamefont {F.}~\bibnamefont
  {Maci{\`a}}}, \bibinfo {author} {\bibfnamefont {L.~M.}\ \bibnamefont
  {Gonz{\'a}lez}}, \bibinfo {author} {\bibfnamefont {J.~L.}\ \bibnamefont
  {{Cercos-Pita}}}, \ and\ \bibinfo {author} {\bibfnamefont {A.}~\bibnamefont
  {{Souto-Iglesias}}},\ }\bibfield  {title} {\enquote {\bibinfo {title} {A
  {{Boundary Integral SPH Formulation}}: {{Consistency}} and {{Applications}}
  to {{ISPH}} and {{WCSPH}}},}\ }\href {\doibase 10.1143/PTP.128.439}
  {\bibfield  {journal} {\bibinfo  {journal} {Progress of Theoretical Physics}\
  }\textbf {\bibinfo {volume} {128}},\ \bibinfo {pages} {439--462} (\bibinfo
  {year} {2012})}\BibitemShut {NoStop}%
\bibitem [{\citenamefont {Brookshaw}(1985)}]{brookshaw_method_1985}%
  \BibitemOpen
  \bibfield  {author} {\bibinfo {author} {\bibfnamefont {L.}~\bibnamefont
  {Brookshaw}},\ }\bibfield  {title} {\enquote {\bibinfo {title} {A {Method} of
  {Calculating} {Radiative} {Heat} {Diffusion} in {Particle} {Simulations}},}\
  }\href {\doibase 10.1017/S1323358000018117} {\bibfield  {journal} {\bibinfo
  {journal} {Publications of the Astronomical Society of Australia}\ }\textbf
  {\bibinfo {volume} {6}},\ \bibinfo {pages} {207--210} (\bibinfo {year}
  {1985})}\BibitemShut {NoStop}%
\bibitem [{\citenamefont {Morris}, \citenamefont {Fox},\ and\ \citenamefont
  {Zhu}(1997)}]{morris-lowRe-97}%
  \BibitemOpen
  \bibfield  {author} {\bibinfo {author} {\bibfnamefont {J.~P.}\ \bibnamefont
  {Morris}}, \bibinfo {author} {\bibfnamefont {P.~J.}\ \bibnamefont {Fox}}, \
  and\ \bibinfo {author} {\bibfnamefont {Y.}~\bibnamefont {Zhu}},\ }\bibfield
  {title} {\enquote {\bibinfo {title} {Modeling low reynolds number
  incompressible flows using {SPH}},}\ }\href {\doibase 10.1006/jcph.1997.5776}
  {\bibfield  {journal} {\bibinfo  {journal} {Journal of Computational
  Physics}\ }\textbf {\bibinfo {volume} {136}},\ \bibinfo {pages} {214--226}
  (\bibinfo {year} {1997})}\BibitemShut {NoStop}%
\bibitem [{\citenamefont {Cleary}\ and\ \citenamefont
  {Monaghan}(1999)}]{cleary1999conduction}%
  \BibitemOpen
  \bibfield  {author} {\bibinfo {author} {\bibfnamefont {P.~W.}\ \bibnamefont
  {Cleary}}\ and\ \bibinfo {author} {\bibfnamefont {J.~J.}\ \bibnamefont
  {Monaghan}},\ }\bibfield  {title} {\enquote {\bibinfo {title} {Conduction
  modelling using smoothed particle hydrodynamics},}\ }\href {\doibase
  10.1006/jcph.1998.6118} {\bibfield  {journal} {\bibinfo  {journal} {Journal
  of Computational Physics}\ }\textbf {\bibinfo {volume} {148}},\ \bibinfo
  {pages} {227--264} (\bibinfo {year} {1999})}\BibitemShut {NoStop}%
\bibitem [{\citenamefont {Fatehi}\ and\ \citenamefont
  {Manzari}(2011)}]{fatehi_error_2011}%
  \BibitemOpen
  \bibfield  {author} {\bibinfo {author} {\bibfnamefont {R.}~\bibnamefont
  {Fatehi}}\ and\ \bibinfo {author} {\bibfnamefont {M.}~\bibnamefont
  {Manzari}},\ }\bibfield  {title} {\enquote {\bibinfo {title} {Error
  estimation in smoothed particle hydrodynamics and a new scheme for second
  derivatives},}\ }\href {\doibase 10.1016/j.camwa.2010.11.028} {\bibfield
  {journal} {\bibinfo  {journal} {Computers \& Mathematics with Applications}\
  }\textbf {\bibinfo {volume} {61}},\ \bibinfo {pages} {482--498} (\bibinfo
  {year} {2011})}\BibitemShut {NoStop}%
\bibitem [{\citenamefont {Nugent}\ and\ \citenamefont
  {Posch}(2000)}]{nugent_liquid_2000}%
  \BibitemOpen
  \bibfield  {author} {\bibinfo {author} {\bibfnamefont {S.}~\bibnamefont
  {Nugent}}\ and\ \bibinfo {author} {\bibfnamefont {H.~A.}\ \bibnamefont
  {Posch}},\ }\bibfield  {title} {\enquote {\bibinfo {title} {Liquid drops and
  surface tension with smoothed particle applied mechanics},}\ }\href {\doibase
  10.1103/PhysRevE.62.4968} {\bibfield  {journal} {\bibinfo  {journal}
  {Physical Review E}\ }\textbf {\bibinfo {volume} {62}},\ \bibinfo {pages}
  {4968--4975} (\bibinfo {year} {2000})},\ \bibinfo {note} {publisher: American
  Physical Society}\BibitemShut {NoStop}%
\bibitem [{\citenamefont {Biriukov}\ and\ \citenamefont
  {Price}(2019)}]{biriukovStableAnisotropicHeat2019a}%
  \BibitemOpen
  \bibfield  {author} {\bibinfo {author} {\bibfnamefont {S.}~\bibnamefont
  {Biriukov}}\ and\ \bibinfo {author} {\bibfnamefont {D.~J.}\ \bibnamefont
  {Price}},\ }\bibfield  {title} {\enquote {\bibinfo {title} {Stable
  anisotropic heat conduction in smoothed particle hydrodynamics},}\ }\href
  {\doibase 10.1093/mnras/sty3413} {\bibfield  {journal} {\bibinfo  {journal}
  {Monthly Notices of the Royal Astronomical Society}\ }\textbf {\bibinfo
  {volume} {483}},\ \bibinfo {pages} {4901--4909} (\bibinfo {year}
  {2019})}\BibitemShut {NoStop}%
\bibitem [{\citenamefont {Gomez-Gesteria}\ \emph {et~al.}(2010)\citenamefont
  {Gomez-Gesteria}, \citenamefont {Rogers}, \citenamefont {Dalrymple},\ and\
  \citenamefont {Crespo}}]{wcsph-state-of-the-art-2010}%
  \BibitemOpen
  \bibfield  {author} {\bibinfo {author} {\bibfnamefont {M.}~\bibnamefont
  {Gomez-Gesteria}}, \bibinfo {author} {\bibfnamefont {B.~D.}\ \bibnamefont
  {Rogers}}, \bibinfo {author} {\bibfnamefont {R.~A.}\ \bibnamefont
  {Dalrymple}}, \ and\ \bibinfo {author} {\bibfnamefont {A.~J.}\ \bibnamefont
  {Crespo}},\ }\bibfield  {title} {\enquote {\bibinfo {title} {State-of-the-art
  classical {SPH} for free-surface flows},}\ }\href {\doibase
  10.1080/00221686.2010.9641242} {\bibfield  {journal} {\bibinfo  {journal}
  {Journal of Hydraulic Research}\ }\textbf {\bibinfo {volume} {84}},\ \bibinfo
  {pages} {6--27} (\bibinfo {year} {2010})}\BibitemShut {NoStop}%
\bibitem [{\citenamefont {Xu}, \citenamefont {Stansby},\ and\ \citenamefont
  {Laurence}(2009)}]{acc_stab_xu:jcp:2009}%
  \BibitemOpen
  \bibfield  {author} {\bibinfo {author} {\bibfnamefont {R.}~\bibnamefont
  {Xu}}, \bibinfo {author} {\bibfnamefont {P.}~\bibnamefont {Stansby}}, \ and\
  \bibinfo {author} {\bibfnamefont {D.}~\bibnamefont {Laurence}},\ }\bibfield
  {title} {\enquote {\bibinfo {title} {Accuracy and stability in incompressible
  sph ({ISPH}) based on the projection method and a new approach},}\ }\href
  {\doibase 10.1016/j.jcp.2009.05.032} {\bibfield  {journal} {\bibinfo
  {journal} {Journal of Computational Physics}\ }\textbf {\bibinfo {volume}
  {228}},\ \bibinfo {pages} {6703--6725} (\bibinfo {year} {2009})}\BibitemShut
  {NoStop}%
\bibitem [{\citenamefont {Lind}\ \emph
  {et~al.}(2012{\natexlab{a}})\citenamefont {Lind}, \citenamefont {Xu},
  \citenamefont {Stansby},\ and\ \citenamefont
  {Rogers}}]{lind2012incompressible}%
  \BibitemOpen
  \bibfield  {author} {\bibinfo {author} {\bibfnamefont {S.~J.}\ \bibnamefont
  {Lind}}, \bibinfo {author} {\bibfnamefont {R.}~\bibnamefont {Xu}}, \bibinfo
  {author} {\bibfnamefont {P.~K.}\ \bibnamefont {Stansby}}, \ and\ \bibinfo
  {author} {\bibfnamefont {B.~D.}\ \bibnamefont {Rogers}},\ }\bibfield  {title}
  {\enquote {\bibinfo {title} {Incompressible smoothed particle hydrodynamics
  for free-surface flows: A generalised diffusion-based algorithm for stability
  and validations for impulsive flows and propagating waves},}\ }\href
  {\doibase 10.1016/j.jcp.2011.10.027} {\bibfield  {journal} {\bibinfo
  {journal} {Journal of Computational Physics}\ }\textbf {\bibinfo {volume}
  {231}},\ \bibinfo {pages} {1499--1523} (\bibinfo {year}
  {2012}{\natexlab{a}})}\BibitemShut {NoStop}%
\bibitem [{\citenamefont {Sun}\ \emph {et~al.}(2019{\natexlab{a}})\citenamefont
  {Sun}, \citenamefont {Colagrossi}, \citenamefont {Marrone}, \citenamefont
  {Antuono},\ and\ \citenamefont {Zhang}}]{sun_consistent_2019}%
  \BibitemOpen
  \bibfield  {author} {\bibinfo {author} {\bibfnamefont {P.}~\bibnamefont
  {Sun}}, \bibinfo {author} {\bibfnamefont {A.}~\bibnamefont {Colagrossi}},
  \bibinfo {author} {\bibfnamefont {S.}~\bibnamefont {Marrone}}, \bibinfo
  {author} {\bibfnamefont {M.}~\bibnamefont {Antuono}}, \ and\ \bibinfo
  {author} {\bibfnamefont {A.-M.}\ \bibnamefont {Zhang}},\ }\bibfield  {title}
  {\enquote {\bibinfo {title} {A consistent approach to particle shifting in
  the $\delta$ - {Plus} -{SPH} model},}\ }\href {\doibase
  10.1016/j.cma.2019.01.045} {\bibfield  {journal} {\bibinfo  {journal}
  {Computer Methods in Applied Mechanics and Engineering}\ }\textbf {\bibinfo
  {volume} {348}},\ \bibinfo {pages} {912--934} (\bibinfo {year}
  {2019}{\natexlab{a}})}\BibitemShut {NoStop}%
\bibitem [{\citenamefont {Oger}\ \emph {et~al.}(2016)\citenamefont {Oger},
  \citenamefont {Marrone}, \citenamefont {Le~Touzé},\ and\ \citenamefont
  {de~Leffe}}]{oger_ale_sph_2016}%
  \BibitemOpen
  \bibfield  {author} {\bibinfo {author} {\bibfnamefont {G.}~\bibnamefont
  {Oger}}, \bibinfo {author} {\bibfnamefont {S.}~\bibnamefont {Marrone}},
  \bibinfo {author} {\bibfnamefont {D.}~\bibnamefont {Le~Touzé}}, \ and\
  \bibinfo {author} {\bibfnamefont {M.}~\bibnamefont {de~Leffe}},\ }\bibfield
  {title} {\enquote {\bibinfo {title} {{SPH} accuracy improvement through the
  combination of a quasi-{Lagrangian} shifting transport velocity and
  consistent {ALE} formalisms},}\ }\href {\doibase 10.1016/j.jcp.2016.02.039}
  {\bibfield  {journal} {\bibinfo  {journal} {Journal of Computational
  Physics}\ }\textbf {\bibinfo {volume} {313}},\ \bibinfo {pages} {76--98}
  (\bibinfo {year} {2016})}\BibitemShut {NoStop}%
\bibitem [{\citenamefont {{Marrone}}\ \emph {et~al.}(2011)\citenamefont
  {{Marrone}}, \citenamefont {{Antuono}}, \citenamefont {{Colagrossi}},
  \citenamefont {{Colicchio}}, \citenamefont {{Le Touz{\'e}}},\ and\
  \citenamefont {{Graziani}}}]{marrone-deltasph:cmame:2011}%
  \BibitemOpen
  \bibfield  {author} {\bibinfo {author} {\bibfnamefont {S.}~\bibnamefont
  {{Marrone}}}, \bibinfo {author} {\bibfnamefont {M.}~\bibnamefont
  {{Antuono}}}, \bibinfo {author} {\bibfnamefont {A.}~\bibnamefont
  {{Colagrossi}}}, \bibinfo {author} {\bibfnamefont {G.}~\bibnamefont
  {{Colicchio}}}, \bibinfo {author} {\bibfnamefont {D.}~\bibnamefont {{Le
  Touz{\'e}}}}, \ and\ \bibinfo {author} {\bibfnamefont {G.}~\bibnamefont
  {{Graziani}}},\ }\bibfield  {title} {\enquote {\bibinfo {title}
  {{$\delta$-SPH} model for simulating violent impact flows},}\ }\href
  {\doibase 10.1016/j.cma.2010.12.016} {\bibfield  {journal} {\bibinfo
  {journal} {Computer Methods in Applied Mechanics and Engineering}\ }\textbf
  {\bibinfo {volume} {200}},\ \bibinfo {pages} {1526--1542} (\bibinfo {year}
  {2011})}\BibitemShut {NoStop}%
\bibitem [{\citenamefont {Antuono}\ \emph {et~al.}(2010)\citenamefont
  {Antuono}, \citenamefont {Colagrossi}, \citenamefont {Marrone},\ and\
  \citenamefont {Molteni}}]{antuono-deltasph:cpc:2010}%
  \BibitemOpen
  \bibfield  {author} {\bibinfo {author} {\bibfnamefont {M.}~\bibnamefont
  {Antuono}}, \bibinfo {author} {\bibfnamefont {A.}~\bibnamefont {Colagrossi}},
  \bibinfo {author} {\bibfnamefont {S.}~\bibnamefont {Marrone}}, \ and\
  \bibinfo {author} {\bibfnamefont {D.}~\bibnamefont {Molteni}},\ }\bibfield
  {title} {\enquote {\bibinfo {title} {Free-surface flows solved by means of
  {SPH} schemes with numerical diffusive terms},}\ }\href {\doibase
  10.1016/j.cpc.2009.11.002} {\bibfield  {journal} {\bibinfo  {journal}
  {Computer Physics Communications}\ }\textbf {\bibinfo {volume} {181}},\
  \bibinfo {pages} {532 -- 549} (\bibinfo {year} {2010})}\BibitemShut {NoStop}%
\bibitem [{\citenamefont {Ramachandran}\ and\ \citenamefont
  {Puri}(2019)}]{edac-sph:cf:2019}%
  \BibitemOpen
  \bibfield  {author} {\bibinfo {author} {\bibfnamefont {P.}~\bibnamefont
  {Ramachandran}}\ and\ \bibinfo {author} {\bibfnamefont {K.}~\bibnamefont
  {Puri}},\ }\bibfield  {title} {\enquote {\bibinfo {title} {Entropically
  damped artificial compressibility for {SPH}},}\ }\href {\doibase
  10.1016/j.compfluid.2018.11.023} {\bibfield  {journal} {\bibinfo  {journal}
  {Computers and Fluids}\ }\textbf {\bibinfo {volume} {179}},\ \bibinfo {pages}
  {579--594} (\bibinfo {year} {2019})}\BibitemShut {NoStop}%
\bibitem [{\citenamefont {Clausen}(2013)}]{Clausen2013}%
  \BibitemOpen
  \bibfield  {author} {\bibinfo {author} {\bibfnamefont {J.~R.}\ \bibnamefont
  {Clausen}},\ }\bibfield  {title} {\enquote {\bibinfo {title} {{Entropically
  damped form of artificial compressibility for explicit simulation of
  incompressible flow}},}\ }\href {\doibase 10.1103/PhysRevE.87.013309}
  {\bibfield  {journal} {\bibinfo  {journal} {Physical Review E}\ }\textbf
  {\bibinfo {volume} {87}},\ \bibinfo {pages} {013309--1--013309--12} (\bibinfo
  {year} {2013})}\BibitemShut {NoStop}%
\bibitem [{\citenamefont {Chaniotis}, \citenamefont {Poulikakos},\ and\
  \citenamefont {Koumoutsakos}(2002)}]{remeshed_sph:jcp:2002}%
  \BibitemOpen
  \bibfield  {author} {\bibinfo {author} {\bibfnamefont {A.}~\bibnamefont
  {Chaniotis}}, \bibinfo {author} {\bibfnamefont {D.}~\bibnamefont
  {Poulikakos}}, \ and\ \bibinfo {author} {\bibfnamefont {P.}~\bibnamefont
  {Koumoutsakos}},\ }\bibfield  {title} {\enquote {\bibinfo {title} {Remeshed
  smoothed particle hydrodynamics for the simulation of viscous and heat
  conducting flows},}\ }\href {\doibase 10.1006/jcph.2002.7152} {\bibfield
  {journal} {\bibinfo  {journal} {Journal of Computational Physics}\ }\textbf
  {\bibinfo {volume} {182}},\ \bibinfo {pages} {67 -- 90} (\bibinfo {year}
  {2002})}\BibitemShut {NoStop}%
\bibitem [{\citenamefont {Hieber}\ and\ \citenamefont
  {Koumoutsakos}(2008)}]{hieberImmersedBoundaryMethod2008}%
  \BibitemOpen
  \bibfield  {author} {\bibinfo {author} {\bibfnamefont {S.~E.}\ \bibnamefont
  {Hieber}}\ and\ \bibinfo {author} {\bibfnamefont {P.}~\bibnamefont
  {Koumoutsakos}},\ }\bibfield  {title} {\enquote {\bibinfo {title} {An
  immersed boundary method for smoothed particle hydrodynamics of
  self-propelled swimmers},}\ }\href {\doibase 10.1016/j.jcp.2008.06.017}
  {\bibfield  {journal} {\bibinfo  {journal} {Journal of Computational
  Physics}\ }\textbf {\bibinfo {volume} {227}},\ \bibinfo {pages} {8636--8654}
  (\bibinfo {year} {2008})}\BibitemShut {NoStop}%
\bibitem [{\citenamefont {Nasar}\ \emph {et~al.}(2019)\citenamefont {Nasar},
  \citenamefont {Rogers}, \citenamefont {Revell}, \citenamefont {Stansby},\
  and\ \citenamefont {Lind}}]{nasar2019}%
  \BibitemOpen
  \bibfield  {author} {\bibinfo {author} {\bibfnamefont {A.}~\bibnamefont
  {Nasar}}, \bibinfo {author} {\bibfnamefont {B.}~\bibnamefont {Rogers}},
  \bibinfo {author} {\bibfnamefont {A.}~\bibnamefont {Revell}}, \bibinfo
  {author} {\bibfnamefont {P.}~\bibnamefont {Stansby}}, \ and\ \bibinfo
  {author} {\bibfnamefont {S.}~\bibnamefont {Lind}},\ }\bibfield  {title}
  {\enquote {\bibinfo {title} {Eulerian weakly compressible smoothed particle
  hydrodynamics ({{SPH}}) with the immersed boundary method for thin slender
  bodies},}\ }\href {\doibase 10.1016/j.jfluidstructs.2018.11.005} {\bibfield
  {journal} {\bibinfo  {journal} {Journal of Fluids and Structures}\ }\textbf
  {\bibinfo {volume} {84}},\ \bibinfo {pages} {263--282} (\bibinfo {year}
  {2019})}\BibitemShut {NoStop}%
\bibitem [{\citenamefont {Sun}\ \emph {et~al.}(2019{\natexlab{b}})\citenamefont
  {Sun}, \citenamefont {Colagrossi}, \citenamefont {Marrone}, \citenamefont
  {Antuono},\ and\ \citenamefont {Zhang}}]{sun2019consistent}%
  \BibitemOpen
  \bibfield  {author} {\bibinfo {author} {\bibfnamefont {P.}~\bibnamefont
  {Sun}}, \bibinfo {author} {\bibfnamefont {A.}~\bibnamefont {Colagrossi}},
  \bibinfo {author} {\bibfnamefont {S.}~\bibnamefont {Marrone}}, \bibinfo
  {author} {\bibfnamefont {M.}~\bibnamefont {Antuono}}, \ and\ \bibinfo
  {author} {\bibfnamefont {A.-M.}\ \bibnamefont {Zhang}},\ }\bibfield  {title}
  {\enquote {\bibinfo {title} {A consistent approach to particle shifting in
  the $\delta$-plus-{SPH} model},}\ }\href {\doibase 10.1016/j.cma.2019.01.045}
  {\bibfield  {journal} {\bibinfo  {journal} {Computer Methods in Applied
  Mechanics and Engineering}\ }\textbf {\bibinfo {volume} {348}},\ \bibinfo
  {pages} {912--934} (\bibinfo {year} {2019}{\natexlab{b}})}\BibitemShut
  {NoStop}%
\bibitem [{\citenamefont {Lind}\ \emph
  {et~al.}(2012{\natexlab{b}})\citenamefont {Lind}, \citenamefont {Xu},
  \citenamefont {Stansby},\ and\ \citenamefont
  {Rogers}}]{diff_smoothing_sph:lind:jcp:2009}%
  \BibitemOpen
  \bibfield  {author} {\bibinfo {author} {\bibfnamefont {S.}~\bibnamefont
  {Lind}}, \bibinfo {author} {\bibfnamefont {R.}~\bibnamefont {Xu}}, \bibinfo
  {author} {\bibfnamefont {P.}~\bibnamefont {Stansby}}, \ and\ \bibinfo
  {author} {\bibfnamefont {B.}~\bibnamefont {Rogers}},\ }\bibfield  {title}
  {\enquote {\bibinfo {title} {Incompressible smoothed particle hydrodynamics
  for free-surface flows: A generalised diffusion-based algorithm for stability
  and validations for impulsive flows and propagating waves},}\ }\href
  {\doibase 10.1016/j.jcp.2011.10.027} {\bibfield  {journal} {\bibinfo
  {journal} {Journal of Computational Physics}\ }\textbf {\bibinfo {volume}
  {231}},\ \bibinfo {pages} {1499 -- 1523} (\bibinfo {year}
  {2012}{\natexlab{b}})}\BibitemShut {NoStop}%
\bibitem [{\citenamefont {Dehnen}\ and\ \citenamefont
  {Aly}(2012)}]{dehnen-aly-paring-instability-mnras-2012}%
  \BibitemOpen
  \bibfield  {author} {\bibinfo {author} {\bibfnamefont {W.}~\bibnamefont
  {Dehnen}}\ and\ \bibinfo {author} {\bibfnamefont {H.}~\bibnamefont {Aly}},\
  }\bibfield  {title} {\enquote {\bibinfo {title} {Improving convergence in
  smoothed particle hydrodynamics simulations without pairing instability},}\
  }\href {\doibase 10.1111/j.1365-2966.2012.21439.x} {\bibfield  {journal}
  {\bibinfo  {journal} {Monthly Notices of the Royal Astronomical Society}\
  }\textbf {\bibinfo {volume} {425}},\ \bibinfo {pages} {1068--1082} (\bibinfo
  {year} {2012})}\BibitemShut {NoStop}%
\bibitem [{\citenamefont {Gresho}\ and\ \citenamefont
  {Chan}(1990)}]{gresho1990theory}%
  \BibitemOpen
  \bibfield  {author} {\bibinfo {author} {\bibfnamefont {P.~M.}\ \bibnamefont
  {Gresho}}\ and\ \bibinfo {author} {\bibfnamefont {S.~T.}\ \bibnamefont
  {Chan}},\ }\bibfield  {title} {\enquote {\bibinfo {title} {On the theory of
  semi-implicit projection methods for viscous incompressible flow and its
  implementation via a finite element method that also introduces a nearly
  consistent mass matrix. part 2: Implementation},}\ }\href {\doibase
  10.1002/fld.1650110510} {\bibfield  {journal} {\bibinfo  {journal}
  {International journal for numerical methods in fluids}\ }\textbf {\bibinfo
  {volume} {11}},\ \bibinfo {pages} {621--659} (\bibinfo {year}
  {1990})}\BibitemShut {NoStop}%
\bibitem [{\citenamefont {Di}\ \emph {et~al.}(2005)\citenamefont {Di},
  \citenamefont {Li}, \citenamefont {Tang},\ and\ \citenamefont
  {Zhang}}]{diMovingMeshFinite2005}%
  \BibitemOpen
  \bibfield  {author} {\bibinfo {author} {\bibfnamefont {Y.}~\bibnamefont
  {Di}}, \bibinfo {author} {\bibfnamefont {R.}~\bibnamefont {Li}}, \bibinfo
  {author} {\bibfnamefont {T.}~\bibnamefont {Tang}}, \ and\ \bibinfo {author}
  {\bibfnamefont {P.}~\bibnamefont {Zhang}},\ }\bibfield  {title} {\enquote
  {\bibinfo {title} {Moving {{Mesh Finite Element Methods}} for the
  {{Incompressible Navier}}--{{Stokes Equations}}},}\ }\href {\doibase
  10.1137/030600643} {\bibfield  {journal} {\bibinfo  {journal} {SIAM Journal
  on Scientific Computing}\ }\textbf {\bibinfo {volume} {26}},\ \bibinfo
  {pages} {1036--1056} (\bibinfo {year} {2005})}\BibitemShut {NoStop}%
\bibitem [{\citenamefont {Ramachandran}(2018)}]{pr:automan:2018}%
  \BibitemOpen
  \bibfield  {author} {\bibinfo {author} {\bibfnamefont {P.}~\bibnamefont
  {Ramachandran}},\ }\bibfield  {title} {\enquote {\bibinfo {title} {automan: A
  python-based automation framework for numerical computing},}\ }\href
  {\doibase 10.1109/MCSE.2018.05329818} {\bibfield  {journal} {\bibinfo
  {journal} {Computing in Science \& Engineering}\ }\textbf {\bibinfo {volume}
  {20}},\ \bibinfo {pages} {81--97} (\bibinfo {year} {2018})}\BibitemShut
  {NoStop}%
\bibitem [{\citenamefont {Monaghan}(2005)}]{monaghan-review:2005}%
  \BibitemOpen
  \bibfield  {author} {\bibinfo {author} {\bibfnamefont {J.~J.}\ \bibnamefont
  {Monaghan}},\ }\bibfield  {title} {\enquote {\bibinfo {title} {{Smoothed
  Particle Hydrodynamics}},}\ }\href {\doibase 10.1088/0034-4885/68/8/R01}
  {\bibfield  {journal} {\bibinfo  {journal} {{Reports on Progress in
  Physics}}\ }\textbf {\bibinfo {volume} {68}},\ \bibinfo {pages} {1703--1759}
  (\bibinfo {year} {2005})}\BibitemShut {NoStop}%
\bibitem [{\citenamefont {Wendland}(1995)}]{wendland_piecewise_1995}%
  \BibitemOpen
  \bibfield  {author} {\bibinfo {author} {\bibfnamefont {H.}~\bibnamefont
  {Wendland}},\ }\bibfield  {title} {\enquote {\bibinfo {title} {Piecewise
  polynomial, positive definite and compactly supported radial functions of
  minimal degree},}\ }\href {\doibase 10.1007/BF02123482} {\bibfield  {journal}
  {\bibinfo  {journal} {Advances in Computational Mathematics}\ }\textbf
  {\bibinfo {volume} {4}},\ \bibinfo {pages} {389--396} (\bibinfo {year}
  {1995})}\BibitemShut {NoStop}%
\bibitem [{\citenamefont {Randles}\ and\ \citenamefont
  {Libersky}(1996)}]{randles1996smoothed}%
  \BibitemOpen
  \bibfield  {author} {\bibinfo {author} {\bibfnamefont {P.}~\bibnamefont
  {Randles}}\ and\ \bibinfo {author} {\bibfnamefont {L.~D.}\ \bibnamefont
  {Libersky}},\ }\bibfield  {title} {\enquote {\bibinfo {title} {Smoothed
  particle hydrodynamics: some recent improvements and applications},}\ }\href
  {\doibase 10.1016/S0045-7825(96)01090-0} {\bibfield  {journal} {\bibinfo
  {journal} {Computer methods in applied mechanics and engineering}\ }\textbf
  {\bibinfo {volume} {139}},\ \bibinfo {pages} {375--408} (\bibinfo {year}
  {1996})}\BibitemShut {NoStop}%
\bibitem [{\citenamefont {Negi}\ and\ \citenamefont
  {Ramachandran}(2021)}]{negi2019improved}%
  \BibitemOpen
  \bibfield  {author} {\bibinfo {author} {\bibfnamefont {P.}~\bibnamefont
  {Negi}}\ and\ \bibinfo {author} {\bibfnamefont {P.}~\bibnamefont
  {Ramachandran}},\ }\bibfield  {title} {\enquote {\bibinfo {title} {Algorithms
  for uniform particle initialization in domains with complex boundaries},}\
  }\href {\doibase 10.1016/j.cpc.2021.108008} {\bibfield  {journal} {\bibinfo
  {journal} {Computer Physics Communications}\ }\textbf {\bibinfo {volume}
  {265}},\ \bibinfo {pages} {108008} (\bibinfo {year} {2021})}\BibitemShut
  {NoStop}%
\bibitem [{\citenamefont {Colagrossi}\ \emph {et~al.}(2012)\citenamefont
  {Colagrossi}, \citenamefont {Bouscasse}, \citenamefont {Antuono},\ and\
  \citenamefont {Marrone}}]{colagrossi2012particle}%
  \BibitemOpen
  \bibfield  {author} {\bibinfo {author} {\bibfnamefont {A.}~\bibnamefont
  {Colagrossi}}, \bibinfo {author} {\bibfnamefont {B.}~\bibnamefont
  {Bouscasse}}, \bibinfo {author} {\bibfnamefont {M.}~\bibnamefont {Antuono}},
  \ and\ \bibinfo {author} {\bibfnamefont {S.}~\bibnamefont {Marrone}},\
  }\bibfield  {title} {\enquote {\bibinfo {title} {Particle packing algorithm
  for {SPH} schemes},}\ }\href {\doibase 10.1016/j.cpc.2012.02.032} {\bibfield
  {journal} {\bibinfo  {journal} {Computer Physics Communications}\ }\textbf
  {\bibinfo {volume} {183}},\ \bibinfo {pages} {1641--1653} (\bibinfo {year}
  {2012})}\BibitemShut {NoStop}%
\bibitem [{\citenamefont {Sun}, \citenamefont {Colagrossi},\ and\ \citenamefont
  {Zhang}(2018)}]{sun2018numerical}%
  \BibitemOpen
  \bibfield  {author} {\bibinfo {author} {\bibfnamefont {P.-N.}\ \bibnamefont
  {Sun}}, \bibinfo {author} {\bibfnamefont {A.}~\bibnamefont {Colagrossi}}, \
  and\ \bibinfo {author} {\bibfnamefont {A.-M.}\ \bibnamefont {Zhang}},\
  }\bibfield  {title} {\enquote {\bibinfo {title} {Numerical simulation of the
  self-propulsive motion of a fishlike swimming foil using the $\delta$+-{SPH}
  model},}\ }\href {\doibase 10.1016/j.taml.2018.02.007} {\bibfield  {journal}
  {\bibinfo  {journal} {Theoretical and Applied Mechanics Letters}\ }\textbf
  {\bibinfo {volume} {8}},\ \bibinfo {pages} {115--125} (\bibinfo {year}
  {2018})}\BibitemShut {NoStop}%
\bibitem [{\citenamefont {Zhang}, \citenamefont {Hu},\ and\ \citenamefont
  {Adams}(2017)}]{zhang_hu_adams17}%
  \BibitemOpen
  \bibfield  {author} {\bibinfo {author} {\bibfnamefont {C.}~\bibnamefont
  {Zhang}}, \bibinfo {author} {\bibfnamefont {X.~Y.~T.}\ \bibnamefont {Hu}}, \
  and\ \bibinfo {author} {\bibfnamefont {N.~A.}\ \bibnamefont {Adams}},\
  }\bibfield  {title} {\enquote {\bibinfo {title} {A generalized
  transport-velocity formulation for smoothed particle hydrodynamics},}\ }\href
  {\doibase 10.1016/j.jcp.2017.02.016} {\bibfield  {journal} {\bibinfo
  {journal} {Journal of Computational Physics}\ }\textbf {\bibinfo {volume}
  {337}},\ \bibinfo {pages} {216--232} (\bibinfo {year} {2017})}\BibitemShut
  {NoStop}%
\bibitem [{\citenamefont {Shao}\ and\ \citenamefont
  {Lo}(2003)}]{isph:shao:lo:awr:2003}%
  \BibitemOpen
  \bibfield  {author} {\bibinfo {author} {\bibfnamefont {S.}~\bibnamefont
  {Shao}}\ and\ \bibinfo {author} {\bibfnamefont {E.~Y.}\ \bibnamefont {Lo}},\
  }\bibfield  {title} {\enquote {\bibinfo {title} {Incompressible {SPH} method
  for simulating newtonian and non-newtonian flows with a free surface},}\
  }\href {\doibase 10.1016/S0309-1708(03)00030-7} {\bibfield  {journal}
  {\bibinfo  {journal} {Advances in Water Resources}\ }\textbf {\bibinfo
  {volume} {26}},\ \bibinfo {pages} {787 -- 800} (\bibinfo {year}
  {2003})}\BibitemShut {NoStop}%
\bibitem [{\citenamefont {Sun}\ \emph {et~al.}(2017)\citenamefont {Sun},
  \citenamefont {Colagrossi}, \citenamefont {Marrone},\ and\ \citenamefont
  {Zhang}}]{sun2017deltaplus}%
  \BibitemOpen
  \bibfield  {author} {\bibinfo {author} {\bibfnamefont {P.}~\bibnamefont
  {Sun}}, \bibinfo {author} {\bibfnamefont {A.}~\bibnamefont {Colagrossi}},
  \bibinfo {author} {\bibfnamefont {S.}~\bibnamefont {Marrone}}, \ and\
  \bibinfo {author} {\bibfnamefont {A.}~\bibnamefont {Zhang}},\ }\bibfield
  {title} {\enquote {\bibinfo {title} {The $\delta$plus-{SPH} model: Simple
  procedures for a further improvement of the {SPH} scheme},}\ }\href {\doibase
  10.1016/j.cma.2016.10.028} {\bibfield  {journal} {\bibinfo  {journal}
  {Computer Methods in Applied Mechanics and Engineering}\ }\textbf {\bibinfo
  {volume} {315}},\ \bibinfo {pages} {25--49} (\bibinfo {year}
  {2017})}\BibitemShut {NoStop}%
\bibitem [{\citenamefont {Oger}\ \emph {et~al.}(2007)\citenamefont {Oger},
  \citenamefont {Doring}, \citenamefont {Alessandrini},\ and\ \citenamefont
  {Ferrant}}]{oger_improved_2007}%
  \BibitemOpen
  \bibfield  {author} {\bibinfo {author} {\bibfnamefont {G.}~\bibnamefont
  {Oger}}, \bibinfo {author} {\bibfnamefont {M.}~\bibnamefont {Doring}},
  \bibinfo {author} {\bibfnamefont {B.}~\bibnamefont {Alessandrini}}, \ and\
  \bibinfo {author} {\bibfnamefont {P.}~\bibnamefont {Ferrant}},\ }\bibfield
  {title} {\enquote {\bibinfo {title} {An improved {SPH} method: {Towards}
  higher order convergence},}\ }\href {\doibase 10.1016/j.jcp.2007.01.039}
  {\bibfield  {journal} {\bibinfo  {journal} {Journal of Computational
  Physics}\ }\textbf {\bibinfo {volume} {225}},\ \bibinfo {pages} {1472--1492}
  (\bibinfo {year} {2007})}\BibitemShut {NoStop}%
\bibitem [{\citenamefont {Chorin}(1967)}]{chorin_numerical_1967}%
  \BibitemOpen
  \bibfield  {author} {\bibinfo {author} {\bibfnamefont {A.~J.}\ \bibnamefont
  {Chorin}},\ }\bibfield  {title} {\enquote {\bibinfo {title} {A numerical
  method for solving incompressible viscous flow problems},}\ }\href {\doibase
  10.1016/0021-9991(67)90037-X} {\bibfield  {journal} {\bibinfo  {journal}
  {Journal of Computational Physics}\ }\textbf {\bibinfo {volume} {2}},\
  \bibinfo {pages} {12--26} (\bibinfo {year} {1967})}\BibitemShut {NoStop}%
\bibitem [{\citenamefont {Fatehi}\ \emph {et~al.}(2019)\citenamefont {Fatehi},
  \citenamefont {Rahmat}, \citenamefont {Tofighi}, \citenamefont {Yildiz},\
  and\ \citenamefont {Shadloo}}]{fatehi-2019}%
  \BibitemOpen
  \bibfield  {author} {\bibinfo {author} {\bibfnamefont {R.}~\bibnamefont
  {Fatehi}}, \bibinfo {author} {\bibfnamefont {A.}~\bibnamefont {Rahmat}},
  \bibinfo {author} {\bibfnamefont {N.}~\bibnamefont {Tofighi}}, \bibinfo
  {author} {\bibfnamefont {M.}~\bibnamefont {Yildiz}}, \ and\ \bibinfo {author}
  {\bibfnamefont {M.~S.}\ \bibnamefont {Shadloo}},\ }\bibfield  {title}
  {\enquote {\bibinfo {title} {Density-based smoothed particle hydrodynamics
  methods for incompressible flows},}\ }\href {\doibase
  10.1016/j.compfluid.2019.02.018} {\bibfield  {journal} {\bibinfo  {journal}
  {Computers \& Fluids}\ }\textbf {\bibinfo {volume} {185}},\ \bibinfo {pages}
  {22--33} (\bibinfo {year} {2019})}\BibitemShut {NoStop}%
\bibitem [{\citenamefont {Ramachandran}, \citenamefont {Muta},\ and\
  \citenamefont {Ramakrishna}(2021)}]{ramachandran_dual-time_2021}%
  \BibitemOpen
  \bibfield  {author} {\bibinfo {author} {\bibfnamefont {P.}~\bibnamefont
  {Ramachandran}}, \bibinfo {author} {\bibfnamefont {A.}~\bibnamefont {Muta}},
  \ and\ \bibinfo {author} {\bibfnamefont {M.}~\bibnamefont {Ramakrishna}},\
  }\bibfield  {title} {\enquote {\bibinfo {title} {Dual-time smoothed particle
  hydrodynamics for incompressible fluid simulation},}\ }\href {\doibase
  10.1016/j.compfluid.2021.105031} {\bibfield  {journal} {\bibinfo  {journal}
  {Computers \& Fluids}\ }\textbf {\bibinfo {volume} {227}},\ \bibinfo {pages}
  {105031} (\bibinfo {year} {2021})}\BibitemShut {NoStop}%
\bibitem [{\citenamefont {Ramachandran}\ \emph {et~al.}(2021)\citenamefont
  {Ramachandran}, \citenamefont {Bhosale}, \citenamefont {Puri}, \citenamefont
  {Negi}, \citenamefont {Muta}, \citenamefont {Dinesh}, \citenamefont {Menon},
  \citenamefont {Govind}, \citenamefont {Sanka}, \citenamefont {Sebastian},
  \citenamefont {Sen}, \citenamefont {Kaushik}, \citenamefont {Kumar},
  \citenamefont {Kurapati}, \citenamefont {Patil}, \citenamefont {Tavker},
  \citenamefont {Pandey}, \citenamefont {Kaushik}, \citenamefont {Dutt},\ and\
  \citenamefont {Agarwal}}]{pysph2020}%
  \BibitemOpen
  \bibfield  {author} {\bibinfo {author} {\bibfnamefont {P.}~\bibnamefont
  {Ramachandran}}, \bibinfo {author} {\bibfnamefont {A.}~\bibnamefont
  {Bhosale}}, \bibinfo {author} {\bibfnamefont {K.}~\bibnamefont {Puri}},
  \bibinfo {author} {\bibfnamefont {P.}~\bibnamefont {Negi}}, \bibinfo {author}
  {\bibfnamefont {A.}~\bibnamefont {Muta}}, \bibinfo {author} {\bibfnamefont
  {A.}~\bibnamefont {Dinesh}}, \bibinfo {author} {\bibfnamefont
  {D.}~\bibnamefont {Menon}}, \bibinfo {author} {\bibfnamefont
  {R.}~\bibnamefont {Govind}}, \bibinfo {author} {\bibfnamefont
  {S.}~\bibnamefont {Sanka}}, \bibinfo {author} {\bibfnamefont {A.~S.}\
  \bibnamefont {Sebastian}}, \bibinfo {author} {\bibfnamefont {A.}~\bibnamefont
  {Sen}}, \bibinfo {author} {\bibfnamefont {R.}~\bibnamefont {Kaushik}},
  \bibinfo {author} {\bibfnamefont {A.}~\bibnamefont {Kumar}}, \bibinfo
  {author} {\bibfnamefont {V.}~\bibnamefont {Kurapati}}, \bibinfo {author}
  {\bibfnamefont {M.}~\bibnamefont {Patil}}, \bibinfo {author} {\bibfnamefont
  {D.}~\bibnamefont {Tavker}}, \bibinfo {author} {\bibfnamefont
  {P.}~\bibnamefont {Pandey}}, \bibinfo {author} {\bibfnamefont
  {C.}~\bibnamefont {Kaushik}}, \bibinfo {author} {\bibfnamefont
  {A.}~\bibnamefont {Dutt}}, \ and\ \bibinfo {author} {\bibfnamefont
  {A.}~\bibnamefont {Agarwal}},\ }\bibfield  {title} {\enquote {\bibinfo
  {title} {{{PySPH}}: {{A Python-based Framework}} for {{Smoothed Particle
  Hydrodynamics}}},}\ }\href {\doibase 10.1145/3460773} {\bibfield  {journal}
  {\bibinfo  {journal} {ACM Transactions on Mathematical Software}\ }\textbf
  {\bibinfo {volume} {47}},\ \bibinfo {pages} {34:1--34:38} (\bibinfo {year}
  {2021})}\BibitemShut {NoStop}%
\bibitem [{\citenamefont {Maci{\'a}}\ \emph {et~al.}(2011)\citenamefont
  {Maci{\'a}}, \citenamefont {Antuono}, \citenamefont {Gonz{\'a}lez},\ and\
  \citenamefont {Colagrossi}}]{maciaTheoreticalAnalysisNoSlip2011}%
  \BibitemOpen
  \bibfield  {author} {\bibinfo {author} {\bibfnamefont {F.}~\bibnamefont
  {Maci{\'a}}}, \bibinfo {author} {\bibfnamefont {M.}~\bibnamefont {Antuono}},
  \bibinfo {author} {\bibfnamefont {L.~M.}\ \bibnamefont {Gonz{\'a}lez}}, \
  and\ \bibinfo {author} {\bibfnamefont {A.}~\bibnamefont {Colagrossi}},\
  }\bibfield  {title} {\enquote {\bibinfo {title} {Theoretical {{Analysis}} of
  the {{No}}-{{Slip Boundary Condition Enforcement}} in {{SPH Methods}}},}\
  }\href {\doibase 10.1143/PTP.125.1091} {\bibfield  {journal} {\bibinfo
  {journal} {Progress of Theoretical Physics}\ }\textbf {\bibinfo {volume}
  {125}},\ \bibinfo {pages} {1091--1121} (\bibinfo {year} {2011})}\BibitemShut
  {NoStop}%
\bibitem [{\citenamefont {Adami}, \citenamefont {Hu},\ and\ \citenamefont
  {Adams}(2012)}]{Adami2012}%
  \BibitemOpen
  \bibfield  {author} {\bibinfo {author} {\bibfnamefont {S.}~\bibnamefont
  {Adami}}, \bibinfo {author} {\bibfnamefont {X.}~\bibnamefont {Hu}}, \ and\
  \bibinfo {author} {\bibfnamefont {N.}~\bibnamefont {Adams}},\ }\bibfield
  {title} {\enquote {\bibinfo {title} {{A generalized wall boundary condition
  for smoothed particle hydrodynamics}},}\ }\href {\doibase
  10.1016/j.jcp.2012.05.005} {\bibfield  {journal} {\bibinfo  {journal}
  {Journal of Computational Physics}\ }\textbf {\bibinfo {volume} {231}},\
  \bibinfo {pages} {7057--7075} (\bibinfo {year} {2012})}\BibitemShut {NoStop}%
\bibitem [{\citenamefont {Muta}\ and\ \citenamefont
  {Ramachandran}(2021)}]{mutaadaptive2021}%
  \BibitemOpen
  \bibfield  {author} {\bibinfo {author} {\bibfnamefont {A.}~\bibnamefont
  {Muta}}\ and\ \bibinfo {author} {\bibfnamefont {P.}~\bibnamefont
  {Ramachandran}},\ }\bibfield  {title} {\enquote {\bibinfo {title} {Efficient
  and {{Accurate Adaptive Resolution}} for {{Weakly-Compressible SPH}}},}\
  }\href@noop {} {\bibfield  {journal} {\bibinfo  {journal} {arXiv:2107.01276
  [physics]}\ } (\bibinfo {year} {2021})},\ \Eprint
  {http://arxiv.org/abs/2107.01276} {arXiv:2107.01276 [physics]} \BibitemShut
  {NoStop}%
\bibitem [{\citenamefont
  {Monaghan}(2000)}]{sph:tensile-instab:monaghan:jcp2000}%
  \BibitemOpen
  \bibfield  {author} {\bibinfo {author} {\bibfnamefont {J.~J.}\ \bibnamefont
  {Monaghan}},\ }\bibfield  {title} {\enquote {\bibinfo {title} {{SPH} without
  a tensile instability},}\ }\href {\doibase 10.1006/jcph.2000.6439} {\bibfield
   {journal} {\bibinfo  {journal} {Journal of Computational Physics}\ }\textbf
  {\bibinfo {volume} {159}},\ \bibinfo {pages} {290--311} (\bibinfo {year}
  {2000})}\BibitemShut {NoStop}%
\bibitem [{\citenamefont {Nguyen}\ \emph {et~al.}(2008)\citenamefont {Nguyen},
  \citenamefont {Rabczuk}, \citenamefont {Bordas},\ and\ \citenamefont
  {Duflot}}]{nguyen_meshless_2008}%
  \BibitemOpen
  \bibfield  {author} {\bibinfo {author} {\bibfnamefont {V.~P.}\ \bibnamefont
  {Nguyen}}, \bibinfo {author} {\bibfnamefont {T.}~\bibnamefont {Rabczuk}},
  \bibinfo {author} {\bibfnamefont {S.}~\bibnamefont {Bordas}}, \ and\ \bibinfo
  {author} {\bibfnamefont {M.}~\bibnamefont {Duflot}},\ }\bibfield  {title}
  {\enquote {\bibinfo {title} {Meshless methods: {A} review and computer
  implementation aspects},}\ }\href {\doibase 10.1016/j.matcom.2008.01.003}
  {\bibfield  {journal} {\bibinfo  {journal} {Mathematics and Computers in
  Simulation}\ }\textbf {\bibinfo {volume} {79}},\ \bibinfo {pages} {763--813}
  (\bibinfo {year} {2008})}\BibitemShut {NoStop}%
\bibitem [{\citenamefont {Muta}, \citenamefont {Ramachandran},\ and\
  \citenamefont {Negi}(2020)}]{muta2020efficient}%
  \BibitemOpen
  \bibfield  {author} {\bibinfo {author} {\bibfnamefont {A.}~\bibnamefont
  {Muta}}, \bibinfo {author} {\bibfnamefont {P.}~\bibnamefont {Ramachandran}},
  \ and\ \bibinfo {author} {\bibfnamefont {P.}~\bibnamefont {Negi}},\
  }\bibfield  {title} {\enquote {\bibinfo {title} {An efficient, open source,
  iterative {ISPH} scheme},}\ }\href@noop {} {\bibfield  {journal} {\bibinfo
  {journal} {Computer Physics Communications}\ ,\ \bibinfo {pages} {107283}}
  (\bibinfo {year} {2020})}\BibitemShut {NoStop}%
\end{thebibliography}%

\end{document}